\documentclass[12pt,a4paper]{article}

\usepackage{amsbsy,amsfonts,amsmath,amssymb,amsthm}
\usepackage{array}
\usepackage{bm}
\usepackage{color}
\usepackage{enumerate}
\usepackage{mathrsfs}
\usepackage{latexsym}
\usepackage[colorlinks=true]{hyperref}
\hypersetup{urlcolor=blue, citecolor=blue, linkcolor=blue}

\definecolor{wineRed}{rgb}{0.7,0,0.3}
\definecolor{grandBleu}{rgb}{0,0,0.8}
\definecolor{darkGreen}{rgb}{0,0.4,0}
\definecolor{blueViolet}{rgb}{0.4,0,1.0}
\definecolor{bloodOrange}{rgb}{0.85,0.05,0}
\definecolor{mycolor}{rgb}{0.8,0,0.2}
\definecolor{}{rgb}{0.8,0,0.2}

\usepackage[textsize=small]{todonotes}
\usepackage{cite}

\usepackage{comment}

\DeclareMathAlphabet{\mathpzc}{OT1}{pzc}{m}{it}

\usepackage[textsize=small]{todonotes}
\numberwithin{equation}{section}

\theoremstyle{plain}
\newtheorem{mTh}{Main Theorem} 
\newtheorem{lem}{Lemma}
\newtheorem{keyLem}{Key-Lemma}
\newtheorem{prop}{Proposition}
\newtheorem{cor}{Corollary}

\theoremstyle{definition}
\newtheorem{defn}{Definition}
\newtheorem{rem}{Remark}

\newtheorem{ex}{Example}

\def\N{\mathbb{N}}
\def\R{\mathbb{R}}

\def\sH{\mathscr{H}}

\def\ds{\displaystyle}
\def\ts{\textstyle}
\def\Sgn{\mathop{\mathrm{Sgn}}\nolimits}

\newcommand{\jump}[1]{\ensuremath{[\hspace{-0.34ex}[#1]\hspace{-0.34ex}]} }

\begin{document}
\thispagestyle{plain}
\begin{center}
    \textbf{\Large Temperature Control of PDE Constrained Optimization Problems Governed by {K}obayashi--{W}arren--{C}arter Type \\[0.75ex] Models of Grain Boundary Motions}\footnotemark[1]
\end{center}
    \bigskip
\begin{center}
    \textsc{Harbir Antil} 
    \\[1ex]
    {Department of Mathematical Sciences and the Center for Mathematics and Artificial Intelligence,  (CMAI), George Mason University, Fairfax, VA 22030, USA}
    \\[0ex]
    ({\ttfamily hantil@gmu.edu})
\end{center}
\begin{center}
    \textsc{Shodai Kubota$^*$}
    \\[1ex]
    {Department of Mathematics and Informatics, \\ Graduate School of Science and Engineering, Chiba University, \\ 1-33, Yayoi-cho, Inage-ku, 263-8522, Chiba, Japan}
    \\[0ex]
    ({\ttfamily skubota@chiba-u.jp})
\end{center}
\begin{center}
    \textsc{Ken Shirakawa}
    \\[1ex]
    {Department of Mathematics, Faculty of Education, Chiba University, \\ 1-33, Yayoi-cho, Inage-ku, 263-8522, Chiba, Japan}
    \\[0ex]
    ({\ttfamily sirakawa@faculty.chiba-u.jp})
\end{center}
\begin{center}
    \textsc{and}
\end{center}
\begin{center}
\textsc{Noriaki Yamazaki}
\\[1ex]
{Department of Mathematics, Faculty of Engineering, Kanagawa University, \\ 3-27-1, Rokkakubashi, Kanagawa-ku, Yokohama, 221-8686, Japan}
\\[0ex]
({\ttfamily noriaki@kanagawa-u.ac.jp})
\end{center}

\footnotetext[1]{
The work of the third author supported by Grant-in-Aid for Scientific Research (C) No. 16K05224 and No. 20K03672, JSPS. The work of the forth author supported by Grant-in-Aid for Scientific Research (C) No. 20K03665 , JSPS. In addition, the work of the first and the third authors is partially supported by the Air Force Office of Scientific Research (AFOSR) under Award NO: FA9550-19-1-0036 and NSF grants DMS-1818772 and DMS-1913004.  
\\[1ex]
\\[1ex]
}
\bigskip

\noindent
{\bf Abstract.} In this paper, we consider a class of optimal control problems governed by state-equations of Kobayashi--Warren--Carter type. The control is given by physical temperature. The focus is on problems in dimensions less than equal to 4.  The results are divided in four Main Theorems, concerned with: solvability and parameter-dependence of state-equations and optimal control problems; the first order necessary optimality conditions for these regularized optimal control problems. Subsequently, we derive the limiting systems and optimality conditions and study their well-posedness.
\newpage

\section*{Introduction}
Let $ (0, T) $ be a time-interval with a constant $ 0 < T < \infty $, and let $N \in \{ 2, 3, 4\}$ denotes the spatial dimension. 
Let $\Omega \subset \R^N $ be a bounded domain with a Lipschitz boundary $ \Gamma := \partial \Omega $, and let $ n_\Gamma $ be the unit outer normal on $ \Gamma $. Besides, we set $ Q := (0, T) \times \Omega $ and $ \Sigma := (0, T) \times \Gamma $, and we define $ H := L^2(\Omega) $ with norm $| \cdot |_H$, $ V := H^1(\Omega) $, $ V_0 := H_0^1(\Omega) $, and  $ \mathscr{H} := L^2(0, T; L^2(\Omega)) $,  as the base spaces for this work. Moreover, we set:
$$
    \begin{array}{c}
        \jump{\kappa^0, \kappa^1} :=  \left\{ \begin{array}{l|l}
                    \tilde{u} \in \mathscr{H} & \kappa^0 \leq \tilde{u} \leq \kappa^1 \mbox{ a.e. in $ Q $}
                \end{array} \right\},
                \\[1ex]
        \mbox{for arbitrary measurable obstacles $ \kappa^\ell : Q \longrightarrow [-\infty, \infty] $, \quad $ \ell = 0, 1 $,}
    \end{array}
$$
and define a family of functional classes $ \mathfrak{K} \subset 2^{\mathscr{H}} $, as follows:
\begin{equation}\label{obs_K}
\mathfrak{K} := \left\{ \begin{array}{l|l} 
        {K} \subset \mathscr{H} & \parbox{8cm}{
            $ {K} = \jump{\kappa^0, \kappa^1} $ for some measurable obstacles $ \kappa^\ell : Q \longrightarrow [-\infty, \infty] $, $ \ell = 0, 1 $, such that $ \kappa^0 \leq \kappa^1 $ a.e. in $ Q $ (i.e. $ K \ne \emptyset $) 
        }
    \end{array}\right\}.
\end{equation}
\medskip

In this paper, we consider a class of optimal control problems, denoted by (\hyperlink{OP}{OP})$_\varepsilon^K $, which are labeled by constants $ \varepsilon \geq 0 $ and functional classes $ K = \jump{\kappa^0, \kappa^1} \in \mathfrak{K} $, with the obstacles $ \kappa^\ell : Q \longrightarrow [-\infty, \infty] $, $ \ell = 0, 1 $. For every $ \varepsilon \geq 0 $ and $ K = \jump{\kappa^0, \kappa^1} \in \mathfrak{K} $, the optimal control problem (\hyperlink{OP}{OP})$_{\varepsilon}^K$ is prescribed as follows:
\begin{description}
    \item[\textmd{(\hypertarget{OP}{OP})$_{\varepsilon}^K$}]Find a pair of functions $ [u^*, v^*] \in [\mathscr{H}]^2 $, called the \emph{optimal control}, such that
        \begin{align*}
            [u^* &, v^*] \in \mathscr{U}_\mathrm{ad}^K := \left\{ \begin{array}{l|l}
                [\tilde{u}, \tilde{v}] \in [\mathscr{H}]^2 & \tilde{u} \in K
        \end{array} \right\},
            \\[0.5ex]
            \ds \mbox{and}~~& \mathcal{J}_{\varepsilon}(u^*, v^*) = 
            \min \left\{ \begin{array}{l|l} 
                \mathcal{J}_{\varepsilon}(u, v) & [u, v] \in \mathscr{U}_\mathrm{ad}^K
        \end{array} \right\},
    \end{align*}
        where $ \mathcal{J}_{\varepsilon} = \mathcal{J}_{\varepsilon}(u, v) $ is a cost functional on $ [\mathscr{H}]^2 $, defined as follows:
        \begin{align}\label{J}
            \mathcal{J}_{\varepsilon} : [u, &\,  v] \in [\mathscr{H}]^2 \mapsto \mathcal{J}_{\varepsilon}(u, v) 
            \nonumber
            \\
            := &~ \frac{M_\eta }{2} \int_0^T |(\eta -\eta_\mathrm{ad})(t)|_H^2 \, dt  +\frac{M_\theta }{2} \int_0^T |(\theta -\theta_\mathrm{ad})(t)|_H^2 \, dt
            \\
            &~ +\frac{M_u }{2} \int_0^T |u(t)|_H^2 \, dt +\frac{M_v }{2} \int_0^T |v(t)|_H^2 \, dt \quad \in [0, \infty),  
            \nonumber
        \end{align}
        with $ [\eta, \theta] \in [\mathscr{H}]^2 $ solving the state-system, denoted by (S)$_{\varepsilon}$: 
\end{description}

~~~~(S)$ _{\varepsilon} $
\vspace{-0ex}
\begin{equation}\label{1}
\left\{ \parbox{11cm}{
    $ \partial_{t} \eta -\Delta \eta +g(\eta) +\alpha'(\eta) \sqrt{\varepsilon^2 +|\nabla \theta|^2} = M_u u $ \quad in $ Q $,
\\[1ex]
$ \nabla \eta(t, x) \cdot n_\Gamma = 0 $, \quad $ (t, x) \in \Sigma $, 
\\[1ex]
$ \eta(0, x) = \eta_{0}(x) $, \quad $ x \in \Omega $;
}\right. 
\end{equation}
\begin{equation}\label{2}
\left\{\parbox{11cm}{
    $ \displaystyle \alpha_{0}(t, x) \partial_{t} \theta -\mbox{div} \left( \alpha(\eta) \frac{\nabla\theta}{\sqrt{\varepsilon^2 +|\nabla \theta|^2}} +\nu^{2}\nabla\theta\right) = M_v v $ \quad in $Q$,
\\[1ex]
$ \theta(t, x) = 0 $, \quad $ (t, x) \in \Sigma $, 
\\[1ex]
$ \theta(0, x) = \theta_{0}(x)$, \quad $x \in \Omega $.
}\right. 
\end{equation}

The state-system (S)$_{\varepsilon}$ is based on a phase field model of grain boundary motion, known as \emph{Kobayashi--Warren--Carter system} (cf. \cite{MR1752970, MR1794359}). In this context, the unknowns $ \eta \in \mathscr{H} $ and $ \theta \in \mathscr{H} $ are order parameters that indicate the \emph{orientation order} and \emph{orientation angle} of the polycrystal body, respectively. Besides, $ [\eta_0, \theta_0] \in V \times V_0 $ is an \emph{initial pair}, i.e. a pair of initial data of $ [\eta, \theta] $. The \emph{forcing pair} $ [u, v] \in [\mathscr{H}]^2 $ denotes the control variables that can control the profile of solution $ [\eta, \theta] \in [\mathscr{H}]^2 $ to (S)$_{\varepsilon}$. Additionally, $ 0 < \alpha_0 \in W^{1, \infty}(Q) $ and $ 0 < \alpha \in C^2(\R) $ are given functions to reproduce the mobilities of grain boundary motions. Finally, $ g \in W_\mathrm{loc}^{1, \infty}(\R) $ is a perturbation for the orientation order $ \eta $, and $ \nu > 0 $ is a fixed constant to relax the diffusion of the orientation angle $ \theta $. 
\medskip

The first part \eqref{1} of the state-system (S)$_{\varepsilon}$ is the initial-boundary value problem of an Allen--Cahn type equation, so that the forcing term $ u $ can be regarded as a \emph{temperature control} of the grain boundary formation. Also, the second problem \eqref{2} is the initial-boundary value problem to reproduce crystalline micro-structure of polycrystal, and the case of $ \varepsilon = 0 $ is the closest to the original setting adopted by Kobayashi et al \cite{MR1752970, MR1794359}. Indeed, when $ \varepsilon = 0 $, the quasi-linear diffusion as in \eqref{2} is described in a singular form $ -\mbox{div} \bigl( \alpha(\eta) \frac{\nabla \theta}{|\nabla \theta|} +\nu^2 \nabla \theta \bigr) $, and it is known that this type of singularity is effective to reproduce the \emph{facet}, i.e. the locally uniform (constant) phase in each oriented grain (cf. \cite{MR1752970, MR1794359,MR2746654,MR2101878,MR2436794,MR2033382,MR2223383,MR3038131,MR3670006,MR3268865,MR2836557,MR1865089,MR1712447,MR3951294}). Hence, the systems (S)$_{\varepsilon}$, for positive $ \varepsilon $, can be regarded as \emph{regularized approximating systems}, that are to approach to the physically realistic situation (S)$_{0}$, in the limit $ \varepsilon \downarrow 0$.
\medskip

Meanwhile, in the optimal control problem (\hyperlink{OP}{OP})$_{\varepsilon}^K$, the class $ K = \jump{\kappa^0, \kappa^1} \in \mathfrak{K} $ is to constrain the range of temperature control $ u $, and the obstacles $ \kappa^\ell : Q \longrightarrow [-\infty, \infty] $, $ \ell = 0, 1 $, indicate the \emph{control bounds} of the temperature. The pair of functions $  [\eta_\mathrm{ad}, \theta_\mathrm{ad}] \in [\mathscr{H}]^2 $ is a given \emph{admissible target profile} of $ [\eta, \theta] \in [\mathscr{H}]^2 $. Moreover, $ M_\eta \geq 0 $, $ M_\theta \geq 0 $, $ M_u \geq 0 $, and $ M_v \geq 0 $ are fixed constants. 
\medskip

The objective of this paper is to significantly extend the results of our previous work \cite{MR4218112}, which dealt with: 
\begin{description}
    \item[\textmd{~~\hypertarget{sharp1}{$\sharp \, 1$})}]key-properties of the state-systems (S)$ _{\varepsilon} $ with 1-dimensional domain $\Omega \subset (0,1)$; 
        \vspace{-1ex}
    \item[\textmd{~~\hypertarget{sharp2}{$\sharp \, 2$})}]mathematical analysis of the optimal control problem (\hyperlink{OP}{OP})$ _{\varepsilon}^K $, for $ \varepsilon \geq 0 $, but with 1-dimensional domain $\Omega \subset (0,1)$ without any control constraints, i.e., $ K = \jump{\kappa^0, \kappa^1} \in \mathfrak{K} $ ($=\mathscr{H} $);
\end{description}
In light of this, the novelty of this work is in:
\begin{description}
    \item[\textmd{~~\hypertarget{sharp3}{$\sharp \, 3$})}]the development of a mathematical analysis to obtain optimal controls of grain boundaries under
        the higher dimensional setting $ N \in \{2, 3, 4\} $ of the spatial domain, 
        and the temperature constraint $ K = \jump{\kappa^0, \kappa^1} \in \mathfrak{K} $.
\end{description}
In addition, the presence of constraints $ K = \jump{\kappa^0, \kappa^1} \in \mathfrak{K} $ makes the mathematical analysis further challenging. Notice that such constraints are meaningful from a practical point of view.
We further emphasize that in the main part of this work, the $ L^\infty $-boundedness of $ \eta $ will be essential, and the main results will be valid under the following assumption on the data:
\begin{description}
    \item[\textmd{(\hypertarget{rs0}{r.s.0})}]
        $ \varepsilon > 0 $, $  [\eta_0, \theta_0] \in D_0 := \bigl( V \cap L^\infty(\Omega) \bigr) \times V_0 $, and $ K \in \mathfrak{K}_0 $, where
        \begin{equation}\label{obs_K0} 
            \mathfrak{K}_0 :=  \mathfrak{K} \cap 2^{L^\infty(Q)} = \left\{ \begin{array}{l|l}
                K & \parbox{5cm}{
                    $ K = \jump{\kappa^0, \kappa^1} \in \mathfrak{K} $ such that $ \kappa^\ell \in L^\infty(Q) $, $ \ell = 0, 1 $
                }
            \end{array} \right\}.
        \end{equation} 
\end{description}
Hence, in general cases of constraints $ K \in \mathfrak{K} $ (including no constraint case), we will be forced to adopt some limiting (approximating) approach on the basis of the results under the restricted situation (\hyperlink{rs0}{r.s.0}). 
\medskip

Now, in view of \hyperlink{sharp1}{$\sharp \, 1$})--\hyperlink{sharp3}{$\sharp \, 3$}), we set the goal of this paper to prove four Main Theorems, summarized as follows.
\begin{description}
    \item[{\boldmath Main Theorem \ref{mainTh01}:}]Mathematical results concerning the following items.
    \vspace{-0.5ex}
\item[{\boldmath~~(\hyperlink{I-A}{I-A})(Solvability of state-systems):}]Existence and uniqueness for the state-system (S)$_{\varepsilon}$, for every $ \varepsilon \geq 0 $, initial pair $ w_0 = [\eta_0, \theta_0] \in V \times V_0 $, and forcing pair $ [u, v] \in [\mathscr{H}]^2 $.
    \vspace{-0.5ex}
\item[{\boldmath~~(\hyperlink{I-B}{I-B})(Continuous dependence on data among state-systems):}]Continuous dependence of solutions to the systems (S)$_{\varepsilon}$, with respect to 
 the constant $ \varepsilon \geq 0 $, initial pair $ [\eta_0, \theta_0] \in V \times V_0 $, and forcing pair $ [u, v] \in [\mathscr{H}]^2 $. 

\item[{\boldmath Main Theorem \ref{mainTh02}:}]Mathematical results concerning the following items.
    \vspace{-0.5ex}
\item[{\boldmath~~(\hyperlink{II-A}{II-A})(Solvability of optimal control problems):}]Existence for the optimal control problem (\hyperlink{OP}{OP})$_{\varepsilon}^K$, for every constant $ \varepsilon \geq 0 $, initial pair $ [\eta_0, \theta_0] \in V \times V_0 $, and constraint $ K = \jump{\kappa^0, \kappa^1} \in \mathfrak{K} $.
    \vspace{-0.5ex}
\item[{\boldmath~~(\hyperlink{II-B}{II-B})(Parameter dependence of optimal controls):}]Some semi-continuous dependence of  the optimal controls, with respect to the constant $ \varepsilon \geq 0 $, initial pair $ [\eta_0, \theta_0] \in V \times V_0 $, and constraint $ K = \jump{\kappa^0, \kappa^1} \in \mathfrak{K} $.

\item[{\boldmath Main Theorem \ref{mainTh03}:}]mathematical results concerning the following items.
    \vspace{-0.5ex}
\item[{\boldmath~~(\hyperlink{III-A}{III-A})(Necessary optimality conditions under (\hyperlink{rs0}{r.s.0})):}] Derivation of first order necessary optimality conditions for (\hyperlink{OP}{OP})$_{\varepsilon}^K$, via adjoint method, under the restricted situation (\hyperlink{rs0}{r.s.0}). 
    \vspace{-0.5ex}
\item[{\boldmath~~(\hyperlink{III-B}{III-B})(Specific parameter dependence under (\hyperlink{rs0}{r.s.0})):}]Strong parameter dependence of optimal controls, which is specifically obtained under (\hyperlink{rs0}{r.s.0}).

\item[{\boldmath Main Theorem \ref{mainTh04} (Limiting optimality conditions):}]The optimality conditions which are obtained as approximation limits of the necessary conditions under (\hyperlink{rs0}{r.s.0}). 
\end{description}

This paper is organized as follows. The Main Theorems are stated in Section 3, after the preliminaries in Section 1, and the auxiliary lemmas in Section 2. The part after Section 3 will be divided in Sections 4--7, and these four Sections will be devoted to the proofs of the respective four Main Theorems 1--4.

\section{Preliminaries} 

We begin by prescribing the notations used throughout this paper. 
\medskip

\noindent
\underline{\textbf{\textit{Basic notations.}}} 
For arbitrary $ r_0 $, $ s_0 \in [-\infty, \infty]$, we define:
\begin{equation*}
r_0 \vee s_0 := \max\{r_0, s_0 \}\ \mbox{and}\ r_0 \wedge s_0 := \min\{r_0, s_0 \},
\end{equation*}
and in particular, we set:
\begin{equation*}
    [r]^+ := r \vee 0 \ \mbox{and}\ [r]^- :=  -(r \wedge 0), \mbox{ for any $ r \in \R $.}
\end{equation*}

For any dimension $ d \in \N $, we denote by $\mathcal{L}^d$ the $d$-dimensional Lebesgue measure. The measure theoretical phrases, such as ``a.e.'', ``$dt$'', ``$dx$'', and so on, are all with respect to the Lebesgue measure in each corresponding dimension.
\\

\noindent
\underline{\textbf{\textit{Abstract notations.}}}
For an abstract Banach space $ X $, we denote by $ |\cdot|_{X} $ the norm of $ X $, and denote by $ \langle \cdot, \cdot \rangle_X $ the duality pairing between $ X $ and its dual $ X^* $. In particular, when $ X $ is a Hilbert space, we denote by $ (\cdot,\cdot)_{X} $ the inner product of $ X $. Moreover, when there is no possibility of confusion, we uniformly denote by $ |\cdot| $ the norm of Euclidean spaces, and for any dimension $ d \in \N $, we write the inner product (scalar product) of $ \R^d $, as follows:
    \begin{align*}
        & y \cdot \tilde{y} = \sum_{i = 1}^d y_i \tilde{y}_i, 
        \mbox{ for all } y = [y_1, \dots, y_d], ~ \tilde{y} = [\tilde{y}_1, \dots, \tilde{y}_d] \in \R^d.
    \end{align*}

For any subset $ A $ of a Banach space $ X $, let $ \chi_A : X \longrightarrow \{0, 1\} $ be the characteristic function of $ A $, i.e.:
    \begin{equation*}
        \chi_A: w \in X \mapsto \chi_A(w) := \begin{cases}
            1, \mbox{ if $ w \in A $,}
            \\[0.5ex]
            0, \mbox{ otherwise.}
        \end{cases}
    \end{equation*}

For two Banach spaces $ X $ and $ Y $,  we denote by $  \mathscr{L}(X; Y)$ the Banach space of bounded linear operators from $ X $ into $ Y $, and in particular, we let $ \mathscr{L}(X) := \mathscr{L}(X; X) $. 

For Banach spaces $ X_1, \dots, X_d $, with $ 1 < d \in \mathbb{N} $, let $ X_1 \times \dots \times X_d $ be the product Banach space endowed with the norm $ |\cdot|_{X_1 \times \cdots \times X_d} := |\cdot|_{X_1} + \cdots +|\cdot|_{X_d} $. However, when all $ X_1, \dots, X_d $ are Hilbert spaces, $ X_1 \times \dots \times X_d $ denotes the product Hilbert space endowed with the inner product $ (\cdot, \cdot)_{X_1 \times \cdots \times X_d} := (\cdot, \cdot)_{X_1} + \cdots +(\cdot, \cdot)_{X_d} $ and the norm $ |\cdot|_{X_1 \times \cdots \times X_d} := \bigl( |\cdot|_{X_1}^2 + \cdots +|\cdot|_{X_d}^2 \bigr)^{\frac{1}{2}} $. In particular, when all $ X_1, \dots,  X_d $ coincide with a Banach space $ Y $, we write:
\begin{equation*}
    [Y]^d := \overbrace{Y \times \cdots \times Y}^{\mbox{$d$ times}}.
\end{equation*}
Additionally, for any transform (operator) $ \mathcal{T} : X \longrightarrow Y $, we let:
\begin{equation*}
    \mathcal{T}[w_1, \dots, w_d] := \bigl[ \mathcal{T} w_1, \dots, \mathcal{T} w_d \bigl] \mbox{ in $ [Y]^d $, \quad for any $ [w_1, \dots, w_d] \in [X]^d $.}
\end{equation*}

\noindent
\underline{\textbf{\textit{Specific notations of this paper.}}} 
As is mentioned in the previous section, let $ (0, T) \subset \R$ be a bounded time-interval with a finite constant $ T > 0 $, and let $N \in \{ 2, 3, 4\}$ be a constant of spatial dimension. Let $\Omega \subset \R^N $ be a fixed spatial bounded domain with a smooth boundary $ \Gamma := \partial\Omega $. We denote by $n_\Gamma$ the unit outward normal vector on $\Gamma$. Besides, we set $ Q := (0, T) \times \Omega $ and $ \Sigma := (0, T) \times \Gamma $. Especially, we denote by $ \partial_t $, $ \nabla $, and $ \mathrm{div} $ the distributional time-derivative, the distributional gradient, and distributional divergence, respectively. 
\medskip

On this basis, we define  
\begin{equation*} \begin{cases}
    H := L^2(\Omega) \mbox{ and } \mathscr{H} := L^2(0, T; H),
\\
    V := H^1(\Omega) \mbox{ and } \mathscr{V} := L^2(0, T; V),
\\
    V_0 := H_0^1(\Omega) \mbox{ and } \mathscr{V}_0 := L^2(0, T; V_0),
\\
    \mathscr{X} := L^\infty(Q) \times \mathscr{H}.
\end{cases}
\end{equation*}
Also, we identify the Hilbert spaces $ H $ and $ \mathscr{H} $ with their dual spaces. Based on the identifications, we have the following relationships of continuous embeddings:
\begin{equation*}
\begin{cases}
    V \subset H = H^* \subset V^* \mbox{ and } \mathscr{V} \subset \mathscr{H} = \mathscr{H}^* \subset \mathscr{V}^*,
\\
    V_0 \subset H = H^* \subset V_0^* \mbox{ and } \mathscr{V}_0 \subset \mathscr{H} = \mathscr{H}^* \subset \mathscr{V}_0^*,
\end{cases}
\end{equation*}
among the Hilbert spaces $ H $, $ V $, $ V_0 $, $ \mathscr{H} $, $ \mathscr{V} $, and $ \mathscr{V}_0 $, and the respective dual spaces $ H^* $, $ V^* $, $ V_0^* $, $ \mathscr{H}^* $, $ \mathscr{V}^* $, and $ \mathscr{V}_0^* $. Additionally, in this paper, we define the topology of the Hilbert space $ V_0 $ by using the following inner product:
\begin{equation*}
    (w, \tilde{w})_{V_0} := (\nabla w, \nabla \tilde{w})_{[H]^N}, \mbox{ for all $ w, \tilde{w} \in V_0 $.}
\end{equation*}

\begin{rem}(cf. \cite[Remark 3]{MR3888633})\label{Rem.Prelim01}
     Due to the restriction $N \in \{2, 3, 4\}$ of spatial dimension, we can suppose the continuous embedding $ V \subset L^4(\Omega) $, and we can easily check that:
    \begin{itemize}
        \item[(i)] if $ 0 \leq \check{\mu} \in H $ and $ \check{p} \in V $, then $ \sqrt{\check{\mu}} \check{p} \in H $, $ \check{\mu} \check{p} \in V^* $, and
            \begin{equation*}
                \begin{cases}
                    ({\sqrt{\check{\mu}}} \check{p},  \psi)_H \leq C_{V}^{L^4} |\check{\mu}|_H^{\frac{1}{2}} |\check{p}|_V |\psi|_H, \mbox{ for any $ \psi \in H $,}
                    \\[1ex]
                    \langle \check{\mu} \check{p}, \check{\psi} \rangle_V = (\sqrt{\check{\mu}} \check{p}, \sqrt{\check{\mu}} \check{\psi})_H \leq (C_V^{L^4})^2 |\check{\mu}|_H |\check{p}|_V |\check{\psi}|_V, \mbox{ for any $ \check{\psi} \in V $;}
                \end{cases}
            \end{equation*}
        \item[(ii)] if $ 0 \leq \hat{\mu} \in L^\infty(0, T; H) $ and $ \hat{p} \in \mathscr{V} $, then $ \sqrt{\hat{\mu}} \hat{p} \in \mathscr{H} $, $ \hat{\mu} \hat{p} \in \mathscr{V}^*$, and
            \begin{equation*}
                \begin{cases}
                    ({\sqrt{\hat{\mu}}} \hat{p},  \varphi)_{\mathscr{H}} \leq C_{V}^{L^4} |\hat{\mu}|_{L^\infty(0, T; H)}^{\frac{1}{2}} |\hat{p}|_{\mathscr{V}} |\varphi|_{\mathscr{H}}, \mbox{ for any $ \varphi \in \mathscr{H} $,}
                    \\[1ex]
                    \langle \hat{\mu} \hat{p}, \hat{\varphi} \rangle_{\mathscr{V}} = (\sqrt{\hat{\mu}} \hat{p}, \sqrt{\hat{\mu}} \hat{\varphi})_{\mathscr{H}} \leq (C_V^{L^4})^2 |\hat{\mu}|_{L^\infty(0, T; H)} |\hat{p}|_{\mathscr{V}} |\hat{\varphi}|_{\mathscr{V}}, \mbox{ for any $ \hat{\varphi} \in \mathscr{V} $;}
                \end{cases}
            \end{equation*}
                \end{itemize}
                where $ C_V^{L^4} > 0  $ is the constant of embedding $ V \subset L^4(\Omega) $. 
\end{rem}

Finally, we define:
\begin{equation*}
    D := V \times V_0, \mbox{ ~and~  } D_0 := \bigl( V \cap L^\infty(\Omega) \bigr) \times V_0,
\end{equation*}
as the notations to specify the range of the initial pair $ [\eta_0, \theta_0] $ in the state-system.
\bigskip

\noindent
\underline{\textbf{\textit{Notations in convex analysis. (cf. \cite[Chapter II]{MR0348562})}}} 
    Let $ X $ be an abstract Hilbert space $ X $. Then, any closed and convex set $ K \subset X $ defines a single-valued operator $ \mathrm{proj}_K : X \longrightarrow K $, which maps any $ w \in X $ to a point $ \mathrm{proj}_K(w) \in K $, satisfying:
\begin{equation*}
    \bigl| \mathrm{proj}_K(w) -w \bigr|_X = \min \left\{ \begin{array}{l|l}
        |\tilde{w} -w|_X & \tilde{w} \in K
    \end{array} \right\}. 
\end{equation*}
The operator $ \mathrm{proj}_K $ is called the \emph{orthogonal projection} (or \emph{projection} in short) onto $ K $.
\begin{rem}[Key-properties of the projection]\label{Rem.proj01}
    Let $ K $ be a closed and convex set in a Hilbert space $ X $. Then, the following facts hold.
    \begin{description}
        \item[(\hypertarget{Fact1l}{Fact\,1})]The projection $ \mathrm{proj}_K : X \longrightarrow K $ is a nonexpansive operator from $ X $ into itself, i.e.:
    \begin{equation*}
        |\mathrm{proj}_K(w^1) -\mathrm{proj}_K(w^2)|_X \leq |w^1 -w^2|_X, \mbox{ for all $ w^\ell \in X $, $ \ell = 1, 2 $.}
    \end{equation*}
\item[(\hypertarget{Fact2l}{Fact\,2})]$ w_K^\circ = \mathrm{proj}_K(w) $ in $ X $, iff. $ (w -w_K^\circ, \tilde{w} -w_K^\circ)_X \leq 0 $,  for any $ \tilde{w} \in K $.
    \end{description}
\end{rem}
\begin{rem}[Examples of projections]\label{Rem.proj02}
    Based on Remark \ref{Rem.proj01}, we can also see the following facts.
    \begin{description}
        \item[(\hypertarget{Fact3l}{Fact\,3})]If  $ -\infty < r^\ell \leq s^\ell < \infty $, $ \ell = 1, 2 $, then the projections $ \mathrm{proj}_{[r^\ell, s^\ell]} : \R \longrightarrow [r^\ell, s^\ell] $ onto compact intervals $ [r^\ell, s^\ell] \subset \R $ fulfills that:
            \begin{equation*}
                \bigl| \mathrm{proj}_{[r^1, s^1]}(\xi) -\mathrm{proj}_{[r^2, s^2]}(\xi) \bigr| \leq |r^1 -r^2| \vee |s^1 -s^2|, \mbox{ for any $ \xi \in \R $.}
            \end{equation*}
        \item[(\hypertarget{Fact4l}{Fact\,4})]Let $ \mathfrak{K} $ be the class of constraints defined in \eqref{obs_K}, and let $ K = \jump{\kappa^0, \kappa^1} \in \mathfrak{K} $ be the constraint with the obstacles $ \kappa^\ell : Q \longrightarrow [-\infty, \infty] $, $ \ell = 0, 1 $. Then, for the projection $ \mathrm{proj}_K : \mathscr{H} \longrightarrow K $, it holds that:
            \begin{align*}
                \bigl[ \mathrm{proj}_{K} & (u) \bigr](t, x) = \mathrm{proj}_{[\kappa^0(t, x), \kappa^1(t, x)] \cap \R}(u(t, x))
                \\[1ex]
                & = \bigl( \kappa^0 \vee (\kappa^1 \wedge u) \bigr)(t, x) = \left\{ \begin{array}{ll}
                    \kappa^1(t, x), & \mbox{if $ u(t, x) > \kappa^1(t, x) $,}
                    \\[1ex]
                    u(t, x), & \mbox{if $ \kappa^0(t, x) \leq u(t, x) \leq \kappa^1(t, x) $,}
                    \\[1ex]
                    \kappa^0(t, x), & \mbox{if $ u(t, x) < \kappa^0(t, x) $,} 
                \end{array}
                \right.
                \ \ \\
                \\[-3ex]
                & \mbox{ a.e. $(t, x) \in Q$, for any $ u \in \mathscr{H} $.}
            \end{align*}
    \end{description}
\end{rem}

For a proper, lower semi-continuous (l.s.c.), and convex function $ \Psi : X \to (-\infty, \infty] $ on a Hilbert space $ X $, we denote by $ D(\Psi) $ the effective domain of $ \Psi $. Also, we denote by $\partial \Psi$ the subdifferential of $\Psi$. The subdifferential $ \partial \Psi $ corresponds to a weak differential of convex function $ \Psi $, and it is known as a maximal monotone graph in the product space $ X \times X $. The set $ D(\partial \Psi) := \bigl\{ z \in X \ |\ \partial \Psi(z) \neq \emptyset \bigr\} $ is called the domain of $ \partial \Psi $. We often use the notation ``$ [w_{0}, w_{0}^{*}] \in \partial \Psi $ in $ X \times X $\,'', to mean that ``$ w_{0}^{*} \in \partial \Psi(w_{0})$ in $ X $ for $ w_{0} \in D(\partial\Psi) $'', by identifying the operator $ \partial \Psi $ with its graph in $ X \times X $.
\medskip

Next, for Hilbert spaces $X_1, \cdots, X_d$, with $1<d \in \mathbb{N}$, let us consider a proper, l.s.c., and convex function on the product space $X_1 \times \dots \times X_d$:
\begin{equation*}
\widehat{\Psi}: w = [w_1,\cdots,w_d] \in X_1 \times\cdots\times X_d \mapsto \widehat{\Psi}(w)=\widehat{\Psi}(w_1,\cdots,w_d) \in (-\infty,\infty]. 
\end{equation*}
Besides, for any $i \in \{1, \dots, d\}$, we denote by $\partial_{w_i} \widehat{\Psi}:X_1 \times \cdots \times X_d \to X_i$ a set-valued operator, which maps any $w=[w_1,\dots,w_i,\dots,w_d] \in X_1 \times \dots \times X_i \times \dots \times X_d$ to a subset $ \partial_{w_i} \widehat{\Psi}(w) \subset  X_i $, prescribed as follows:
\begin{equation*}
\begin{array}{rl}
\partial_{w_i}\widehat{\Psi}(w)&=\partial_{w_i}\widehat{\Psi}(w_1,\cdots,w_i,\cdots,w_d)
\\[2ex]
&:= \left\{\begin{array}{l|l}\tilde{w}^* \in X_i & \begin{array}{ll}\multicolumn{2}{l}{(\tilde{w}^*,\tilde{w}-w_i)_{X_i} \le \widehat{\Psi}(w_1,\cdots,\tilde{w},\cdots, w_d)}
\\[0.25ex] 
& \quad -\widehat{\Psi}(w_1,\cdots,w_i,\cdots,w_d), \mbox{ for any $\tilde{w} \in X_i$}\end{array}
\end{array}\right\}.
\end{array}
\end{equation*}
As is easily checked, 
\begin{equation}\label{prodSubDif}
    \partial \widehat{\Psi} \subset \bigl[ \partial_{w_1} \widehat{\Psi} \times \cdots \times \partial_{w_d} \widehat{\Psi} \bigr] \mbox{ in $ [X_1 \times \cdots \times X_d]^2 $,}
\end{equation}
where $ \bigl[ \partial_{w_1} \widehat{\Psi} \times \cdots \times \partial_{w_d} \widehat{\Psi} \bigr] : X_1 \times \cdots \times X_d \longrightarrow 2^{X_1 \times \cdots \times X_d} $ is a set-valued operator, defined as:
\begin{equation*}
    \begin{array}{c}
        \ds \bigl[ \partial_{w_1} \widehat{\Psi} \times \cdots \times \partial_{w_d} \widehat{\Psi} \bigr](w) := \partial_{w_1} \widehat{\Psi}(w) \times \cdots \times \partial_{w_d} \widehat{\Psi}(w) \mbox{ in $ X_1 \times \cdots \times X_d $},
        \\[1.5ex]
        \mbox{for any } w = [w_1, \dots, w_d] \in  D\bigl( \bigl[ \partial_{w_1} \widehat{\Psi} \times \cdots \times \partial_{w_d} \widehat{\Psi} \bigr] \bigr) := D(\partial_{w_1} \widehat{\Psi}) \cap \cdots \cap D(\partial_{w_d} \widehat{\Psi}).
    \end{array}
\end{equation*}
But, it should be noted that the converse inclusion of \eqref{prodSubDif} is not true, in general. 

\begin{ex}[Examples of the subdifferential]\label{exConvex}
    As one of the representatives of the subdifferentials, we exemplify the following set-valued signal function $ \Sgn^d: \R^d \rightarrow 2^{\mathbb{R}^d} $, with $ d \in \mathbb{N} $, which is defined as:
\begin{align}\label{Sgn^d}
    \xi = [\xi_1, & \dots, \xi_d] \in \mathbb{R}^d \mapsto \Sgn^d(\xi) = \Sgn^d(\xi_1, \dots, \xi_d) 
    \nonumber
    \\
    & := \left\{ \begin{array}{ll}
            \multicolumn{2}{l}{
                    \ds \frac{\xi}{|\xi|} = \frac{[\xi_1, \dots, \xi_d]}{\sqrt{\xi_1^2 +\cdots +\xi_d^2}}, ~ } \mbox{if $ \xi \ne 0 $,}
                    \\[3ex]
            \mathbb{D}^d, & \mbox{otherwise,}
        \end{array} \right.
    \end{align}
where $ \mathbb{D}^d $ denotes the closed unit ball in $ \mathbb{R}^d $ centered at the origin. Indeed, the set-valued function $ \Sgn^d $ coincides with the subdifferential of the Euclidean norm $ |{}\cdot{}| : \xi \in \mathbb{R}^d \mapsto |\xi| = \sqrt{\xi_1^2 + \cdots +\xi_d^2} \in [0, \infty) $, i.e.:
\begin{equation*}
\partial |{}\cdot{}|(\xi) = \Sgn^d(\xi), \mbox{ for any $ \xi \in D(\partial |{}\cdot{}|) = \mathbb{R}^d $,}
\end{equation*}
and furthermore, it is observed that:
\begin{equation*}
    \partial  |{}\cdot{}|(0) = \mathbb{D}^d \begin{array}{c} \subseteq_{\hspace{-1.25ex}\mbox{\tiny$_/$}}  
    \end{array} [-1, 1]^d 
        = \bigl[ \partial_{\xi_1}  |{}\cdot{}| \times \cdots \times \partial_{\xi_d}  |{}\cdot{}| \bigr](0).
\end{equation*}
\end{ex}
\begin{ex}\label{Rem.f_eps}
    Let $ d \in \mathbb{N} $ be the constant of dimension. For any $ \varepsilon \geq 0 $, let $ f_\varepsilon : \mathbb{R}^d \longrightarrow [0, \infty) $ be a continuous and convex function, defined as:
    \begin{equation}\label{f_eps}
        f_\varepsilon : y \in \mathbb{R}^d \mapsto f_\varepsilon(y) := \sqrt{\varepsilon^2 +|y|^2} \in [0, \infty).
    \end{equation}

    When $ \varepsilon = 0 $, the convex function $ f_0 $ of this case coincides with the $ d $-dimensional Euclidean norm $ |\cdot| $, and hence, the subdifferential $ \partial f_0 $ coincides with the set valued signal function $ \mathrm{Sgn}^d : \R^d \longrightarrow 2^{\mathbb{R}^d} $, defined in \eqref{Sgn^d}. 

    In the meantime, when $ \varepsilon > 0 $, the convex function $ f_\varepsilon  $ belongs to $ C^\infty $-class, and the subdifferential $ \partial f_\varepsilon $ is identified with the (single-valued) usual gradient:
    \begin{equation*}
        \nabla f_\varepsilon : y \in \mathbb{R}^d \mapsto \nabla f_\varepsilon(y) = \frac{y}{\sqrt{\varepsilon^2 +|y|^2}} \in \R^d.
    \end{equation*}

    Moreover, since:
    \begin{align*}
        f_\varepsilon(y) = \bigl| [\varepsilon, y] \bigr|_{\R^{d +1}} & =\bigl| [\varepsilon, y_1, \dots, y_d] \bigr|_{\R^{d +1}}, \mbox{ for all $ [\varepsilon, y] = [\varepsilon, y_1, \dots, y_d] \in \R^{d +1} $,}
        \\
        & \mbox{with $ \varepsilon \geq 0 $ and $ y = [y_1, \dots, y_d] \in \R^d $,}
    \end{align*}
    it will be estimated that:
    \begin{subequations}\label{exM}
        \begin{align}\label{exM00}
            |f_\varepsilon(y) -f_{\tilde{\varepsilon}}(\tilde{y})| \leq & \bigl| [\varepsilon, y] -[\tilde{\varepsilon}, \tilde{y}] \bigr|_{\R^{d +1}} \leq |\varepsilon -\tilde{\varepsilon}| +|y -\tilde{y}|_{\R^d}, 
            \nonumber
            \\
            & \mbox{for all $ \varepsilon, \tilde{\varepsilon} \geq 0  $ and $ y, \tilde{y} \in \R^d $,}
    \end{align}
    \begin{align}\label{exM01}
        & 
        \begin{cases} 
            \ds \bigl| \nabla f_\varepsilon(y) \bigr|_{\R^{d}} = \left| \frac{y}{\bigl| [\varepsilon, y] \bigr|_{\mathbb{R}^{d +1}}} \right|_{\R^{d}} \leq \left| \frac{[\varepsilon, y]}{\bigl| [\varepsilon, y] \bigr|_{\mathbb{R}^{d +1}}} \right|_{\R^{d +1}} = 1,
            \\[3ex]
            \ds \bigl| \nabla f_\varepsilon(y) - \nabla f_{\tilde{\varepsilon}}(\tilde{y}) \bigr|_{\R^d} \leq \left| \frac{[\varepsilon, y]}{\bigl| [\varepsilon, y] \bigr|_{\R^{d +1}}} -\frac{[\tilde{\varepsilon}, \tilde{y}]}{\bigl| [\tilde{\varepsilon}, \tilde{y}] \bigr|_{\R^{d +1}}} \right|_{\R^{d +1}}
        \\
        \ds \qquad \leq \frac{2}{\varepsilon \wedge \tilde{\varepsilon}} \bigl( |\varepsilon -\tilde{\varepsilon}| +|y -\tilde{y}|_{\R^d} \bigr), \mbox{ for all $ \varepsilon, \tilde{\varepsilon} > 0  $ and $ y, \tilde{y} \in \R^d $,}
        \end{cases}
    \end{align}
    and
    \begin{align}\label{exM02}
        & 
        \begin{cases}
            \ds |\nabla^2 f_\varepsilon(y)|_{\R^{d \times d}} \leq \bigl| \nabla^2 \bigl[ |\cdot|_{\R^{d +1}} \bigr]([\varepsilon, y]) \bigr|_{\R^{(d +1) \times (d +1)}} \leq \frac{d +1}{\varepsilon},
            \\[2ex]
            \bigl| \nabla^2 f_\varepsilon(y) - \nabla^2 f_{\tilde{\varepsilon}}(\tilde{y}) \bigr|_{\R^{d \times d}} 
            \\[1ex]
            \ds \qquad \leq \bigl| \nabla^2 \bigl[ |\cdot|_{\R^{d +1}} \bigr]([\varepsilon, y]) -\nabla^2 \bigl[ |\cdot|_{\R^{d +1}} \bigr]([\tilde{\varepsilon}, \tilde{y}])  \bigr|_{\R^{(d +1) \times (d +1)}}
            \\[1ex]
            \ds \qquad \leq \frac{3(d +1)^2}{(\varepsilon \wedge \tilde{\varepsilon})^2} \bigl( |\varepsilon -\tilde{\varepsilon}| +|y -\tilde{y}|_{\R^d} \bigr), \mbox{ for all $ \varepsilon, \tilde{\varepsilon} > 0  $ and $ y, \tilde{y} \in \R^d $.}
        \end{cases}
   \end{align}
    \end{subequations}
\end{ex}
\medskip

Finally, we mention about a notion of functional convergence, known as ``Mosco-convergence''. 
 
\begin{defn}[Mosco-convergence: cf. \cite{MR0298508}]\label{Def.Mosco}
    Let $ X $ be an abstract Hilbert space. Let $ \Psi : X \rightarrow (-\infty, \infty] $ be a proper, l.s.c., and convex function, and let $ \{ \Psi_n \}_{n = 1}^\infty $ be a sequence of proper, l.s.c., and convex functions $ \Psi_n : X \rightarrow (-\infty, \infty] $, $ n = 1, 2, 3, \dots $.  Then, it is said that $ \Psi_n \to \Psi $ on $ X $, in the sense of Mosco, as $ n \to \infty $, iff. the following two conditions are fulfilled:
\begin{description}
    \item[(\hypertarget{M_lb}{M1}) The condition of lower-bound:]$ \ds \varliminf_{n \to \infty} \Psi_n(\check{w}_n) \geq \Psi(\check{w}) $, if $ \check{w} \in X $, $ \{ \check{w}_n  \}_{n = 1}^\infty \subset X $, and $ \check{w}_n \to \check{w} $ weakly in $ X $, as $ n \to \infty $. 
    \item[(\hypertarget{M_opt}{M2}) The condition of optimality:]for any $ \hat{w} \in D(\Psi) $, there exists a sequence \linebreak $ \{ \hat{w}_n \}_{n = 1}^\infty  \subset X $ such that $ \hat{w}_n \to \hat{w} $ in $ X $ and $ \Psi_n(\hat{w}_n) \to \Psi(\hat{w}) $, as $ n \to \infty $.
\end{description}
    As well as, if the sequence of convex functions $ \{ \widehat{\Psi}_\varepsilon \}_{\varepsilon \in \Xi} $ is labeled by a continuous argument $\varepsilon \in \Xi$ with a range $\Xi \subset \mathbb{R}$ , then for any $\varepsilon_{0} \in \Xi$, the Mosco-convergence of $\{ \widehat{\Psi}_\varepsilon \}_{\varepsilon \in \Xi}$, as $\varepsilon \to \varepsilon_{0}$, is defined by those of subsequences $ \{ \widehat{\Psi}_{\varepsilon_n} \}_{n = 1}^\infty $, for all sequences $\{ \varepsilon_n \}_{n=1}^{\infty} \subset \Xi$, satisfying $\varepsilon_{n} \to \varepsilon_{0}$ as $n \to \infty$.
\end{defn}

\begin{rem}\label{Rem.MG}
    Let $ X $, $ \Psi $, and $ \{ \Psi_n \}_{n = 1}^\infty $ be as in Definition~\ref{Def.Mosco}. Then, the following hold.
\begin{description}
    \item[(\hypertarget{Fact5}{Fact\,5})](cf. \cite[Theorem 3.66]{MR0773850} and \cite[Chapter 2]{Kenmochi81}) Let us assume that
    \begin{equation}\label{Mosco01}
    \Psi_n \to \Psi \mbox{ on $ X $, in the sense of  Mosco, as $ n \to \infty $,}
    \vspace{-1ex}
\end{equation}
and
\begin{equation*}
\left\{ ~ \parbox{10cm}{
$ [w, w^*] \in X \times X $, ~ $ [w_n, w_n^*] \in \partial \Psi_n $ in $ X \times X $, $ n \in \N $,
\\[1ex]
$ w_n \to w $ in $ X $ and $ w_n^* \to w^* $ weakly in $ X $, as $ n \to \infty $.
} \right.
\end{equation*}
Then, it holds that:
\begin{equation*}
[w, w^*] \in \partial \Psi \mbox{ in $ X \times X $, and } \Psi_n(w_n) \to \Psi(w) \mbox{, as $ n \to \infty $.}
\end{equation*}
    \item[(\hypertarget{Fact6}{Fact\,6})](cf. \cite[Lemma 4.1]{MR3661429} and \cite[Appendix]{MR2096945}) Let $ d \in \mathbb{N} $ denote dimension constant, and let $  S \subset \R^d $ be a bounded open set. Then, under the Mosco-convergence as in \eqref{Mosco01}, a sequence $ \{ \widehat{\Psi}_n^S \}_{n = 1}^\infty $ of proper, l.s.c., and convex functions on $ L^2(S; X) $, defined as:
        \begin{equation*}
            w \in L^2(S; X) \mapsto \widehat{\Psi}_n^S(w) := \left\{ \begin{array}{ll}
                    \multicolumn{2}{l}{\ds \int_S \Psi_n(w(t)) \, dt,}
                    \\[1ex]
                    & \mbox{ if $ \Psi_n(w) \in L^1(S) $,}
                    \\[2.5ex]
                    \infty, & \mbox{ otherwise,}
                \end{array} \right. \mbox{for $ n = 1, 2, 3, \dots $;}
        \end{equation*}
        converges to a proper, l.s.c., and convex function $ \widehat{\Psi}^S $ on $ L^2(S; X) $, defined as:
        \begin{equation*}
            z \in L^2(S; X) \mapsto \widehat{\Psi}^S(z) := \left\{ \begin{array}{ll}
                    \multicolumn{2}{l}{\ds \int_S \Psi(z(t)) \, dt, \mbox{ if $ \Psi(z) \in L^1(S) $,}}
                    \\[2ex]
                    \infty, & \mbox{ otherwise;}
                \end{array} \right. 
        \end{equation*}
        on $ L^2(S; X) $, in the sense of Mosco, as $ n \to \infty $. 
\end{description}
\end{rem}
\begin{ex}[Example of Mosco-convergence]\label{Rem.ExMG}
    Let $ d \in \N $ be the constant of dimension, and let $ \{ f_\varepsilon \}_{\varepsilon \geq 0} \subset C(\R^d) $ be the sequence of nonexpansive convex functions, as in \eqref{f_eps} and \eqref{exM}. Then, the uniform estimate \eqref{exM00} immediately leads to:
    \begin{equation*}
        f_\varepsilon \to f_{\varepsilon_0} \mbox{ on $ \mathbb{R}^d $, in the sense of Mosco, as $ \varepsilon \to \varepsilon_0 $, for any $ \varepsilon_0 \geq 0 $.}
    \end{equation*}
\end{ex}

\section{Auxiliary results}

In this Section, we prepare some auxiliary results for our study. The auxiliary results are stated in the following two Subsections.
\bigskip

\begin{description}
    \item[\textmd{$\S$\,2.1}]Abstract theory for the state-system (S)$_\varepsilon$;
        \vspace{-1ex}
    \item[\textmd{$\S$\,2.2}]Mathematical theory for the linearized system of (S)$_\varepsilon$. 
\end{description}

\subsection{Abstract theory for the state-system (S)$_\varepsilon$}

In this Subsection, we refer to \cite[Appendix]{MR4218112} to overview the abstract theory of nonlinear evolution equation in an abstract Hilbert space $X$, which enables us to handle the state-systems (S)$_\varepsilon$, for all $ \varepsilon \geq 0 $, in a unified fashion. 
\medskip

The general theory consists of the following two Propositions. 

\begin{prop}[{cf. \cite[Lemma 8.1]{MR4218112}}]\label{Lem.CP}
    Let $ \{ \mathcal{A}_0(t) \, | \, t \in [0, T] \} \subset \mathscr{L}(X) $ be a class of time-dependent bounded linear operators, let $ \mathcal{G}_0 : X \longrightarrow X $ be a given nonlinear operator, and let $ \Psi_0 : X \longrightarrow [0, \infty] $ be a non-negative, proper, l.s.c., and convex function, fulfilling the following conditions: 
    \begin{description}
        \item[\textmd{\it(\hypertarget{cp0}{cp.0})}]$ \mathcal{A}_0(t) \in \mathscr{L}(X) $ is positive and selfadjoint, for any $ t \in [0, T] $, and it holds that
            \begin{equation*}
                (\mathcal{A}_0(t) w, w)_X \geq \kappa_0 |w|_X^2 ,\ \mbox{for any}\ w \in X,
            \end{equation*}
            with some constant $\kappa_{0} \in (0, 1)$, independent of $t \in [0, T]$ and $w \in X$.
        \item[\textmd{\it(\hypertarget{cp1}{cp.1})}]$ \mathcal{A}_0 : [0, T] \longrightarrow \mathscr{L}(X) $ is Lipschitz continuous, so that $ \mathcal{A}_0 $ admits the (strong) time-derivative $ \mathcal{A}_0'(t) \in \mathscr{L}(X) $ a.e. in $ (0, T) $, and  
            \begin{equation*}
                A_T^* := \mathrm{ess} \sup_{\hspace{-3ex}t \in (0, T)} \left\{ \max \{ |\mathcal{A}_0(t)|_{\mathscr{L}(X)}, |\mathcal{A}_0'(t)|_{\mathscr{L}(X)} \} \right\} < \infty;
            \end{equation*}
        \item[\textmd{\it(\hypertarget{cp2}{cp.2})}]$ \mathcal{G}_0 : X \longrightarrow X $ is a Lipschitz continuous operator, and $ \mathcal{G}_0 $ has a $ C^1 $-potential functional $ \widehat{\mathcal{G}}_0 : X \longrightarrow \mathbb{R} $, so that the G\^{a}teaux derivative $ \widehat{\mathcal{G}}_0'(w) \in X^* $ $ (= X) $ at any $ w \in X $ coincides with $ \mathcal{G}_0(w) \in X $;  
        \item[\textmd{\it(\hypertarget{cp3}{cp.3})}]$ \Psi_0 \geq 0 $ on $ X $, and the sublevel set $ \bigl\{ w \in X \, \bigl| \, \Psi_0(w) \leq r \bigr\} $ is compact in $ X $, for any $ r \geq 0 $.
    \end{description}
    Then, for any initial data $ w_0 \in D(\Psi_0) $ and a forcing term $ \mathfrak{f}_0 \in L^2(0, T; X) $, the following Cauchy problem of evolution equation:
    \begin{equation*}
        (\mathrm{CP})~~
        \begin{cases}
            \mathcal{A}_0(t) w'(t) +\partial \Psi_0(w(t)) +\mathcal{G}_0(w(t)) \ni \mathfrak{f}_0(t) \mbox{ in $ X $, ~ $ t \in (0, T) $,}
            \\
            w(0) = w_0 \mbox{ in $ X $;}
        \end{cases}
    \end{equation*}
    admits a unique solution $ w \in L^2(0, T; X) $, in the sense that:
    \begin{equation*}
        w \in W^{1, 2}(0, T; X), ~  \Psi_0(w) \in L^\infty(0, T),
    \end{equation*}
    and
    \begin{equation*}
        \begin{array}{c}
            \displaystyle \bigl( \mathcal{A}_0(t)w'(t) +\mathcal{G}_0(w(t)) -\mathfrak{f}_0(t), w(t) -\varpi \bigr)_X +\Psi_0(w(t)) \leq \Psi_0(\varpi), 
            \\[1ex]
            \mbox{for any $ \varpi \in D(\Psi_0) $, a.e. $ t \in (0, T) $.} 
        \end{array}
    \end{equation*}
    Moreover, both $ t \in [0, T] \mapsto \Psi_0(w(t)) \in [0, \infty) $ and $ t \in [0, T] \mapsto \widehat{\mathcal{G}}_0(w(t)) \in \mathbb{R} $ are absolutely continuous functions in time, and
    \begin{equation*}
        \begin{array}{c}
            \displaystyle |\mathcal{A}_0(t)^{\frac{1}{2}}w'(t)|_X^2 +\frac{d}{dt} \left( \Psi_0(w(t)) +\widehat{\mathcal{G}}_0(w(t)) \right) = (\mathfrak{f}_0(t), w'(t))_X,
            \\[1ex]
            \mbox{for a.e. $ t \in (0, T) $.}
        \end{array}
    \end{equation*}
\end{prop}
\begin{prop}[{cf. \cite[Lemma 8.2]{MR4218112}}]\label{Lem.CP02}
    Under the notations $ \mathcal{A}_0 $, $ \mathcal{G}_0 $, $ \Psi_0 $, and assumptions (\hyperlink{cp0}{cp.0})--(\hyperlink{cp3}{cp.3}), as in the previous Proposition \ref{Lem.CP}, let us fix $ w_0 \in D(\Psi_0) $ and $ \mathfrak{f}_0 \in L^2(0, T; X) $, and take the unique solution $ w \in L^2(0, T; X) $ to the Cauchy problem (CP). Let $ \{ \Psi_n \}_{n = 1}^\infty $, $ \{ w_{0, n} \}_{n = 1}^\infty $, and $ \{ \mathfrak{f}_n \}_{n = 1}^\infty $ be, respectively, a sequence of proper, l.s.c., and convex functions on $ X $, a sequence of initial data in $ X $, and a sequence of forcing terms in $ L^2(0, T; X) $, such that:
\begin{description}
    \item[\textmd{\it(\hypertarget{cp4}{cp.4})}]$ \Psi_n \geq 0 $ on $ X $, for $ n = 1, 2, 3, \dots $, and the union $ \bigcup_{n = 1}^\infty \bigl\{  w \in X \, \bigl| \, \Psi_n(w) \leq r \bigr\} $ of sublevel sets is relatively compact in $ X $, for any $ r \geq 0 $;
    \item[\textmd{\it(\hypertarget{cp5}{cp.5})}]$ \Psi_n $ converges to $ \Psi_0 $ on $ X $, in the sense of Mosco, as $ n \to \infty $; 
    \item[\textmd{\it(\hypertarget{cp6}{cp.6})}]$ \sup_{n \in \mathbb{N}} \Psi_n(w_{0, n}) < \infty $, and $ w_{0, n} \to w_0 $ in $ X $, as $ n \to \infty $;
    \item[\textmd{\it(\hypertarget{cp7}{cp.7})}]$ \mathfrak{f}_n \to \mathfrak{f}_0 $ weakly in $ L^2(0, T; X) $, as $ n \to \infty $.
\end{description}
For any $ n \in \N $, let $ w_n \in L^{2}(0, T; X)$ be the solution to the Cauchy problem (CP), for the initial data $ w_{0, n} \in D(\Psi_n) $ and forcing term $ \mathfrak{f}_n \in L^2(0, T; X) $. Then, 
    \begin{equation*}
        \begin{array}{c}
            \displaystyle w_n \to w\mbox{ in $ C([0, T]; X) $, weakly in $ W^{1, 2}(0, T; X) $,}
            \\[1ex]
            \displaystyle \int_0^T \Psi_n(w_n(t)) \, dt \to \int_0^T \Psi_0(w(t)) \, dt, 
            \mbox{ as $ n \to \infty $,}
        \end{array}
    \end{equation*}
    and
    \begin{equation*}
        \bigl| \Psi_0(w) \bigr|_{C([0, T])} \leq  \sup_{n \in \mathbb{N}} \, \bigl| \Psi_n(w_n) \bigr|_{C([0, T])} < \infty.
    \end{equation*}
\end{prop}

In this paper, the readers are recommended  to see \cite[Appendix]{MR4218112} for the detailed proofs of the above Propositions \ref{Lem.CP} and \ref{Lem.CP02}. Roughly summarized, these Propositions can be obtained by means of modified (mixed and reduced) methods of the existing theories, such as \cite{MR0348562,MR2582280,Kenmochi81}.


\subsection{Mathematical theory for the linearized system of (S)$_\varepsilon$}

In this Subsection, we recall the previous work \cite{MR3888633}, and set up some auxiliary results. In what follows, we let $ \mathscr{Y} := \mathscr{V} \times \mathscr{V}_0 $, with the dual $ \mathscr{Y}^* := \mathscr{V}^* \times \mathscr{V}_0^* $. 
Note that $\mathscr{Y}$ is a Hilbert space which is endowed with a uniform convex topology, based on the inner product for product space, as in the Preliminaries (see the paragraph of Abstract notations).

Besides, we define:
    \begin{equation*}
        \mathscr{Z} := \bigl( W^{1, 2}(0, T; V^*) \cap \mathscr{V} \bigr) \times \bigl( W^{1, 2}(0, T; V_0^*) \cap \mathscr{V}_0 \bigr),
    \end{equation*}
    as a Banach space, endowed with the norm:
    \begin{align*}
        | [\tilde{p}, & \tilde{z}] |_{\mathscr{Z}} := |[\tilde{p}, \tilde{z}]|_{[C([0, T]; H)]^2} +\bigl( |[\tilde{p}, \tilde{z}]|_{\mathscr{Y}}^2 +|[\partial_t \tilde{p}, \partial_t \tilde{z}]|_{\mathscr{Y}^*}^2 \bigr)^{\frac{1}{2}}, \mbox{ for $ [\tilde{p}, \tilde{z}] \in \mathscr{Z} $.}
    \end{align*}

Based on this, let us consider the following linear system of parabolic initial-boundary value problem, denoted by (P):
\bigskip

(P)
\vspace{-0ex}
\begin{equation*}
\left\{\parbox{11cm}{
    $\partial_{t}p - \Delta p + \mu(t, x)p + \lambda(t, x)p + \omega(t, x) \cdot \nabla z = h(t, x)$, $(t, x) \in Q,$
\\[1ex]
$\nabla p(t, x)\cdot n_\Gamma = 0$, $(t, x) \in \Sigma$,
\\[1ex]
$p(0, x) = p_{0}(x)$, $x \in \Omega$;
}\right. 
\end{equation*}
\begin{equation*}
\left\{ \parbox{11cm}{
    $a(t, x)\partial_{t}z + b(t, x)z- \mbox{div} \bigl( A(t, x)\nabla z + \nu^{2} \nabla z + \omega(t, x)p \bigr) $ 
    \\
    $ ~~~~~~= k(t, x)$, $(t, x) \in Q $,
    \\[1ex]
    $ z(t, x) = 0 $, $ (t, x) \in \Sigma $,
    \\[1ex]
    $ z(0, x) = z_{0}(x) $, $ x \in \Omega $.
} \right. 
\end{equation*}
\noindent
This system is studied in \cite{MR3888633} as a key-problem for the G$\hat{\mbox{a}}$teaux differential of the cost $\mathcal{J}_{\varepsilon}$ for $ \varepsilon > 0 $. 
In the context, $[a, b, \mu, \lambda, \omega, A] \in [\mathscr{H}]^{6}$ is a given sextuplet of functions which belongs to a subclass $\mathscr{S} \subset [\mathscr{H}]^{6}$, defined as:
\begin{eqnarray}\label{P01}
    \mathscr{S} := \left\{ \begin{array}{l|l}
        [\tilde{a}, \tilde{b}, \tilde{\mu}, \tilde{\lambda}, \tilde{\omega}, \tilde{A}] \in [\mathscr{H}]^{6} & \hspace{-2ex} \parbox{8cm}{
            \vspace{-2ex}
            \begin{itemize}
                \item $ \tilde{a} \in W^{1, \infty}(Q) $ and $ \log \tilde{a} \in L^{\infty}(Q) $,
                \vspace{-1ex}
                \item $ [\tilde{b}, \tilde{\lambda}]\in [L^{\infty}(Q)]^{2} $,
                \vspace{-1ex}
                \item $ \tilde{\mu} \in L^{\infty}(0, T; H) $ with $ \tilde{\mu} \geq 0 $ a.e. in  $Q$, 
                \vspace{-1ex}
                \item $ \tilde{\omega}\in [L^{\infty}(Q)]^{N} $,
                \vspace{-1ex}
            \item $ \tilde{A} \in [L^{\infty}(Q)]^{N\times N} $, and the value $ \tilde{A}(t, x) \in \mathbb{R}^{N\times N} $ is positive and symmetric matrix,
 for a.e. $ (t, x) \in Q $
            \vspace{-2ex}
            \end{itemize}
    }
    \end{array}
\right\}.
\end{eqnarray}
Also, $[p_{0}, z_{0}] \in [H]^{2}$  and $[h, k] \in \mathscr{Y}^* $ are, respectively, an initial pair and forcing pair, in the system (P).
\medskip

Now, we refer to the previous work \cite{MR3888633}, to recall the key-properties of the system (P), in forms of Propositions. 

\begin{prop}[cf. {\cite[Main Theorem 1 (I-A)]{MR3888633}}]\label{Prop(I-A)}
    For any sextuplet $[a, b, \mu, \lambda, \omega, A] \in \mathscr{S}$, any initial pair $[p_{0}, z_{0}] \in [H]^{2}$, and any forcing pair $[h, k] \in \mathscr{Y}^* $, the system (P) admits a unique solution, in the sense that: 
    \begin{equation}\label{ASY1-01}
        \begin{cases}
            p\in W^{1, 2}(0, T; V^{*})\cap L^{2}(0, T; V) \subset C([0, T]; H),
            \\
            z \in W^{1, 2}(0, T; V_{0}^{*})\cap L^{2}(0, T; V_{0}) \subset C([0, T]; H);
        \end{cases}
    \end{equation}
    \begin{align}\label{ASY1-02}
        \ds \left< \partial_{t}p \right. & \hspace{-0.75ex} \left. (t),  \varphi \right>_{V} + (\nabla p(t), \nabla \varphi)_{[H]^N} +\left< \mu(t)p(t), \varphi \right>_V  
        \nonumber
        \\[0.5ex]
        & + (\lambda(t)p(t) + \omega(t) \cdot \nabla z(t), \varphi)_{H} = \left<h(t), \varphi\right>_{V}, 
        \\[0.5ex]
        & \mbox{ for any $ \varphi \in V $, a.e. $ t \in (0, T) $, subject to $ p(0) = p_0 $ in $ H $;}
        \nonumber
    \end{align}
    and
    \begin{align}\label{ASY1-03}
        \langle \partial_{t}z  & (t), a(t)\psi \rangle_{V_{0}} + (b(t)z(t), \psi)_{H} 
        \nonumber
        \\[0.5ex]
        & +\bigl( A(t)\nabla z(t) + \nu^{2}\nabla z(t) +p(t)\omega(t), \nabla\psi \bigr)_{[H]^N} = \langle k(t), \psi \rangle_{V_{0}},
        \\[0.5ex]
        & \mbox{for any $ \psi \in V_{0} $, a.e. $ t \in (0, T) $, subject to $ z(0) = z_0 $ in $ H $.}
        \nonumber
        \end{align}
\end{prop}
\begin{prop}[cf. {\cite[Main Theorem 1 (I-B)]{MR3888633}}]\label{Prop(I-B)}
    For each $ \ell \in \{ 1, 2 \} $, let us take arbitrary $ [a^\ell, b^\ell, \mu^\ell, \lambda^\ell, \omega^\ell, A^\ell] \in \mathscr{S} $, $ [p_0^\ell, z_0^\ell] \in [H]^2 $, and $ [h^\ell, k^\ell] \in \mathscr{Y}^* $, and let us denote by $ [p^\ell, z^\ell] \in [\mathscr{H}]^2 $ the solution to (P), corresponding to the sextuplet $  [a^\ell, b^\ell, \mu^\ell, \lambda^\ell, \omega^\ell, A^\ell] $, initial pair $ [p_0^\ell, z_0^\ell] $, and forcing pair $ [h^\ell, k^\ell] $. Besides, let $ C_0^* = C_0^*(a^1, b^1, \lambda^1, \omega^1) $ be a positive constant, depending on $a^1, b^1, \lambda^1,$ and $\omega^1$, which is defined as:
    \begin{align}\label{C_0*}
        C_0^* := & \frac{9(1 +\nu^2)}{1 \wedge \nu^2 \wedge \inf a^1(Q)} \cdot \bigl( 1 +(C_{V}^{L^4})^2 +(C_{V}^{L^4})^4 +(C_{V_0}^{L^4})^2 \bigr)
        \nonumber
        \\
        & \quad \cdot 
        \bigl( 1 +|a^1|_{W^{1, \infty}(Q)} +|b^1|_{L^\infty(Q)} +|\lambda^1|_{L^\infty(Q)} +|\omega^1|_{[L^\infty(Q)]^N}^2 \bigr),
    \end{align}
    with use of the constants $ C_{V}^{L^4} > 0 $ and $ C_{V_0}^{L^4} > 0  $ of the respective embeddings $ V \subset L^4(\Omega) $ and $ V_0 \subset L^4(\Omega) $.
    Then, it is estimated that:
    \begin{align}\label{est_I-B} \frac{d}{dt} & \bigl( |(p^1 -p^2)(t)|_H^2 +|\sqrt{a^1(t)}(z^1 -z^2)(t)|_H^2 \bigr)
        \nonumber
        \\
        & \quad +\bigl( |(p^1 -p^2)(t)|_V^2 +\nu^2 |(z^1 -z^2)(t)|_{V_0}^2 \bigr)
        \nonumber
        \\[1ex]
        \leq & 3 C_0^* \bigl( |(p^1 -p^2)(t)|_H^2 +|\sqrt{a^1(t)}(z^1 -z^2)(t)|_H^2 \bigr)
        \\
        & \quad +2C_0^*  \bigl( |(h^1 -h^2)(t)|_{V^*}^2 +|(k^1 -k^2)(t)|_{V_0^*}^2 +R_0^*(t) \bigr),
        \nonumber
        \\[1ex]
        & \mbox{for a.e. $ t \in (0, T) $;}
        \nonumber
    \end{align}
    where
    \begin{align*}
        & R_0^* (t) := |\partial_t z^2(t)|_{V_0^*}^2 \bigl( |a^1 -a^2|_{C(\overline{Q})}^2 +|\nabla (a^1 -a^2)(t)|_{[L^4(\Omega)]^N}^2 \bigr)
        \\
        & \quad +|p^2(t)|_{V}^2 \bigl( |(\mu^1 -\mu^2)(t)|_H^2 +|(\omega^1 -\omega^2)(t)|_{[L^4(\Omega)]^N}^2 \bigr)
        \\
        & \quad +|z^2(t)|_{V_0}^2 \bigl( |(b^1 -b^2)(t)|_{L^4(\Omega)}^2 +|p^2(t)(\lambda^1 -\lambda^2)(t)|_{H}^2 \bigr)
        \\
        & \quad +|\nabla z^2(t) (\omega^1 -\omega^2)(t)|_{H}^2 +|(A^1 -A^2)(t) \nabla z^2(t)|_{[H]^N}^2,
            \\
        & \mbox{for a.e. $ t \in (0, T) $.}
    \end{align*}
\end{prop}
\begin{prop}[cf. {\cite[Corollary 1]{MR3888633}}]\label{ASY_Cor.1}
    For any $ [a, b, \mu, \lambda, \omega, A] \in \mathscr{S} $, let us denote by $\mathcal{P} = \mathcal{P}(a, b, \mu, \lambda, \omega, A) : [H]^2 \times \mathscr{Y}^* \longrightarrow \mathscr{Z} $ a linear operator, which maps any pair of data $ \bigl[ [p_0, z_0], [h, k] \bigr] \in [H]^2 \times \mathscr{Y}^* $ to the solution $ [p, z] \in \mathscr{Z} $ to the corresponding linear system (P), for the sextuplet $ [a, b, \mu, \lambda, \omega, A] $, initial pair $ [p_0, z_0] $, and forcing pair $ [h, k] $. Then, for any sextuplet $ [a, b, \mu, \lambda, \omega, A] \in \mathscr{S} $, there exist positive constants $ M_0^* = M_0^*(a, b, \mu, \lambda, \omega, A) $ and $ M_1^* = M_1^*(a, b, \mu, \lambda, \omega, A) $, depending on $ a $, $ b $, $ \mu $, $ \lambda$, $\omega $, and $A$, such that:
    \begin{equation*}
        \begin{array}{c}
            M_0^* \bigl| \bigl[ [p_0, z_0], [h, k] \bigr] \bigr|_{[H]^2 \times \mathscr{Y}^*} \leq | [p, z] |_{\mathscr{Z}} \leq M_1^* \bigl| \bigl[ [p_0, z_0], [h, k] \bigr] \bigr|_{[H]^2 \times \mathscr{Y}^*},
            \\[1ex]
            \mbox{for all $ [p_0, z_0] \in [H]^2 $, $ [h, k] \in \mathscr{Y}^* $,}
            \\[1ex]
            \mbox{and $ [p, z] = \mathcal{P}(a, b, \mu, \lambda, \omega, A)\bigl[ [p_0, z_0], [h, k] \bigr] \in \mathscr{Z} $,}
        \end{array}
    \end{equation*}
    i.e. the operator $ \mathcal{P} = \mathcal{P}(a, b, \mu, \lambda, \omega, A) $ is an isomorphism between the Hilbert space $ [H]^2 \times \mathscr{Y}^* $ and the Banach space $ \mathscr{Z} $. 
\end{prop}
\begin{prop}[cf. {\cite[Corollary 2]{MR3888633}}]\label{ASY_Cor.2}
    Let us assume: 
    \begin{equation*}
        [a, b, \mu, \lambda, \omega, A] \in \mathscr{S}, ~ \{ [a_n, b_n, \mu_n, \lambda_n, \omega_n, A_n] \}_{n = 1}^\infty \subset \mathscr{S}, 
    \end{equation*}
    \begin{align}\label{ASY01}
        & \hspace{-4ex}[a_n, \partial_t a_n, \nabla a_n, b_n, \lambda_n, \omega_n, A_n] \to [a, \partial_t a, \nabla a, b, \lambda, \omega, A] \ \mbox{weakly-$*$ in }
        \nonumber
        \\
        & \mbox{\small $L^\infty(Q) \times L^\infty(Q) \times [L^\infty(Q)]^N \times L^\infty(Q) \times L^\infty(Q) \times [L^\infty(Q)]^N \times [L^\infty(Q)]^{N \times N} $,}\nonumber
        \\
        & \mbox{and in the pointwise sense a.e. in  $ Q $,} \ \mbox{as}\ n \to \infty, 
    \end{align}
    and
    \begin{equation*}
        \begin{cases}
            \mu_n \to \mu \mbox{ weakly-$*$ in $ L^\infty(0, T; H) $,}
            \\[0ex]
            \mu_n(t) \to \mu(t) \mbox{ in $ H $, for a.e. $ t \in (0, T) $,}
        \end{cases}
        \mbox{as $ n \to \infty $.}
    \end{equation*}
    Let us assume $ [p_0, z_0] \in [H]^2 $, $ [h, k] \in \mathscr{Y}^* $, and let us denote by $ [p, z] \in [\mathscr{H}]^2 $ the solution to (P), for the initial pair $ [p_0, z_0] $ and forcing pair $ [h, k] $. Also, let us assume $ \{ [p_{0, n}, z_{0, n}] \}_{n = 1}^\infty \subset [H]^2 $, $ \{ [h_n, k_n] \}_{n = 1}^\infty \subset \mathscr{Y}^* $, and for any $ n \in \N $, let us denote by $ [p_n, z_n] \in [\mathscr{H}]^2 $ the solution to (P), for the 
    sextuplet $ [a_n, b_n, \mu_n, \lambda_n, \omega_n, A_n] \in \mathscr{S} $, initial pair $ [p_{0, n}, z_{0, n}] $ and forcing pair $ [h_n, k_n] $. Then, the following two items hold.
    \begin{description}
        \item[\textmd{\em(\hypertarget{A}{A})}]The convergence:
            \begin{equation*}
                \begin{cases}
                    [p_{0, n}, z_{0, n}] \to [p_0, z_0] \mbox{ in $ [H]^2 $},
                    \\[1ex]
                    [h_n, k_n] \to [h, k] \mbox{ in $ \mathscr{Y}^* $,}
                \end{cases}
                \mbox{as $ n \to \infty $,}
            \end{equation*}
            implies the convergence:
            \begin{align*}
                [p_n, z_n] & \to [p, z] \mbox{ in $ [C([0, T]; H)]^2 $, and in $ \mathscr{Y} $, as $ n \to \infty $.}
            \end{align*}
        \item[\textmd{\em(\hypertarget{B}{B})}]The following two convergences:
            \begin{equation*}
                \begin{cases}
                    [p_{0, n}, z_{0, n}] \to [p_0, z_0] \mbox{ weakly in $ [H]^2 $},
                    \\[1ex]
                    [h_n, k_n] \to [h, k] \mbox{ weakly in $ \mathscr{Y}^* $,}
                \end{cases}
                \mbox{as $ n \to \infty $,}
            \end{equation*}
            and
            \begin{align*}
                [p_n, z_n] & \to [p, z] \mbox{ in $ [\mathscr{H}]^2 $, weakly in $ \mathscr{Y} $,}
                \\
                & \mbox{and weakly in $ W^{1, 2}(0, T; V^*) \times W^{1, 2}(0, T; V_0^*) $, as $ n \to \infty $,}
            \end{align*}
        are equivalent each other.
    \end{description}
\end{prop}
\begin{rem}\label{Rem.Prev01}
    In the previous work \cite{MR3888633}, one of the essential requirements is to use the continuous embedding $  V \subset L^4(\Omega) $, as in Remark \ref{Rem.Prelim01}, which is satisfied under the restriction $ N \leq 4 $ of the spatial dimension $ N \in \mathbb{N} $. Therefore, under the assumption $ N \in \{2, 3, 4\} $ of this paper, the Propositions \ref{Prop(I-A)}--\ref{ASY_Cor.2} will be applicable, although the previous results as in \cite{MR3888633} were obtained under strict assumption $ N \in \{1, 2, 3\} $.
\end{rem}

Finally, we recall an auxiliary result, which was indirectly obtained in the proof of \cite[Key-Lemma 2]{MR3888633}.
\begin{lem}\label{Lem.ax01}
    Let us assume that $ \hat{\mu} \in L^\infty(0, T; H) $, $ \{ \hat{\mu}_n \}_{n = 1}^\infty \subset L^\infty(0, T; H) $, $ \hat{p} \in \mathscr{V} $, $ \{ \hat{p}_n \}_{n = 1}^\infty \subset \mathscr{V} $,
    \begin{subequations}\label{ax01}
        \begin{equation}\label{ax01a}
            \hat{\mu} \geq 0 \mbox{ and } \hat{\mu}_n \geq 0, \mbox{ a.e. in $ Q $, $ n = 1, 2, 3, \dots $,}
        \end{equation}
        \begin{equation}\label{ax01b}
            \begin{cases}
                \hat{\mu}_n \to \hat{\mu} \mbox{ weakly-$*$ in $ L^\infty(0, T; H) $,}
                \\[0.5ex]
                \hat{\mu}_n(t) \to \hat{\mu}(t) \mbox{ in $ H $, for a.e. $ t \in (0, T) $,}
            \end{cases}
        \end{equation}
        and
        \begin{equation}\label{ax01c}
            \hat{p}_n \to \hat{p} \mbox{ in $ \mathscr{H} $, and weakly in $ \mathscr{V} $, as $ n \to \infty $.}
        \end{equation}
    \end{subequations}
    Then, it holds that:
    \begin{equation}\label{ax10}
        \hat{\mu}_n \hat{p}_n \to \hat{\mu} \hat{p} \mbox{ weakly in $ \mathscr{V}^* $, as $ n \to \infty $.}
    \end{equation}
\end{lem}
\paragraph{Proof.}{
    From \eqref{ax01} and Remark \ref{Rem.Prelim01}, we can see that:
    \begin{equation}\label{ax02}
        \sup_{n \in \N} | \hat{\mu}_n \hat{p}_n |_{\mathscr{V}^*} \leq (C_V^{L^4})^2 \sup_{n \in \N} |\hat{\mu}_n|_{L^\infty(0, T; H)} |\hat{p}_n|_{\mathscr{V}} < \infty,
    \end{equation}
    with the use of the constant $ C_V^{L^4} > 0 $ of embedding $ V \subset L^4(\Omega) $. This implies that:
\begin{description}
    \item[\textmd{(\hypertarget{*0}{$\star$\,0})}]the sequence $ \{ \hat{\mu}_n \hat{p}_n \}_{n = 1}^\infty $ is weakly compact in $ \mathscr{V}^* $.
\end{description}
Also, with \eqref{ax01b} and the dominated convergence theorem \cite[Theorem 10 on page 36]{MR0492147} in mind, we can derive that:
\begin{equation}\label{ax03}
    \hat{\mu}_n \to \hat{\mu} \mbox{ in $ \mathscr{H} $, as $ n \to \infty $. }
\end{equation}

Now, on the basis of (\hyperlink{*0}{$\star$\,0}), let us take any $ \hat{q}^* \in \mathscr{V}^* $, such that $ \hat{q}^* \in \mathscr{V}^* $ is a weak limit of a subsequence of $ \{ \hat{\mu}_n \hat{p}_n \}_{n = 1}^\infty $ (not relabeled), i.e.:
\begin{equation}\label{ax04}
    \hat{\mu}_n\hat{p}_n \to \hat{q}^* \mbox{ weakly in $ \mathscr{V}^* $, as $ n \to \infty $.}
\end{equation}
Besides, by taking subsequences if necessary, \eqref{ax01c} and \eqref{ax03} enable us to say that:
\begin{align}
    \hat{\mu}_n \to \hat{\mu} \mbox{ and } & \hat{p}_n \to \hat{p} \mbox{ in the pointwise sense,}
    \nonumber
    \\
    &\mbox{a.e. in $ Q $, as $ n \to \infty $.}
    \label{ax05}
\end{align}
Additionally, by \eqref{ax01} and Remark \ref{Rem.Prelim01}, we can compute that:
\begin{subequations}\label{ax06}
    \begin{align}
        \bigl| \sqrt{\hat{\mu}_n(t)} & \hat{\varphi}(t) -\sqrt{\hat{\mu}(t)} \hat{\varphi}(t) \bigr|_H^2 \leq \bigl| \sqrt{|\hat{\mu}_n -\hat{\mu}|(t)} \, \hat{\varphi}(t) \bigr|_H^2
        \nonumber
        \\
        & \leq (C_V^{L^4})^2 |(\hat{\mu}_n -\hat{\mu})(t)|_H |\hat{\varphi}(t)|_V^2 \to 0, \mbox{ as $ n \to \infty $,}
        \label{ax06a}
    \end{align}
    \begin{align}
        & \sup_{n \in \N} \bigl| \sqrt{\hat{\mu}_n(t)} \hat{\varphi}(t) -\sqrt{\hat{\mu}(t)} \hat{\varphi}(t) \bigr|_H^2 
        \nonumber
        \\
        & \qquad \leq (C_V^{L^4})^2 \sup_{n \in \N} |(\hat{\mu}_n -\hat{\mu})(t)|_H |\hat{\varphi}(t)|_V^2 < \infty,
        \label{ax06b}
        \\
        & \mbox{for any $ \hat{\varphi} \in \mathscr{V} $, and a.e.  $ t \in (0, T) $,}
        \nonumber
    \end{align}
    and
    \begin{equation}\label{ax06c}
        \sup_{n \in \N} \bigl| \sqrt{\hat{\mu}_n}\hat{p}_n\bigr|_{\mathscr{H}}^2 \leq (C_V^{L^4})^2 \sup_{n \in \N} \bigl\{ |\hat{\mu}_n|_{L^\infty(0, T; H)} |\hat{p}_n|_{\mathscr{V}}^2 \bigr\} < \infty. 
    \end{equation}
\end{subequations}
Taking into account \eqref{ax06}, Remark \ref{Rem.Prelim01}, Lions's lemma \cite[Lemma 1.3 on page 12]{MR0259693}, and the dominated convergence theorem \cite[Theorem 10 on page 36]{MR0492147}, one can observe that:
\begin{equation*}
    \begin{cases}
        \sqrt{\hat{\mu}_n} \varphi \to \sqrt{\hat{\mu}} \varphi \mbox{ in $ \mathscr{H} $,}
        \\
        \sqrt{\hat{\mu}_n} \hat{p}_n \to \sqrt{\hat{\mu}} \hat{p} \mbox{ weakly in $ \mathscr{H} $, as $ n \to \infty $,}
    \end{cases}
\end{equation*}
and therefore, 
\begin{align}
    \langle \hat{\mu}_n \hat{p}_n, \hat{\varphi} \rangle_{\mathscr{V}} = (\sqrt{\hat{\mu}_n} \hat{p}_n, &~ \sqrt{\hat{\mu}_n} \hat{\varphi})_{\mathscr{H}} \to (\sqrt{\hat{\mu}} \hat{p}, \sqrt{\hat{\mu}} \hat{\varphi})_{\mathscr{H}} = \langle \hat{\mu} \hat{p}, \hat{\varphi} \rangle_{\mathscr{V}}
    \nonumber
    \\
    &~ \mbox{as $ n \to \infty $, for any $ \hat{\varphi} \in \mathscr{V} $.}
    \label{ax07}
\end{align}

\eqref{ax04} and \eqref{ax07} imply the uniqueness of the weak limit $ \hat{q}^* = \hat{\mu} \hat{p} $ of subsequences of $ \{ \hat{\mu}_n \hat{p}_n \}_{n = 1}^\infty $ in $ \mathscr{V}^* $. Hence, invoking the separability of the Hilbert space $ \mathscr{V}^* $, we conclude the weak convergence \eqref{ax10} with non-necessity of subsequences.  \qed
}

\section{Main Theorems}

We begin by setting up the assumptions needed in our Main Theorems. All Main Theorems are discussed under the following assumptions.  
\begin{description}
\item[\textmd{(\hypertarget{A1l}{A1})}]
    Let $\nu>0$ be a fixed constant. Let $ [\eta_\mathrm{ad}, \theta_\mathrm{ad}] \in [\mathscr{H}]^2 $ be a fixed pair of the \emph{admissible target profile}. 
\item[\textmd{(\hypertarget{A2l}{A2})}]
    For any $ \varepsilon \geq 0 $, let $ f_\varepsilon : \R^N \longrightarrow [0, \infty) $ be the convex function, defined in \eqref{f_eps}.
    \item[\textmd{(\hypertarget{A3l}{A3})}]
        Let $g : \mathbb{R} \longrightarrow  \mathbb{R}$ be a $C^{1}$-function, which is Lipschitz continuous on $\mathbb{R}$. Also, $g$ has a nonnegative primitive $ 0 \leq G \in C^{2}(\mathbb{R})$, i.e. the derivative $ G'= \frac{dG}{d\eta} $ coincides with $ g $ on $ \R $. Moreover, $g$ satisfies that:
        \begin{equation*}
            \liminf_{\xi \downarrow -\infty}g(\xi) = -\infty \mbox{ and } \limsup_{\xi \uparrow \infty}g(\xi) = \infty.
        \end{equation*}
\item[\textmd{(\hypertarget{A4l}{A4})}]
    Let $\alpha:\mathbb{R} \longrightarrow (0, \infty)$ and $\alpha_{0}: Q \longrightarrow (0, \infty)$ be Lipschitz continuous functions, such that:
        \begin{itemize}
            \item $\alpha \in C^2(\R)$, with the first derivative $ \alpha' = \frac{d \alpha}{d \eta} $ and the second one $ \alpha'' = \frac{d^2 \alpha}{d\eta^2} $;
            \item $\alpha'(0) = 0$, $ \alpha'' \geq 0 $ on $ \R $, and $\alpha \alpha'$ is Lipschitz continuous on $\mathbb{R}$; 
            \item $ \alpha \geq \delta_* $ on $ \R $, and $ \alpha_0 \geq \delta_* $ on $ \overline{Q} $, for some constant $ \delta_* \in (0, 1) $.
        \end{itemize}
\item[\textmd{(\hypertarget{A5l}{A5})}]
    Let $ \mathfrak{K} $ and $ \mathfrak{K}_0 $ be the classes of constraints given in \eqref{obs_K} and \eqref{obs_K0}, respectively, and for any constraint $ K = [\![\kappa^0, \kappa^1]\!] \in \mathfrak{K} $, with the measurable obstacles $ \kappa^\ell : Q \longrightarrow [-\infty, \infty] $, $ \ell = 0, 1 $, let $ \mathscr{U}_\mathrm{ad}^K \subset [\mathscr{H}]^2 $ be a class of admissible controls $ [u, v] $, which is defined as:
        \begin{align*}
            \mathscr{U}_\mathrm{ad}^K ~& := \left\{ \begin{array}{l|l}
                [\tilde{u}, \tilde{v}] \in [\mathscr{H}]^2 & \tilde{u} \in K
            \end{array} \right\}
            \\
            & =  \left\{ \begin{array}{l|l}
                [\tilde{u}, \tilde{v}] \in [\mathscr{H}]^2 & \kappa^0 \leq \tilde{u} \leq \kappa^1 \mbox{ a.e. in $ Q $}
            \end{array} \right\}.
        \end{align*}
\end{description}
Moreover, the following extra assumption will be adopted to verify the dependence of optimal controls with respect to the constraint $ K = \jump{\kappa^0, \kappa^1} \in \mathfrak{K} $. 
\begin{description}
\item[\textmd{(\hypertarget{A6l}{A6})}]
    The constraint $ K = \jump{\kappa^0, \kappa^1} \in \mathfrak{K} $ satisfies that: 
        \begin{align*}
            & \kappa^\ell \in L^1\bigl(Q \setminus |\kappa^\ell|^{-1}(\infty) \bigr), \mbox{ with}
            \\
            \mbox{$ |\kappa^\ell|^{-1}(\infty) := $} & \mbox{$ \left\{ \begin{array}{l|l}
                (t, x) \in Q & |\kappa^\ell|(t, x) = \infty
            \end{array} \right\} $ , for $ \ell = 0, 1 $,}
        \end{align*}
        and $ \{ K_n \}_{n = 1}^\infty =  \{\jump{\kappa_n^0, \kappa_n^1} \}_{n = 1}^\infty \subset \mathfrak{K} $ is a sequence of constraints 
        such that:
\begin{align*}
    \kappa_n^\ell(t, x) & \to \kappa^\ell(t, x) ~ (\in [-\infty, \infty]) \mbox{ as $ n \to \infty $,}
    \\
    & \mbox{for a.e. $ (t, x) \in Q $, and $ \ell = 0, 1 $,}
\end{align*}
\begin{align*}
    \int_{Q \setminus |\kappa^\ell|^{-1}(\infty)} |\kappa_n^\ell -\kappa^\ell| \, dxdt \to 0 \mbox{ as $ n \to \infty $, for $ \ell = 0, 1 $,}
\end{align*}
        and moreover, $ \bigcap_{n=1}^{\infty}K_n \neq \emptyset $, i.e. there exists $ \bar{\kappa} \in \mathscr{H} $ satisfying
        \begin{center}
        \parbox{10cm}{
        \begin{equation*}
            \begin{cases}
                \kappa_n^0 \leq \bar{\kappa} \leq \kappa_n^1 \mbox{ a.e. in $ Q $, for $ n = 1, 2, 3, \dots $,}
                \\
                \kappa^0 \leq \bar{\kappa} \leq \kappa^1 \mbox{ a.e. in $ Q $.}
            \end{cases}
        \end{equation*}
    }
        \end{center}
\end{description}
\begin{rem}\label{Rem.alp''}
    The assumption (\hyperlink{A4l}{A4}) leads to the boundedness of the second derivative $ \alpha'' $ of $ \alpha $. In fact, from the Lipschitz continuity of $ \alpha $ and $ \alpha \alpha' $, one can see that:
    \begin{align*}
        |\alpha''(\eta)| \leq \frac{1}{\delta_*} \bigl( \bigl|{\ts \frac{d}{d\eta} (\alpha \alpha')} \bigr|_{L^\infty(\R)} +|\alpha'|_{L^\infty(\R)}^2 \bigr) < \infty,  \mbox{ for any $ \eta \in \R $.}
    \end{align*}
\end{rem}

Now, the Main Theorems of this paper are stated as follows.
\begin{mTh}\label{mainTh01}
    Under the assumptions (\hyperlink{A1l}{A1})--(\hyperlink{A4l}{A4}), let us fix a constant $ \varepsilon \geq 0 $, an initial pair $ [\eta_0, \theta_0] \in D $, 
    and a forcing pair $ [u, v] \in [\mathscr{H}]^2 $. Then, the following hold.
    \begin{description}
        \item[\textmd{(\hypertarget{I-A}{I-A})}]The state-system (S)$_{\varepsilon}$ admits a unique solution $[\eta, \theta] \in [\mathscr{H}]^{2}$, in the sense that:
        \begin{equation}\label{S1}
            \begin{cases}
                \eta \in  W^{1, 2}(0, T; H) \cap L^{\infty}(0, T; V) \cap L^2(0, T; H^2(\Omega)) \subset C([0, T]; H) ,
            \\
                \theta \in  W^{1, 2}(0, T; H) \cap L^{\infty}(0, T; V_0) \subset C([0, T]; H),
        \end{cases}
        \end{equation}
        \begin{equation}\label{S2}
                \begin{array}{c}
                    \ds\bigl( \partial_{t}\eta(t), \varphi \bigr)_{H} + \bigl( \nabla\eta(t), \nabla\varphi \bigr)_{[H]^N} +\bigl( g(\eta(t)), \varphi \bigr)_{H} 
                    \\[1ex]
                    +\bigl( \alpha'(\eta(t))f_{\varepsilon}(\nabla\theta(t)), \varphi \bigr)_{H}  = \bigl( M_u u(t), \varphi \bigr)_{H}, 
                    \\[1ex]
                    \mbox{for any $ \varphi \in V $, a.e. $ t \in (0, T) $, subject to } \eta(0) = \eta_0 \mbox{ in $ H $;}
                \end{array}
        \end{equation}
and
        \begin{equation}\label{S3}
                \begin{array}{c}
                    \ds\bigl( \alpha_{0}(t)\partial_{t}\theta(t), \theta(t)-\psi \bigr)_{H} + \nu^2 \bigl( \nabla\theta(t), \nabla(\theta(t)-\psi) \bigr)_{[H]^N} 
                    \\[1ex]
                    \ds +\int_{\Omega}\alpha(\eta(t))f_{\varepsilon}(\nabla\theta(t))dx \leq \int_{\Omega}\alpha(\eta(t))f_{\varepsilon}(\nabla\psi)dx
                    \\[2ex]
                     +\bigl( M_v v(t), \theta(t)-\psi \bigr)_{H}, \mbox{ for any } \psi \in V_{0},
                     \\[1ex]
                     \mbox{a.e. } t \in (0, T), \mbox{ subject to $ \theta(0) = \theta_0 $ in $ H $.} 
                \end{array}
        \end{equation}
        In particular,  if $ \eta_0 \in L^\infty(\Omega) $ and $ u \in L^\infty(Q) $, then $ \eta \in L^\infty(Q) $. 
    \item[\textmd{(\hypertarget{I-B}{I-B})}]Let $\{\varepsilon_{n} \}_{n=1}^{\infty} \subset [0, \infty)$, $\{[\eta_{0, n}, \theta_{0, n}] \}_{n=1}^{\infty} \subset D $, and $\{ [u_{n}, v_{n}]\}_{n=1}^{\infty} \subset [\mathscr{H}]^2$ be given sequences such that:
\begin{equation}\label{w.i01}
    \varepsilon_{n} \to \varepsilon,\ [\eta_{0, n}, \theta_{0, n}] \to [\eta_{0}, \theta_{0}] \mbox{ weakly in } V \times V_{0}, 
\end{equation}
\begin{equation}\label{w.i02}
    \mbox{and } [M_u u_{n}, M_v v_{n}] \to [M_u u, M_v v] \ \mbox{weakly in}\ [\mathscr{H}]^{2}, \mbox{ as } n \to \infty.
\end{equation}
            In addition, let $[\eta, \theta]$ be the unique solution to (S)$_{\varepsilon}$, for the initial pair $ [\eta_0, \theta_0] $ and forcing pair $[u, v]$. Also, for any $n \in \N$, let $[\eta_{n}, \theta_{n}]$ be the unique solution to (S)$_{\varepsilon_n}$, for the initial pair $ [\eta_{0, n}, \theta_{0, n}] $ and forcing pair $ [u_{n}, v_{n}] $. Then, it holds that:
\begin{align}\label{mThConv}
    [\eta_{n}, & \theta_{n}] \to [\eta, \theta] \mbox{ in } [C([0, T]; H)]^2, \mbox{ in } \mathscr{Y}, \mbox{ weakly in $ [W^{1, 2}(0, T; H)]^2 $,}
    \\
    & \mbox{and weakly-$*$ in $ L^\infty(0, T; V) \times L^\infty(0, T; V_0) $, as $ n \to \infty $.}
    \nonumber
\end{align}
In particular, if:
\begin{align}\label{w.i03}
    \begin{cases}
        \{\eta_{0, n}\}_{n = 1}^\infty \subset L^\infty(\Omega), ~ \{ u_n \}_{n = 1}^\infty \subset L^\infty(Q), 
        \\[0.5ex]
        \sup_{n \in \N} |\eta_{0, n}|_{L^\infty(\Omega)} \vee \sup_{n \in \N} |u_n|_{L^\infty(Q)} < \infty,
    \end{cases}
\end{align}
    then 
        \begin{align}\label{mThConv00}
            \eta_n \to \eta \ \mbox{weakly-$*$ in $ L^\infty(Q) $, as $ n \to \infty $.}
        \end{align}
    \end{description}
\end{mTh}
\begin{rem}\label{Rem.mTh01Conv}
    As a consequence of \eqref{mThConv} and Remark \ref{Rem.alp''}, we further find a subsequence $ \{ n_i \}_{i = 1}^\infty \subset \{n\} $, such that:
    \begin{align*}
        [\eta_{n_i}, \theta & _{n_i} ] \to [\eta, \theta], ~ [\nabla \eta_{n_i}, \nabla \theta_{n_i}] \to [\nabla \eta, \nabla \theta], 
        \\
        & \mbox{in the pointwise sense a.e. in $ Q $,}
    \end{align*}
    \begin{align*}
        \alpha''( \eta_{n_i}) f_{\varepsilon_{n_i}} ( & \nabla \theta_{n_i}) \to \alpha''(\eta) f_\varepsilon(\nabla \theta) \mbox{ weakly-$*$ in $ L^\infty(0, T; H) $,}
        \\
        & \mbox{and in the pointwise sense a.e. in $ Q $,}
    \end{align*}
    and
    \begin{equation*}
        \begin{array}{l}
            [\eta_{n_i}(t), \theta_{n_i}(t)] \to [\eta(t), \theta(t)] \mbox{ in $ V \times V_0 $,}
            \\[1ex]
            \qquad \qquad \mbox{and } \alpha''(\eta_{n_i}(t)) f_{\varepsilon_{n_i}}(\nabla \theta_{n_i}(t)) \to \alpha''(\eta(t)) f_\varepsilon(\nabla \theta(t)) \mbox{ in $ H $,}
            \\[1ex]
            \qquad \qquad \qquad \qquad \mbox{for a.e. $ t \in (0, T) $, as $ i \to \infty $.}
        \end{array}
    \end{equation*}
\end{rem}
\begin{mTh}\label{mainTh02}
    Let us assume (\hyperlink{A1l}{A1})--(\hyperlink{A5l}{A5}). Let us fix any constant $ \varepsilon \geq 0 $, any initial data $ [\eta_0, \theta_0] \in D $, and any constraint $ K = \jump{\kappa^0, \kappa^1} \in \mathfrak{K} $. Then, the following two items hold.
    \begin{description}
        \item[\textmd{(\hypertarget{II-A}{II-A})}]The problem (\hyperlink{OP}{OP})$_{\varepsilon}^K$ has at least one optimal control $[u^{*}, v^{*}] \in \mathscr{U}_\mathrm{ad}^K $, so that:
\begin{equation*}
    \mathcal{J}_{\varepsilon}(u^{*}, v^{*}) = \min \left\{ \begin{array}{l|l}
        \mathcal{J}_{\varepsilon}(u, v) & [u, v] \in \mathscr{U}_\mathrm{ad}^K 
    \end{array} \right\}.
\end{equation*}
        \item[\textmd{(\hypertarget{II-B}{II-B})}]Let us assume the extra assumption (\hyperlink{A6l}{A6}), for the sequence of constraints $ \{ K_n \}_{n = 1}^\infty = \{ \jump{\kappa_n^0, \kappa_n^1} \}_{n = 1}^\infty \subset \mathfrak{K} $, and let us take the sequences $\{\varepsilon_{n} \}_{n=1}^{\infty} \subset [0, \infty)$ and $\{[\eta_{0, n}, \theta_{0, n}] \}_{n=1}^{\infty} $ $ \subset D $ as in \eqref{w.i01}. In addition, for any $n \in \mathbb{N}$, let $[u_{n}^*, v_{n}^*] \in \mathscr{U}_\mathrm{ad}^{K_n}$ be the optimal control of (\hyperlink{OP}{OP})$_{\varepsilon_{n}}^{K_n}$ in the case when the initial pair of corresponding state system (S)$_{\varepsilon_n}$ is given by $ [\eta_{0, n}, \theta_{0, n}] $. Then, there exist a subsequence $ \{ n_{i} \}_{i = 1}^{\infty} \subset \{ n \} $ and a pair of functions $[u^{**}, v^{**}] \in \mathscr{U}_\mathrm{ad}^K$, such that: 
\begin{equation*}
    \left\{ \hspace{-2ex} \parbox{12cm}{
        \vspace{-1ex}
        \begin{itemize}
            \item $ [M_u u_{n_{i}}^*, M_v v_{n_{i}}^*] \to [M_u u^{**}, M_v v^{**}] $  weakly in $ [\mathscr{H}]^{2} $,  as $ i \to \infty $,

            \item $ [u^{**}, v^{**}] $ is an optimal control of (\hyperlink{OP}{OP})$_{\varepsilon}^K$.
        \vspace{-1ex}
        \end{itemize}
    } \right.
\end{equation*}
\end{description}
\end{mTh}
\begin{mTh}
    \label{mainTh03}
    In addition to the assumptions (\hyperlink{A1l}{A1})--(\hyperlink{A5l}{A5}), let us suppose the restricted situation (\hyperlink{rs0}{r.s.0}) as in the Introduction, i.e.:
    \begin{description}
        \item[\textmd{(\hyperlink{rs0}{r.s.0})}]$ \varepsilon > 0 $, $ [\eta_0, \theta_0] \in D_0 $,  and $ K = \jump{\kappa^0, \kappa^1} \in \mathfrak{K}_0 ~(= \mathfrak{K} \cap 2^{L^\infty(Q)}) $.
    \end{description}
    Let $[u^{*}, v^{*}] \in \mathscr{U}_\mathrm{ad}^K$ be an optimal control of (\hyperlink{OP}{OP})$_{\varepsilon}^K$, and let $ [\eta_\varepsilon^*, \theta_\varepsilon^*] $ be the solution to (S)$_{\varepsilon}$, for the initial pair $[\eta_0, \theta_0] $ and forcing pair $ [u^{*}, v^{*}]$. Then, the following two items hold.
    \begin{description}
        \item[\textmd{\textit{(\hypertarget{III-A}{III-A})}}](Necessary condition for (\hyperlink{OP}{OP})$_\varepsilon^K$ under $ \varepsilon > 0 $ and $ K \in \mathfrak{K}_0 $) 
            For the optimal control $ [u^*, v^*] \in \mathscr{U}_\mathrm{ad}^K $ of (\hyperlink{OP}{OP})$_\varepsilon^K$, it holds that:
    \begin{subequations}\label{Thm.5-00}
      \begin{equation}\label{Thm.5-001}
        \ds M_u \left(u^* -\mathrm{proj}_{K}(-p_\varepsilon^*) \right)= 0, \mbox{in}\ \mathscr{H},     \end{equation}
      \begin{equation}\label{Thm.5-002}
        \ds M_v (v^* +  z_\varepsilon^*) = 0\ \mbox{in}\ \mathscr{H}.
      \end{equation} 
    \end{subequations}
    
    In this context,  $[p_\varepsilon^*, z_\varepsilon^*]$ is a unique solution to the following variational system:
\begin{align}\label{Thm.5-01}
    -\bigl< \partial_{t} p_\varepsilon^*(t), & \varphi \bigr>_{V} + \bigl( \nabla p_\varepsilon^*(t), \nabla \varphi \bigr)_{[H]^N} + \bigl< [\alpha''(\eta_\varepsilon^*)f_{\varepsilon}(\nabla\theta_\varepsilon^*)](t) p_\varepsilon^*(t), \varphi \bigr>_{V}
    \nonumber
    \\[0.5ex]
    & +\bigl( g'(\eta_\varepsilon^*(t)) p_\varepsilon^*(t), \varphi \bigr)_{H} +\bigl( [\alpha'(\eta_\varepsilon^*) \nabla f_{\varepsilon}(\nabla\theta_\varepsilon^*)](t) \cdot \nabla z_\varepsilon^*(t), \varphi \bigr)_{H}
    \\[0.5ex]
    = & \bigl( M_\eta (\eta_\varepsilon^*-\eta_{\mbox{\scriptsize ad}})(t), \varphi \bigr)_{H}, \mbox{ for any $ \varphi \in V $,  and a.e. $ t \in (0, T) $;}
    \nonumber
\end{align}
and
\begin{align}\label{Thm.5-02}
    -\bigl< \partial_{t} & \bigl( \alpha_{0}z_\varepsilon^* \bigr)(t), \psi \bigr>_{V_{0}}  +\bigl( [\alpha(\eta_\varepsilon^*)\nabla^2 f_{\varepsilon}(\nabla\theta_\varepsilon^*)](t)\nabla z_\varepsilon^*(t) + \nu^2 \nabla z_\varepsilon^*(t), \nabla \psi \bigl)_{[H]^N}
    \nonumber
    \\[1ex]
    & +\bigl( p_\varepsilon^*(t)[\alpha'(\eta_\varepsilon^*)\nabla f_{\varepsilon}(\nabla\theta_\varepsilon^*)](t), \nabla \psi \bigr)_{[H]^N} = \bigl( M_\theta (\theta_\varepsilon^*-\theta_{\mbox{\scriptsize ad}})(t), \psi \bigr)_{H},
    \\[1ex]
    & \mbox{for any}\ \psi \in V_{0},\ \mbox{and a.e.}\ t \in (0, T);
    \nonumber
\end{align}
subject to the terminal condition:
\begin{equation}\label{Thm.5-03}
[p_\varepsilon^*(T), z_\varepsilon^*(T)] = [0, 0]\ \mbox{in}\ [H]^{2}.
\end{equation}
\item[\textmd{\textit{(\hypertarget{III-B}{III-B})}}]Let $ \{ \varepsilon_n \}_{n = 1}^\infty \subset [0, \infty) $ and $ \{ [\eta_{0, n}, \theta_{0, n}] \}_{n = 1}^\infty \subset D $ be sequences as in \eqref{w.i01}. Also, let $ \{ K_n \}_{n = 1}^\infty = \{\jump{\kappa_n^0, \kappa_n^1}\}_{n = 1}^\infty \subset \mathfrak{K} $ be a sequence of constraints, fulfilling (\hyperlink{A6l}{A6}). In addition, let us assume:
    \begin{equation}\label{IIIB00}
        \left\{ \hspace{-2ex} \parbox{9.5cm}{
            \vspace{-1ex}
            \begin{itemize}
                \item $ \{K_n\}_{n = 1}^\infty = \{ \jump{\kappa_n^0, \kappa_n^1} \}_{n = 1}^\infty \subset \mathfrak{K}_0 ~ \bigl(  = \mathfrak{K} \cap 2^{L^\infty(Q)} \bigr) $,
                \item $ \{ [\eta_{0, n}, \theta_{0, n}] \}_{n = 1}^\infty \subset D_0 ~ \bigl( = (V \cap L^\infty(\Omega)) \times V_0 \bigr) $,
                \item $ \ds \sup_{n \in \N} \bigl\{ |\eta_{0, n}|_{L^\infty(\Omega)} \vee |\kappa_n^0|_{L^\infty(Q)} \vee |\kappa_n^1|_{L^\infty(Q)} \bigr\} < \infty $. 
            \vspace{-1ex}
            \end{itemize}
        } \right.
    \end{equation}
            Then, the subsequence $ \{ n_i \}_{i = 1}^\infty \subset \{n\} $ and the limiting optimal control $ [u^{**}, v^{**}] \in \mathscr{U}_\mathrm{ad}^K $ as in Main Theorem \ref{mainTh02} (\hyperlink{II-B}{II-B}) fulfill that:
    \begin{subequations}\label{IIIB01}
        \begin{equation}\label{IIIB01a}
            u^{**} \in L^\infty(Q), ~ M_v v^{**} \in W^{1, 2}(0, T; V_0^*) \cap L^\infty(0, T; V_0) \subset C([0, T]; H),
        \end{equation}
        \begin{equation}\label{IIIB01b}
            \begin{cases}
                M_u u_{n_i}^{*} \to M_u u^{**} \mbox{ in $ \mathscr{H} $,}
                \\[1ex]
                u_{n_i}^* \to u^{**} \mbox{ weakly-$*$ in $ L^\infty(Q) $,}
            \end{cases}
            \mbox{as $ i \to \infty $,}
        \end{equation}
        and
        \begin{align}\label{IIIB01c}
            M_v v_{n_i}^* & \to M_v v^{**} \mbox{ in $ C([0, T]; H) $, in $ \mathscr{V}_0 $,}
            \nonumber
            \\
            & \mbox{weakly in $ W^{1, 2}(0, T; V_0^*) $, as $ i \to \infty $.}
        \end{align}
    \end{subequations}
\end{description}
\end{mTh}
\begin{rem}\label{Rem.mTh03}
    Let $ \mathcal{R}_T \in \mathscr{L}(\mathscr{H}) $ be an isomorphism, defined as:
    \begin{equation*}
        \bigl( \mathcal{R}_T \varphi \bigr)(t) := \varphi(T -t) \mbox{ in $ H $, for a.e. $ t \in (0, T) $.}
    \end{equation*}
    Also, let us fix $ \varepsilon > 0 $, and define a bounded linear operator $ \mathcal{Q}_\varepsilon^* : [\mathscr{H}]^2 \longrightarrow \mathscr{Z} $ as the restriction $\mathcal{P}|_{\{[0, 0]\} \times \mathscr{Y}^{*}} $ of the linear isomorphism $ \mathcal{P} = \mathcal{P}(a, b, \mu, \lambda, \omega, A) : [H]^2 \times \mathscr{Y}^* \longrightarrow \mathscr{Z} $, as in Proposition \ref{ASY_Cor.1}, in the case when: 
    \begin{equation}\label{setRem4}
        \begin{cases}
            [a, b] = \mathcal{R}_T [\alpha_0, -\partial_t \alpha_0] \mbox{ in $ W^{1, \infty}(Q) \times L^\infty(Q) $,}
            \\[1ex]
            \mu = \mathcal{R}_T \bigl[ \alpha''(\eta_\varepsilon^*) f_\varepsilon(\nabla \theta_\varepsilon^*) \bigr] \mbox{ in $ L^\infty(0, T; H) $,}
            \\[1ex]
            \lambda = \mathcal{R}_T \bigl[ g'(\eta_\varepsilon^*)\bigl] \mbox{ in $ L^\infty(Q) $,}
            \\[1ex]
            \omega = \mathcal{R}_T \bigl[ \alpha'(\eta_\varepsilon^*) \nabla f_\varepsilon(\nabla \theta_\varepsilon^*) \bigl] \mbox{ in $ [L^\infty(Q)]^N $,}
            \\[1ex]
            A = \mathcal{R}_T \bigl[ \alpha(\eta_\varepsilon^*)\nabla^2 f_\varepsilon(\nabla \theta_\varepsilon^*) \bigl] \mbox{ in $ [L^\infty(Q)]^{N\times N} $.}
        \end{cases}
    \end{equation}
    On this basis, let us define:
    \begin{equation}\label{Peps*}
            \mathcal{P}_\varepsilon^* := \mathcal{R}_T \circ \mathcal{Q}_\varepsilon^* \circ \mathcal{R}_T \mbox{ in $ \mathscr{L}([\mathscr{H}]^2; \mathscr{Z}) $.}
    \end{equation}
    Then, since the embedding $ V_0 \subset L^4(\Omega) $ and $\alpha_0\in W^{1,\infty}(Q)$ guarantee:
    \begin{equation}\label{productrule}
        \partial_t (\alpha_0 \tilde{z}) = \alpha_0 \partial_t \tilde{z} +\tilde{z} \partial_t \alpha_0 \mbox{ in $ \mathscr{V}_0^* $, for any $ \tilde{z} \in W^{1, 2}(0, T; V^*_0) $,}
    \end{equation}
    we can obtain the unique solution $[p_{\varepsilon}^{*}, z_{\varepsilon}^{*}] \in [\mathscr{H}]^2$ to the variational system \eqref{Thm.5-01}--\eqref{Thm.5-03} as follows:
\begin{equation*}
    [p_\varepsilon^*, z_\varepsilon^*] = \mathcal{P}_\varepsilon^* \bigl[ M_\eta  (\eta_\varepsilon^* -\eta_\mathrm{ad}), M_\theta  (\theta_\varepsilon^* -\theta_\mathrm{ad}) \bigr] \mbox{ in $ \mathscr{Z} $.}
\end{equation*}
\end{rem}
    \begin{mTh}\label{mainTh04}
        Let us assume (\hyperlink{A1l}{A1})--(\hyperlink{A5l}{A5}), and let us assume that the situation is not under (\hyperlink{rs0}{r.s.0}), i.e. it is under:
        \begin{description}
            \item[\textmd{$\neg$(\hypertarget{neg-rs0}{r.s.0})}]either $ \varepsilon = 0 $, or $ [\eta_0, \theta_0] \in D \setminus D_0 $, or $ K = \jump{\kappa^0,\kappa^1} \in \mathfrak{K} \setminus \mathfrak{K}_0 $ is satisfied.
        \end{description}
        Also, let us define a Hilbert space $ \mathscr{W}_0 $ as follows:
        \begin{equation*}
            \mathscr{W}_0 := \left\{ \begin{array}{l|l}
                \psi \in W^{1, 2}(0, T; H) \cap \mathscr{V}_0 & \psi(0) = 0 \mbox{ in $ H $}
            \end{array} \right\}.
        \end{equation*}
        Then, there exists an optimal control $ [u^{\circ}, v^{\circ}] \in \mathscr{U}_\mathrm{ad}^K $ of the problem (\hyperlink{OP}{OP})$_\varepsilon^K $, together with the solution $ [\eta^\circ, \theta^\circ] $ to the system (S)$_\varepsilon$, for the initial pair $[\eta_0, \theta_0] $ and forcing pair $ [u^\circ, v^\circ] \in \mathscr{U}_\mathrm{ad}^K $, and moreover, there exist pairs of functions $ [p^{\circ}, z^{\circ}] \in \mathscr{Y} $, $ [\xi^\circ, \sigma^\circ] \in \mathscr{H} \times [L^\infty(Q)]^N $, and a distribution $ \zeta^\circ \in \mathscr{W}_0^* $, such that:
\begin{subequations}\label{Thm.5-10}
      \begin{equation}\label{Thm.5-101}
        M_u \left(u^\circ -\mathrm{proj}_{K}(-p^\circ) \right)= 0,\ \mbox{in}\ \mathscr{H},
      \end{equation}
      \begin{equation}\label{Thm.5-102}
       M_v( v^\circ +z^\circ) = 0\ \mbox{in}\ \mathscr{H},
      \end{equation}
    \end{subequations}
\begin{subequations}\label{Thm.5-13}
    \begin{equation}\label{Thm.5-131}
        p^\circ \in W^{1, 2}(0, T; V^*) \cap \mathscr{V} \subset C([0, T]; H), 
    \end{equation}
    \begin{equation}\label{Thm.5-132}
        \sigma^\circ \in \partial f_\varepsilon(\nabla \theta^\circ), \mbox{ a.e. in $ Q $;}
    \end{equation}
\end{subequations}
    \begin{align}\label{Thm.5-11}
        \bigl< -\partial_t & p^\circ, \varphi \bigr>_{\mathscr{V}} +\bigl( \nabla p^\circ, \nabla \varphi \bigr)_{L^2(0, T; [H]^N)}  +\bigl< \alpha''(\eta^\circ) f_\varepsilon(\nabla \theta^\circ) p^\circ, \varphi \bigl>_\mathscr{V}
        \nonumber
        \\
        & +\bigl( g'(\eta^\circ)p^\circ +\mbox{$ \alpha'(\eta^\circ)\xi^\circ $}, \varphi \bigr)_{\mathscr{H}} = \bigl( M_\eta  (\eta^\circ -\eta_\mathrm{ad}), \varphi \bigr)_{\mathscr{H}}, 
        \\
        & 
        \mbox{for any $ \varphi \in \mathscr{V} $, subject to $ p^\circ(T) = 0 $ in $ H $;}
        \nonumber
    \end{align}
    and
    \begin{align}\label{Thm.5-12}
        \bigl( \alpha_0 z^\circ & , \partial_t \psi \bigr)_{\mathscr{H}} +\bigl< \zeta^\circ, \psi \bigr>_{\mathscr{W}_0} +\bigl( \nu^2 \nabla z^\circ + \alpha'(\eta^\circ) \sigma^\circ p^\circ, \nabla \psi \bigr)_{L^2(0, T; [H]^N)}
        \nonumber
        \\
        & = \bigl( M_\theta (\theta^\circ -\theta_\mathrm{ad}), \psi \bigr)_{\mathscr{H}}, \mbox{ for any $ \psi \in \mathscr{W}_0 $.}
    \end{align}
        In particular, if $ \varepsilon > 0 $, i.e. the situation is under:
        \begin{description}
            \item[\textmd{(\hypertarget{rs1}{r.s.1})}]$ \varepsilon > 0 $, while either $ [\eta_0, \theta_0] \in D \setminus D_0 $ or $ K = \jump{\kappa^0, \kappa^1} \in \mathfrak{K} \setminus \mathfrak{K}_0 $ is satisfied;
        \end{description}
        then:
        \begin{align}
            & \begin{cases} 
                \sigma^\circ = \nabla f_\varepsilon(\nabla \theta^\circ), \mbox{ a.e. in $ Q $,}
                \\[0.5ex]
                \xi^\circ = \sigma^\circ \cdot \nabla z^\circ =  \nabla f_\varepsilon(\nabla \theta^\circ) \cdot \nabla z^\circ \mbox{ in $ \mathscr{H} $,}
                \\[0.5ex]
                \zeta^\circ = -\mathrm{div} \bigl( \alpha(\eta^\circ) \nabla^2 f_\varepsilon(\nabla \theta^\circ) \nabla z^\circ \bigr) \mbox{ in $ \mathscr{W}_0^* $.}
            \end{cases}
            \label{Thm.5-141}
        \end{align}
\end{mTh}
\begin{rem}\label{Rem.mTh04}
    When $ \varepsilon = 0 $, the inclusion \eqref{Thm.5-132} is equivalent to:
    \begin{equation*}
        \sigma^\circ \in \mathrm{Sgn}^N(\nabla \theta^\circ), \mbox{ a.e. in $ Q $.}
    \end{equation*}
    In the meantime, when $ \varepsilon > 0 $, \eqref{Thm.5-13}--\eqref{Thm.5-141} imply that the pair of functions $ [p^\circ, z^\circ] $ solves the following system:
    \begin{align*}
        & \begin{cases}
            ~ -\partial_t p^\circ -{\mit \Delta} p^\circ +\alpha''(\eta^\circ) f_\varepsilon(\nabla \theta^\circ) p^\circ +g'(\eta^\circ) p^\circ +\alpha'(\eta^\circ)\nabla f_\varepsilon(\nabla \theta^\circ) \cdot \nabla z^\circ
            \\
            \qquad = M_\eta (\eta^\circ -\eta_\mathrm{ad}),
            \\[2ex]
            ~ -\partial_t (\alpha_0 z^\circ) -\mathrm{div} \bigl( \alpha(\eta^\circ) \nabla^2 f_\varepsilon(\nabla \theta^\circ) \nabla z^\circ +\nu^2 \nabla z^\circ +p^\circ \alpha'(\eta^\circ) \nabla f_\varepsilon(\nabla \theta^\circ) \bigr)
            \\
            \qquad = M_\theta(\theta^\circ -\theta_\mathrm{ad}),
        \end{cases}
    \end{align*}
    in the sense of distribution on $ Q $. Note that the above system corresponds to the distributional form of the variational system \eqref{Thm.5-01}--\eqref{Thm.5-03}, as in Main Theorem \ref{mainTh03} (\hyperlink{III-A}{III-A}). 
\end{rem}
\begin{rem}\label{Rem.mTh0304}
    Moreover, in the light of \eqref{Thm.5-001}, \eqref{Thm.5-101}, and Remark \ref{Rem.proj02} (\hyperlink{Fact4l}Fact\,4), we will observe that:
    \begin{align*}
        & 
        \begin{cases}
            M_u u^*(t, x) = M_u \bigl[\mathrm{proj}_{K} (-p_\varepsilon^*)\bigr](t, x) =M_u \bigl(\kappa^0 \vee ( \kappa^1 \wedge (-p_\varepsilon^*)) \bigr)(t, x),
            \\
            M_u u^\circ(t, x) = M_u \bigl[\mathrm{proj}_{K} (-p^\circ)\bigr](t, x)  = M_u \bigl(\kappa^0 \vee ( \kappa^1 \wedge (-p^\circ)) \bigr)(t, x),
        \end{cases} 
        \\
        & \hspace{27ex} \mbox{for a.e. $ (t, x) \in Q $.}
    \end{align*}
\end{rem}

\section{Proof of Main Theorem \ref{mainTh01}} 

In this Section, we give the proof of the first Main Theorem \ref{mainTh01}. Before the proof, we refer to the reformulation method as in \cite{MR3888636}, and consider to reduce the state-system (S)$_{\varepsilon}$ to an evolution equation in the Hilbert space $[H]^{2}$.

Let us fix any $\varepsilon \geq 0$. Besides, for any $ R \geq 0 $, let us define a proper functional $\Phi_{\varepsilon}^R:[H]^{2} \longrightarrow [0, \infty]$, by setting: 
\begin{align}\label{Phi_eps}
    \Phi_{\varepsilon}^R: w & = [\eta, \theta] \in [H]^{2} \mapsto \Phi_{\varepsilon}^R(w) =  \Phi_{\varepsilon}^R(\eta, \theta)
    \nonumber
    \\
    &:= \left\{
        \begin{array}{ll}
            \multicolumn{2}{l}{\displaystyle\frac{1}{2}\int_{\Omega}|\nabla\eta|^{2}dx +\frac{R}{2} \int_\Omega |\eta|^2 \, dx +\frac{1}{2}\int_{\Omega}\left(\nu f_{\varepsilon}(\nabla \theta) + \frac{1}{\nu}\alpha(\eta) \right)^{2}\, dx,}
            \\[2ex]
            & \mbox{if $[\eta, \theta] \in V\times V_{0}$,}
            \\[2ex]
            \infty, & \mbox{otherwise.}
        \end{array}
    \right. 
\end{align}
Note that the assumptions (\hyperlink{A2l}{A2}) and (\hyperlink{A4l}{A4}) guarantee the lower semi-continuity and convexity of $\Phi_{\varepsilon}^R$ on $[H]^{2}$. 
\begin{rem}\label{Rem.mTh01-01}
    As consequences of standard variational methods, we easily check the following facts.
    \begin{description}
        \item[(\hypertarget{Fact7}{Fact\,7})]For the operator $ \partial_\eta \Phi_\varepsilon^R : [H]^2 \longrightarrow 2^{H} $, 
            \begin{equation*}
                D(\partial_\eta \Phi_\varepsilon^R) = \left\{ \begin{array}{l|l}
                    [\tilde{\eta}, \tilde{\theta}] \in D & \parbox{4cm}{$ \tilde{\eta} \in H^2(\Omega) $ subject to $ \nabla \eta \cdot n_\Gamma = 0 $ in $ H^{\frac{1}{2}}(\Gamma) $}
                \end{array} \right\},
            \end{equation*}
            independent of $ R \geq 0 $, and $\partial_\eta \Phi_\varepsilon^R $ is a single-valued operator such that: 
            \begin{equation*}
                \partial_\eta \Phi_\varepsilon^R(w) = \partial_\eta \Phi_\varepsilon(\eta, \theta)
                = -{\mit \Delta} \eta +R\eta +\alpha'(\eta) f_\varepsilon(\nabla \theta) +\frac{1}{\nu^2} \alpha(\eta) \alpha'(\eta)  \mbox{ in $ H $,}
            \end{equation*}
            for all $w = [\eta, \theta] \in D(\partial_\eta \Phi_\varepsilon^R)$, and $ R \geq 0 $.
        \item[(\hypertarget{Fact8}{Fact\,8})]$ \partial_\theta \Phi_\varepsilon^R : [H]^2 \longrightarrow 2^H $ is independent of $ R \geq 0 $, and  $\theta \in D(\partial_\theta \Phi_\varepsilon^R)$, and $ \theta^* \in \partial_\theta \Phi_\varepsilon^R(w) = \partial_\theta \Phi_\varepsilon^R(\eta, \theta) $, iff. $\theta \in V_0$, and
            \begin{equation*}
                \begin{array}{c}
                    \ds (\theta^*, \theta -\psi)_H \geq \int_\Omega \alpha(\eta)f_\varepsilon(\nabla \theta) \, dx -\int_\Omega \alpha(\eta) f_\varepsilon(\nabla \psi)\, dx  
                    \\[2ex]
                    +\nu^2 \bigl( \nabla \theta, \nabla(\theta -\psi) \bigr)_{[H]^N},  \mbox{ for all $ \psi \in V_0 $, and  $ R \geq 0 $.} 
                \end{array}
            \end{equation*}
    \end{description}
    In addition, let us define time-dependent operators $\mathcal{A}(t) \in \mathscr{L}([H]^{2})$, for $t \in [0, T]$, nonlinear operators $\mathcal{G}^R : [H]^{2} \longrightarrow \mathbb [H]^{2}$, for $ R \geq 0 $, by setting:
\begin{align}
    \mathcal{A}(t): {w} = [{\eta}, {\theta}] \in [H & ]^{2} \mapsto \mathcal{A}(t) {w} := [{\eta}, \, \alpha_0(t) {\theta}] \in [H]^{2},
    \label{AA_0}
    \\
    & \mbox{for }t \in [0, T],
    \nonumber
\end{align}
\begin{align}
    \mathcal{G}^R: w = [\eta, \theta] \in [H]^{2} \mapsto \mathcal{G}^R(w) := &  \bigl[ g(\eta) -R \eta -\nu^{-2}\alpha(\eta)\alpha'(\eta), ~ 0  \bigr] \in [H]^{2}, 
    \label{GG}
    \\
    & \mbox{for $ R \geq 0 $,}
    \nonumber
\end{align}
    respectively. Then, based on the above (\hyperlink{Fact7}{Fact\,7}) and (\hyperlink{Fact8}{Fact\,8}), it is verified that the state-system (S)$_\varepsilon$ is equivalent to the following Cauchy problem.
    \begin{equation*}
        \left\{ \parbox{11cm}{
            $ \mathcal{A}(t) w'(t) +\bigl[ \partial_\eta \Phi_\varepsilon^R \times \partial_\theta \Phi_\varepsilon^R \bigr](w(t)) +\mathcal{G}^R(w(t)) \ni \mathfrak{f}(t) $  in $ [H]^2 $, 
            \\[1ex]
            \hspace*{4ex}a.e. $ t \in (0, T) $,
            \\[1.5ex]
            $ w(0) = w_0 $ in $ [H]^2 $.
        } \right.
    \end{equation*}
    In the context, ``\,$'$\,'' is the time-derivative, and
    \begin{equation}\label{w0-f}
        \left\{ \hspace{-3ex} \parbox{11cm}{
            \vspace{-1ex}
            \begin{itemize}
                \item $ w_0 := [\eta_0, \theta_0] \in D $ is the initial data of $ w = [\eta, \theta] $, 
                \item $ \mathfrak{f} := [M_u u, M_v v] \in [\mathscr{H}]^2 $ is the forcing term of the Cauchy problem.  
            \vspace{-1ex}
            \end{itemize}
        }
        \right.
    \end{equation}
\end{rem}

Now, before the proof of Main Theorem \ref{mainTh01}, we prepare the following Key-Lemma and its Corollary.
\begin{keyLem}\label{Lem03-01}
    Let us assume (\hyperlink{A1l}{A1})--(\hyperlink{A4l}{A4}). Then, there exists a positive constant $ R_0 > 0 $ such that: 
\begin{align*}
    \partial \Phi_\varepsilon^{R_0} = \bigl[ \partial_\eta \Phi_\varepsilon^{R_0} \times \partial_\theta \Phi_\varepsilon^{R_0} \bigr] \mbox{ in $ [H]^2 \times [H]^2 $.}
\end{align*}
\end{keyLem}
\paragraph{Proof.}{
    We set:
    \begin{equation}\label{R0}
        R_0 := 1 +\frac{2}{\nu^{2}}|\alpha|^2_{L^\infty(\mathbb{R})},
    \end{equation}
    and prove this $ R_0 $ is the required constant. 

    In the light of \eqref{prodSubDif}, it is immediately verified that:
    \begin{equation*}
        \partial \Phi_\varepsilon^{R_0}  \subset \bigl[ \partial_\eta \Phi_\varepsilon^{R_0} \times \partial_\theta \Phi_\varepsilon^{R_0} \bigr] \mbox{ in $ [H]^2 \times [H]^2 $.}
    \end{equation*}
    Hence, having in mind the maximality of the monotone graph $ \partial \Phi_\varepsilon^{R_0} $ in $ [H]^2 \times [H]^2 $, we can reduce our task to show the monotonicity of $ \bigl[ \partial_\eta \Phi_\varepsilon^{R_0} \times \partial_\theta \Phi_\varepsilon^{R_0} \bigr] $ in $ [H]^2 \times [H]^2 $.  
Let us assume:
\begin{equation*}
    \begin{array}{c}
        \ds [w, w^*] \in \bigl[ \partial_\eta \Phi_\varepsilon^{R_0} \times \partial_\theta \Phi_\varepsilon^{R_0} \bigr] \mbox{ and } [\tilde{w}, \tilde{w}^*] \in \bigl[ \partial_\eta \Phi_\varepsilon^{R_0} \times \partial_\theta \Phi_\varepsilon^{R_0} \bigr] \mbox{ in $ [H]^2 \times [H]^2 $,}
        \\[1ex]
        \mbox{with } w = [\eta, \theta] \in [H]^2, ~ w^* = [\eta^*, \theta^*] \in [H]^2,
        \\[1ex]
        \tilde{w} = [\tilde{\eta}, \tilde{\theta}] \in [H]^2, \mbox{ and } \tilde{w}^* = [\tilde{\eta}^*, \tilde{\theta}^*] \in [H]^2, \mbox{ respectively.}
    \end{array}
\end{equation*}
Then, by using \eqref{exM00}, (\hyperlink{Fact7}{Fact\,7}), (\hyperlink{Fact8}{Fact\,8}), (\hyperlink{A4l}{A4}), and Young's inequality, we compute that:
\begin{subequations}\label{IA0-00}
\begin{align}
    & (w^* -\tilde{w}^*, w -\tilde{w})_{[H]^2} = I_1 + I_2 + I_3, 
\end{align}
with 
\begin{align}
I_1 := & |\nabla(\eta - \tilde{\eta})|_{[H]^N}^2 + R_0|\eta - \tilde{\eta}|_H^2 +\nu^2|\nabla (\theta - \tilde{\theta})|_{[H]^N}^2,
\end{align}
\begin{align}
I_2:= & (\alpha'(\eta)f_\varepsilon(\nabla \theta) - \alpha'(\tilde{\eta})f_\varepsilon(\nabla \tilde{\theta}), \eta-\tilde{\eta})_H \nonumber \\
& \quad +\frac{1}{\nu^2}(\alpha(\eta)\alpha'(\eta)-\alpha(\tilde \eta)\alpha'(\tilde \eta), \eta- \tilde \eta)_H\nonumber \\
      = & \int_{\Omega}f_\varepsilon(\nabla \theta)(\alpha'(\eta) - \alpha'(\tilde{\eta}))(\eta - \tilde{\eta})\, dx 
      \\
        & \quad + \int_{\Omega} \alpha'(\tilde{\eta}) (f_\varepsilon(\nabla \theta)- f_\varepsilon(\nabla \tilde{\theta}))(\eta - \tilde{\eta})\, dx \nonumber 
        \\
        & \ds \quad +\frac{1}{2 \nu^2} \int_{\Omega} \bigl( {\textstyle\frac{d}{d\eta}  [\alpha^2](\eta) -\frac{d}{d\eta} [\alpha^2](\tilde{\eta})}\bigr)(\eta -\tilde{\eta}) \, dx
        \nonumber
        \\
      \geq & - |\alpha'|_{L^\infty(\mathbb{R})}|\eta-\tilde{\eta}|_H|\nabla (\theta - \tilde{\theta})|_{[H]^N} \nonumber \\
       \geq & - \frac{|\alpha'|_{L^\infty(\mathbb{R})}^2}{\nu^2}|\eta-\tilde{\eta}|_H^2 - \frac{\nu^2}{4}|\nabla (\theta - \tilde{\theta})|_{[H]^N}^2,
\end{align}
and
\begin{align}
I_3 := & \int_{\Omega}(\alpha(\eta) - \alpha(\tilde{\eta}))(f_\varepsilon(\nabla \theta) - f_\varepsilon(\nabla\tilde{\theta}))\, dx \nonumber \\ 
       \geq &  - |\alpha'|_{L^\infty(\mathbb{R})}|\eta-\tilde{\eta}|_H|\nabla (\theta - \tilde{\theta})|_{[H]^N} \nonumber \\
       \geq &  - \frac{|\alpha'|_{L^\infty(\mathbb{R})}^2}{\nu^2}|\eta-\tilde{\eta}|_H^2 - \frac{\nu^2}{4}|\nabla (\theta - \tilde{\theta})|_{[H]^N}^2.
\end{align}
\end{subequations}
Due to \eqref{R0}, the inequalities in \eqref{IA0-00} lead to:
\begin{align}
 & (w^* -\tilde{w}^*, w -\tilde{w})_{[H]^2} \geq |\eta - \tilde{\eta}|_V^2 + \frac{\nu^2}{2}|\theta - \tilde{\theta}|_{V_0}^2 \geq 0,
\end{align}
which implies the (strict) monotonicity of the operator $\bigl[ \partial_\eta \Phi_\varepsilon^{R_0} \times \partial_\theta \Phi_\varepsilon^{R_0} \bigr]$ in $ [H]^2 \times [H]^2 $. 
\qed
}
\begin{cor}\label{03-01}
    Under the notations and assumptions as in the previous Key-Lemma \ref{Lem03-01}, it further holds that
\begin{align*}
    \partial \Phi_\varepsilon^{R} = \bigl[ \partial_{\eta} \Phi_\varepsilon^{R} \times \partial_\theta \Phi_\varepsilon^{R} \bigr] \mbox{ in $ [H]^2 \times [H]^2 $, for any $R \geq 0$.}
\end{align*}
\end{cor}
\paragraph{Proof.}{
Let us take arbitrary two constants $0 \leq R, \tilde{R} < \infty$.
Then, from (\hyperlink{Fact7}{Fact\,7}), \cite[Theorem 2.10]{MR2582280}, and \cite[Corollary 2.11]{MR0348562}, we immediately have 
\begin{subequations}\label{sharp1}
\begin{align}
 D(\partial_{\eta} \Phi_{\varepsilon}^{R}) = D(\partial_{\eta} \Phi_{\varepsilon}^{\tilde{R}}) \mbox{ in } V,
\end{align}
    and for any $ w = [\eta, \theta] \in D(\partial_{\eta} \Phi_{\varepsilon}^{R}) = D(\partial_{\eta} \Phi_{\varepsilon}^{\tilde{R}}) $, 
\begin{align}
    \partial_{\eta} \Phi_{\varepsilon}^{R}(w) ~& = 
                    -{\mit \Delta}\eta + \tilde{R}\eta + (R - \tilde{R})\eta+ \alpha'(\eta)f_{\varepsilon}(\nabla \theta) + \nu^{-2}\alpha(\eta)\alpha'(\eta) \nonumber 
                    \\
                 & = \partial_{\eta} \Phi_{\varepsilon}^{\tilde{R}}(w) + (R -\tilde{R})\eta \mbox{ in } H.
\end{align}
\end{subequations}
Also, as a straightforward consequence of (\hyperlink{Fact8}{Fact\,8}), it is seen that:
\begin{align}\label{sharp2}
\partial_{\theta} \Phi_{\varepsilon}^{R} = \partial_{\theta} \Phi_{\varepsilon}^{\tilde{R}} \mbox{ in } H \times H.
\end{align}
In the meantime, invoking \eqref{Phi_eps}, \cite[Theorem 2.10]{MR2582280}, and \cite[Corollary 2.11]{MR0348562}, we will infer that
\begin{subequations}\label{sharp3}
\begin{align}
D(\partial \Phi_{\varepsilon}^{R}) = D(\partial \Phi_{\varepsilon}^{\tilde{R}}) \mbox{ in } D,
\end{align}
and
\begin{align}
\partial \Phi_{\varepsilon}^{R}(w) = \partial \Phi_{\varepsilon}^{\tilde{R}}(w) + (R - \tilde{R})[\eta, 0] \mbox{ in } [H]^2.
\end{align}
\end{subequations}

Now, let us take the constant $R_0 > 0$ obtained in Key-Lemma \ref{Lem03-01}. 
Then, owing to \eqref{sharp1}--\eqref{sharp3}, and Key-Lemma \ref{Lem03-01}, we can compute that
\begin{align}\label{sharp4}
    \bigl[\partial_{\eta} \Phi_{\varepsilon}^{R} & \times \partial_{\theta} \Phi_{\varepsilon}^{R} \bigr](w) = \bigl[\partial_{\eta} \Phi_{\varepsilon}^{R_0} \times \partial_{\theta} \Phi_{\varepsilon}^{R_0} \bigr](w) + (R - R_0)[\eta, 0]
\nonumber
\\
    & \hspace{-1ex} = \partial \Phi_{\varepsilon}^{R_0}(w) + (R -R_0)[\eta, 0] = \partial \Phi_{\varepsilon}^{R}(w) \mbox{ in } [H]^2,
 \\
    & \hspace{-2ex}  \mbox{ for any } w \in D(\partial_{\eta} \Phi_{\varepsilon}^{R} \times \partial_{\theta} \Phi_{\varepsilon}^{R}) = D(\partial_{\eta} \Phi_{\varepsilon}^{R}) \cap D(\partial_{\theta} \Phi_{\varepsilon}^{R}).\nonumber
\end{align}
In the light of \eqref{prodSubDif}, the above \eqref{sharp4} is sufficient to conclude this Corollary.
\qed
\medskip

}

\begin{lem}\label{eta.bdd}
    Let us assume (\hyperlink{A1l}{A1})--(\hyperlink{A4l}{A4}), and fix functions $\bar{\theta} \in L^\infty(0, T; V_0)$, $ \eta_0 \in V $, and $u \in \mathscr{H} $. Then, the initial-boundary value problem:
\begin{align}\label{eta.equ}
    & \begin{cases}
        \partial_t \eta - {\mit \Delta} \eta + g(\eta) + \alpha'(\eta)f_\varepsilon(\nabla \bar{\theta}) = M_u u\ \mbox{a.e. in $ Q $,}
        \\[0.5ex]
        \nabla \eta \cdot n_\Gamma = 0 \mbox{ on $ \Sigma $,}
        \\[0.5ex]
        \eta(0, x) = \eta_0(x), ~ x \in \Omega; 
    \end{cases}
\end{align}
    admits a unique solution $\eta \in W^{1,2}(0, T; H) \cap $ $ L^\infty(0, T; V) \cap L^2(0, T; H^2(\Omega))$, and in particular, if:
    \begin{equation}\label{etaEx00}
        \eta_0 \in L^\infty(\Omega), \mbox{ and } u \in L^\infty(Q),
    \end{equation}
    then it holds that $\eta \in L^\infty(Q)$.
\end{lem}
\paragraph{Proof. }
Let us fix $ \bar{\theta} \in L^\infty(0, T; V_0) $, $ \eta_0 \in V $, and $ u \in \sH $. Then, referring to the general theories of nonlinear evolution equations (e.g. \cite{MR0348562,MR2582280,Kenmochi81}), we immediately find a solution $ \eta \in W^{1, 2}(0, T; H) \cap L^\infty(0, T; V) \cap L^2(0, T; H^2(\Omega)) $, in the variational sense:
\begin{align}
    (\partial_t \eta(t), \varphi)_H ~& +(\nabla \eta(t), \nabla \varphi)_{[H]^N} +\bigl( g(\eta(t)) +\alpha'(\eta(t)) f_\varepsilon(\nabla \bar{\theta}(t)), \varphi \bigr)_H 
    \nonumber
    \\
    & = (M_u u(t), \varphi)_H, \mbox{ for any $ \varphi \in V $, a.e. $ t \in (0, T) $.}
    \label{etaEx01}
\end{align}

Next, we assume $ \eta_0 \in V \cap L^\infty(\Omega) $ and $ u \in L^\infty(Q) $, and verify the $ L^\infty $-regularity of the solution $ \eta $ as in \eqref{etaEx00}. To this end, we invoke the assumption (\hyperlink{A3l}{A3}), and take a large constant $ L_0 > 0 $, such that:
\begin{equation}\label{eta.bdd01}
    L_0 \geq |\eta_0|_{L^\infty(\Omega)}, ~ g(L_0) \geq M_u |u|_{L^\infty(Q)}, \mbox{ and } g(-L_0) \leq -M_u |u|_{L^\infty(Q)}.
\end{equation}
On this basis, we set our remaining task to show that:
\begin{equation}\label{eta.bdd02}
    |\eta|_{L^\infty(Q)} \leq L_0, \mbox{ i.e. } -L_0 \leq \eta \leq L_0 \mbox{ a.e. in $ Q $.}
\end{equation}

Due to \eqref{eta.bdd01} and (\hyperlink{A4l}{A4}), the constants $ L_0 $ and $ -L_0 $ fulfill that:
\begin{subequations}\label{eta.bdd03}
    \begin{equation}\label{eta.bdd03a}
        \partial_t L_0 -{\mit \Delta} L_0 +g(L_0) +\alpha'(L_0) f_\varepsilon(\nabla \bar{\theta}) \geq M_u u(t, x), \mbox{ a.e. $ (t, x) \in Q $,}
    \end{equation}
    and
    \begin{equation}\label{eta.bdd03b}
        \partial_t (-L_0) -{\mit \Delta} (-L_0) +g(-L_0) +\alpha'(-L_0) f_\varepsilon(\nabla \bar{\theta}) \leq M_u u(t, x), \mbox{ a.e. $ (t, x) \in Q $,}
    \end{equation}
\end{subequations}
respectively, together with the initial values $ L_0 $ and $ -L_0 $, and the zero-Neumann boundary conditions. 
\medskip

Now, let us take the difference between PDEs in \eqref{eta.equ} and \eqref{eta.bdd03a} (resp. \eqref{eta.bdd03b} and \eqref{eta.equ}), and multiply the both sides by $ [\eta -L_0]^+ $ (resp. $ [-L_0 -\eta]^+ $). Then, from (\hyperlink{A2l}{A2})--(\hyperlink{A4l}{A4}), it is inferred that:
\begin{align*}
    \frac{1}{2} \frac{d}{dt} \bigl( \bigl| & [\eta -L_0]^+(t) \bigr|_H^2 +\bigl| [-L_0 -\eta]^+(t) \bigr|_H^2 \bigr)
    \\
    & \leq |g'|_{L^\infty(\R)} \bigl( \bigl| [\eta -L_0]^+(t) \bigr|_H^2 +\bigl| [-L_0 -\eta]^+(t) \bigr|_H^2 \bigr), 
    \mbox{ a.e. $ t \in (0, T) $.} 
\end{align*}
Applying Gronwall's lemma, and invoking \eqref{eta.bdd01}, we obtain:
\begin{equation*}
    \bigl| [\eta -L_0]^+(t) \bigr|_H^2 +\bigl| [-L_0 -\eta]^+(t) \bigr|_H^2 \leq 0, \mbox{ a.e. $ t \in (0, T) $,}
\end{equation*}
which implies the validity of \eqref{eta.bdd02}. \hfill \qed
\medskip

\begin{rem}\label{Rem.sols}
Let $\varepsilon \geq 0$ be arbitrary constant.
    Then, as a consequence of (\hyperlink{Fact7}{Fact\,7}), (\hyperlink{Fact8}{Fact\,8}), Key-Lemma \ref{Lem03-01}, Corollary \ref{03-01}, and Lemma \ref{eta.bdd}, we can say that the state-system (S)$_\varepsilon$ is equivalent to the following Cauchy problem of evolution equation, denoted by (\hyperlink{E}{E})$_\varepsilon$.
\begin{description}
    \item[\textmd{(\hypertarget{E}{E})$ _\varepsilon $}]:
\vspace{-4ex}
\end{description}
\begin{align*}
      &  \left\{ \parbox{13cm}{
            $ \mathcal{A}(t) w'(t) + \partial \Phi_\varepsilon^R (w(t)) +\mathcal{G}^R(w(t)) \ni \mathfrak{f}(t) $ in $ [H]^2 $, a.e. $ t \in (0, T) $,
            \\[1.5ex]
            $ w(0) = w_0 $ in $ [H]^2 $,
        } \right.
    \end{align*}
for any $ R \geq 0 $. 
\end{rem}
\medskip

Now, we are ready to prove the Main Theorem \ref{mainTh01}.

\paragraph{Proof of Main Theorem \ref{mainTh01} (I-A).}{
     Let us fix any $R > 0$. Then, under the setting \eqref{Phi_eps}--\eqref{w0-f}, we immediately check that:
    \begin{description}
        \item[\textmd{(\hypertarget{ev0}{ev.0})}]for any $ t \in [0, T] $, $ \mathcal{A}(t) \in \mathscr{L}([H]^2) $ is positive and selfadjoint, and  
            \begin{equation*}
                (\mathcal{A}(t) w, w)_{[H]^2} \geq \delta_* |w|_{[H]^2}^2, \mbox{ for any $ w \in [H]^2 $,}
            \end{equation*}
            with the constant $ \delta_* \in (0, 1) $ as in (\hyperlink{A4l}{A4});
        \item[\textmd{(\hypertarget{ev1}{ev.1})}]$ \mathcal{A} \in W^{1, \infty}(0, T; \mathscr{L}([H]^2)) $, and  
            \begin{equation*}
                A^* := \mathrm{ess} \sup_{\hspace{-3ex}t \in (0, T)} \left\{ \max \{ |\mathcal{A}(t)|_{\mathscr{L}([H]^2)}, |\mathcal{A}'(t)|_{\mathscr{L}([H]^2)} \} \right\} \leq 1 +|\alpha_0|_{W^{1, \infty}(Q)} < \infty;
            \end{equation*}
        \item[\textmd{(\hypertarget{ev2}{ev.2})}]$ \mathcal{G}^R : [H]^2 \longrightarrow [H]^2 $ is a Lipschitz continuous operator with a Lipschitz constant:
            \begin{equation*}
                \mathrm{Lip}(\mathcal{G}) := R +|g'|_{L^\infty(\mathbb{R})} +\nu^{-2} {\textstyle{\bigl| \frac{d}{d\eta}(\alpha \alpha') \bigr|_{L^\infty(\mathbb{R})}},}
            \end{equation*}
            and $ \mathcal{G}^R $ has a $ C^1 $-potential functional
            \begin{align*}
                \widehat{\mathcal{G}}^R: w ~& = [\eta, \theta] \in [H]^2 \mapsto \widehat{\mathcal{G}}^R(w) := \int_\Omega \left( G(\eta) - \frac{R\eta^2}{2} -\frac{\alpha(\eta)^2}{2 \nu^2} \right) \, dx \in \mathbb{R};
            \end{align*}
        \item[\textmd{(\hypertarget{ev3}{ev.3})}]$ \Phi_\varepsilon^R \geq 0 $ on $ [H]^2 $, and the sublevel set $ \bigl\{  \tilde{w} \in [H]^2 \, \bigl| \, \Phi_\varepsilon^R(\tilde{w}) \leq r \bigr\} $ is contained in a compact set $ K_\nu^R(r) $ in $ [H]^2 $, defined as
            \begin{equation*}
               K_\nu^R(r) := \left\{ \begin{array}{l|l}
                   \tilde{w} = [\tilde{\eta}, \tilde{\theta}] \in D &  |\tilde{\eta}|_V^2 + |\tilde{\theta}|_{V_0}^2 \leq \frac{2r}{1 \wedge R \wedge \nu^2}
            \end{array} \right\}, 
        \end{equation*}
        for any $ r \geq 0 $.
    \end{description}
        On account of \eqref{Phi_eps}--\eqref{w0-f} and (\hyperlink{ev0}{ev.0})--(\hyperlink{ev3}{ev.3}), we can apply Proposition \ref{Lem.CP}, as the case when:
    \begin{align*}
        & X = [H]^2,  ~ \mathcal{A}_0 =  \mathcal{A} \mbox{ in $ W^{1, \infty}(0, T; \mathscr{L}([H]^2)) $, }\\
         \mathcal{G}_0 = & \mathcal{G}^R \mbox{ on $ [H]^2 $,  } \Psi_0 = \Phi_\varepsilon^R \mbox{ on $ [H]^2 $, and } \mathfrak{f}_0 = \mathfrak{f} \mbox{ in } [\mathscr{H}]^2,
    \end{align*}
        and we can find a solution $ w = [\eta, \theta] \in [\mathscr{H}]^2 $ to the Cauchy problem (\hyperlink{E}{E})$_\varepsilon$. In the light of  Proposition \ref{Lem.CP} and Remark \ref{Rem.sols}, finding this $ w = [\eta, \theta]$ directly leads to the existence and uniqueness of solution to the state-system (S)$_\varepsilon$. 
\medskip

Moreover, if $\eta_0 \in L^\infty(\Omega)$ and $u \in L^\infty(Q)$,  then the regularity $ \eta \in L^\infty(Q) $ will be immediately seen from Lemma \ref{eta.bdd}.
\hfill    \qed
}

\paragraph{\textbf{Proof of Main Theorem \ref{mainTh01} (I-B).}}{

Under the assumptions and notations as in  Main Theorem \ref{mainTh01}, we first fix a constant $R > 0$, and invoke Remark \ref{Rem.sols} to confirm that the solution $w:= [\eta, \theta] \in [\mathscr{H}]^2$ to (S)$_\varepsilon$ coincides with the solution to the Cauchy problem (\hyperlink{E}{E})$_\varepsilon$, and as well as, the solutions $w_n:= [\eta_n, \theta_n] \in [\mathscr{H}]^2$  to (S)$_{\varepsilon_n}$, $n = 1, 2, 3, \ldots ,$ coincide with the solutions to the Cauchy problems (\hyperlink{E}{E})$_{\varepsilon_n}$ for the initial data $ w_{0, n}:= [\eta_{0, n}, \theta_{0, n}] \in D $, and forcing terms $ \mathfrak{f_n} = [M_u u_n, M_v v_n] \in [\mathscr{H}]^2, n = 1, 2, 3,\ldots $, respectively.

On this basis, we next verify:   
\begin{description}
    \item[\textmd{(\hypertarget{ev4}{ev.4})}]$ \Phi_{\varepsilon_n}^R \geq 0 $ on $ [H] $, for $ n = 1, 2, 3, \dots $, and the union $ \bigcup_{n = 1}^\infty \bigl\{  \tilde{w} \in [H]^2 \, \bigl| \, \Phi_{\varepsilon_n}^R(\tilde{w}) \leq r \bigr\} $ of sublevel sets is contained in the compact set $ K_\nu^R(r) \subset [H]^2 $, as in (\hyperlink{ev3}{ev.3}), for any $ r > 0 $;
    \item[\textmd{(\hypertarget{ev5}{ev.5})}]$ \Phi_{\varepsilon_n}^R \to \Phi_\varepsilon^R $ on $ [H]^2 $, in the sense of Mosco, as $ n \to \infty $, more precisely, the uniform estimate \eqref{exM00} 
        will lead to the corresponding lower bound condition and optimality condition, in the Mosco-convergence of $ \{ \Phi_{\varepsilon_n}^R \}_{n=1}^{\infty} $;
    \item[\textmd{(\hypertarget{ev6}{ev.6})}]$ \sup_{n \in \mathbb{N}} \Phi_{\varepsilon_n}^R(w_{0, n}) < \infty $, and $ w_{0, n} \to w_{0} $ in $ [H]^2 $, as $ n \to \infty $, more precisely, it follows from \eqref{w.i01}, (\hyperlink{A1l}{A1}), and (\hyperlink{A4l}{A4}) that 
        \begin{align*}
            \sup_{n \in \mathbb{N}} \Phi_{\varepsilon_n}^R(w_{0, n}) & \leq \sup_{n \in \mathbb{N}} \left( {\textstyle \frac{1+R}{2}|\eta_{0, n}|_V^2 +\nu^2 ( \mathcal{L}^N(\Omega) +|\theta_{0, n}|_{V_0}^2) +\frac{1}{\nu^2} |\alpha(\eta_{0, n})|_H^2} \right) < \infty,
        \end{align*}
        and the weak convergence of $\{w_{0, n} \}_{n=1}^{\infty}$ in $D = V \times V_0$ and the compactness of embedding $D \subset [H]^2$ imply the strong convergence of $\{w_{0, n} \}_{n=1}^{\infty}$ in $ [H]^2$.
\end{description}

On account of \eqref{w.i01} and (\hyperlink{ev0}{ev.0})--(\hyperlink{ev6}{ev.6}), we can apply Proposition \ref{Lem.CP02}, to show that:
\begin{subequations}\label{convKS}
\begin{equation}\label{convKS01}
    \begin{cases}
        w_n \to w \mbox{ in $ C([0, T]; [H]^2) $} \mbox{ (i.e. in $ [C([0, T]; H)]^2 $)}, 
        \\
        \quad \mbox{ weakly in $ W^{1, 2}(0, T; [H]^2) $} \mbox{ (i.e. weakly in $ [W^{1, 2}(0, T; H)]^2  $),}
        \\[1ex]
        \displaystyle \int_0^T \Phi_{\varepsilon_n}^R(w_n(t)) \, dt \to \int_0^T \Phi_{\varepsilon}^R(w(t)) \, dt,
    \end{cases}
    \mbox{as $ n \to \infty $,}
\end{equation}
\begin{align*}
    \displaystyle \sup_{n \in \mathbb{N}} & |w_n|_{L^\infty(0, T; V) \times L^\infty(0, T; V_0)}^2 \leq 4 \sup_{n \in \mathbb{N}} |w_n|_{L^\infty(0, T; V \times V_0)}^2 
    \\
    & \displaystyle \leq \frac{8}{1 \wedge \nu^2 \wedge R} \sup_{n \in \mathbb{N}} |\Phi_{\varepsilon_n}^R(w_n)|_{L^\infty(0, T)} < \infty, 
\end{align*}
and hence,
\begin{equation}\label{convKS02}
    w_n \to w \mbox{ weakly-$*$ in $ L^\infty(0, T; V) \times L^\infty(0, T; V_0) $, as $ n \to \infty $.}
\end{equation}
\end{subequations}
Furthermore, from \eqref{f_eps}, \eqref{exM00}, \eqref{convKS}, and the assumptions (\hyperlink{A2l}{A2}) and (\hyperlink{A4l}{A4}), one can observe that:
\begin{subequations}\label{convKSs}
\begin{equation}\label{convKS03}
    \begin{cases}
        \displaystyle \varliminf_{n \to \infty} \frac{1}{2} | \nabla \eta_n |_{[\mathscr{H}]^N}^2 \geq \frac{1}{2} |\nabla \eta|_{[\mathscr{H}]^N}^2, \quad \varliminf_{n \to \infty} \frac{R}{2} | \eta_n |_{\mathscr{H}}^2 \geq \frac{R}{2} |\eta|_{\mathscr{H}}^2, 
        \\[2ex]
        \displaystyle \varliminf_{n \to \infty} \frac{\nu^2}{2} |\theta_n |_{\mathscr{V}_0}^2 \geq \frac{\nu^2}{2} |\theta|_{\mathscr{V}_0}^2, \quad \lim_{n \to \infty} \frac{1}{2 \nu^2} |\alpha(\eta_n)|_{\mathscr{H}}^2 = \frac{1}{2 \nu^2} |\alpha(\eta)|_{\mathscr{H}}^2,
    \end{cases}
\end{equation}
and
\begin{align}\label{convKS04}
    \ds \varliminf_{n \to \infty} & \ds \bigl| \alpha(\eta_n) f_{\varepsilon_n}(\nabla \theta_n) \big|_{L^1(Q)} = \varliminf_{n \to \infty} \int_0^T \int_\Omega \alpha(\eta_n(t)) f_{\varepsilon_n}(\nabla \theta_n(t)) \, dx dt
    \nonumber
    \\
    & \geq \varliminf_{n \to \infty} \int_0^T \int_\Omega \alpha(\eta(t)) f_{\varepsilon_n}(\nabla \theta_n(t)) \, dx dt
    \nonumber
    \\
    & \qquad -\lim_{n \to \infty} |\alpha(\eta_n) -\alpha(\eta)|_{\mathscr{H}} \cdot \sup_{n \in \N} \bigl( \sqrt{\mathcal{L}^{N +1}(Q)} \, \varepsilon_n +|\theta_n|_{\mathscr{V}_0} \bigr)
    \nonumber
    \\
    & \geq \varliminf_{n \to \infty} \int_0^T \int_\Omega \alpha(\eta(t)) f_{\varepsilon}(\nabla \theta_n(t)) \, dx dt -|\alpha(\eta)|_{L^1(Q)} \cdot \lim_{n \to \infty} |\varepsilon_n -\varepsilon|
    \nonumber
    \\
    & \geq \int_0^T \int_\Omega \alpha(\eta(t)) f_{\varepsilon}(\nabla \theta(t)) \, dx dt = \bigl| \alpha(\eta) f_{\varepsilon}(\nabla \theta) \big|_{L^1(Q)}.
\end{align}
\end{subequations}
Here, from \eqref{Phi_eps}, it is seen that:
\begin{align}\label{convKS'}
    & \displaystyle \int_0^T {\Phi}_{\tilde{\varepsilon}}^R(\tilde{w}(t)) \, dt =  \int_0^T {\Phi}_{\tilde{\varepsilon}}^R(\tilde{\eta}(t), \tilde{\theta}(t)) \, dt
    \nonumber
    \\
    & \quad = \frac{1}{2} |\nabla\tilde{\eta}|_{[\mathscr{H}]^N}^2 + \frac{R}{2}|\tilde{\eta} |_{\mathscr{H}}^2 +\frac{\nu^2}{2} |\tilde{\theta}|_{\mathscr{V}_0}^2 +\bigl| \alpha(\tilde{\eta}) f_{\tilde{\varepsilon}}(\nabla \tilde{\theta}) \bigr|_{L^1(Q)} +\frac{1}{2 \nu^2} |\alpha(\tilde{\eta})|_{\mathscr{H}}^2 
    \nonumber
    \\[1ex]
    & \qquad +\frac{\nu^2 \tilde{\varepsilon}^2}{2} \mathcal{L}^{N +1}(Q), \mbox{ for all $ \tilde{\varepsilon} \geq 0 $ and $ \tilde{w} = [\tilde{\eta}, \tilde{\theta} ] \in  \mathscr{Y}$.}
\end{align}
Taking into account \eqref{convKS01}, \eqref{convKSs}, and \eqref{convKS'}, we deduce that:
\begin{align}\label{convKS05}
        |\nabla\eta_n|_{[\mathscr{H}]^N}^2 & \, + R|\eta_n|_\mathscr{H}^2 +\nu^2 |\theta_n|_{\mathscr{V}_0}^2 \to |\nabla\eta|_{[\mathscr{H}]^N}^2 + R|\eta|_\mathscr{H}^2 + \nu^2 |\theta|_{\mathscr{V}_0}^2, \nonumber
        \\
        & \, \mbox{and hence, } |[\eta_n, \theta_n]|_{\mathscr{Y}} \to |[\eta, \theta]|_{\mathscr{Y}}, \mbox{ as $ n \to \infty $.}
\end{align}

Since the norm of Hilbert space $ \mathscr{Y} := \mathscr{V} \times \mathscr{V}_0 $ is uniformly convex, the convergences \eqref{convKS02} and \eqref{convKS05} imply the strong convergence:
\begin{subequations}\label{convKS06-07}
\begin{equation}\label{convKS06}
    w_n \to w \mbox{ in $ \mathscr{Y} $, as $ n \to \infty $,}
\end{equation}
and furthermore, it follows from \eqref{exM00} and \eqref{convKS06} that:
\begin{align}\label{convKS07}
    |f_{\varepsilon_n} (\nabla \theta_n) & -f_\varepsilon(\nabla \theta)|_\mathscr{H} \leq  |f_{\varepsilon_n}(\nabla \theta_n) -f_{\varepsilon_n}(\nabla \theta)|_\mathscr{H} +|f_{\varepsilon_n}(\nabla \theta) -f_{\varepsilon}(\nabla \theta)|_\mathscr{H} 
    \nonumber
    \\
    \leq &  |\theta_n -\theta|_{\mathscr{V}_0} +\sqrt{\mathcal{L}^{N +1}(Q)} |\varepsilon_n -\varepsilon| \to 0, \mbox{ as $ n \to \infty $.}
\end{align}
\end{subequations}
The convergences \eqref{convKS} and \eqref{convKS06-07} 
are sufficient to obtain the convergence \eqref{mThConv} as in Main Theorem \ref{mainTh01} (\hyperlink{I-B}{I-B}).

Finally, let us assume \eqref{w.i03} to verify \eqref{mThConv00}. In the light of (\hyperlink{A3l}{A3}), we can take a large constant $ L_* > 0 $, independent of $ n $, such that:
\begin{align}\label{etan.bdd01}
    L_* \geq \sup_{n \in \N} &\, |\eta_{0, n}|_{L^\infty(\Omega)}, ~ g(L_*) \geq M_u \sup_{n \in \N} |u_n|_{L^\infty(Q)}, 
    \nonumber
    \\
    & \mbox{and } g(-L_*) \leq - M_u \sup_{n \in \N} |u_n|_{L^\infty(Q)}.
\end{align}
Then, just as in the derivation of \eqref{eta.bdd02}, we can show that:
\begin{equation}\label{etan.bdd01}
    \sup_{n\in \N} |\eta_n|_{L^\infty(Q)} \leq L_*, \mbox{ i.e. } -L_* \leq \eta_n \leq L_* \mbox{ a.e. in $ Q $, $ n = 1, 2, 3, \dots $.} 
\end{equation}
The convergence \eqref{mThConv}, and the $ L^\infty $-weak-$*$ compactness brought by \eqref{etan.bdd01} lead to the convergence \eqref{mThConv00}.  \hfill \qed
}

\section{Proof of Main Theorem \ref{mainTh02}}

In this section, we prove the second Main Theorem \ref{mainTh02}. 
Before the proof, we prepare the following lemma. 

\begin{lem}\label{Lem.MoscoK}
    Let us assume (\hyperlink{A5l}{A5}) and (\hyperlink{A6l}{A6}), and let us fix the function $ \bar{\kappa} \in \bigcap_{n = 1}^\infty K_n $ as in   (\hyperlink{A6l}{A6}). Besides, let us take any function $ u \in K $, and define a sequence $ \{ u_n \}_{n = 1}^\infty \subset \sH $, by setting:
    \begin{equation*}
        u_n := \mathrm{proj}_{K_n} (u) = \kappa_n^0 \vee (\kappa_n^1 \wedge u) \in K_n \mbox{ in $ \sH $, for $ n = 1, 2, 3, \dots $.}
    \end{equation*}
    Then, it holds that:
    \begin{equation}\label{Mk00}
        u_n \to u \mbox{ in $ \sH $ as $ n \to \infty $.}
    \end{equation}
\end{lem}
\paragraph{Proof.}{
As is easily seen,
\begin{align}\label{un.bar}
       u_n(t, x) &= \left\{ \begin{array}{ll}
            \kappa_n^1(t, x), & \mbox{if $ u(t, x) > \kappa_n^1(t, x) $,} 
            \\[1ex]
            u(t, x), & \mbox{if $ \kappa_n^0(t, x) \leq u(t, x) \leq \kappa_n^1(t, x) $,}
            \\[1ex]
            \kappa^0_n(t, x), & \mbox{if $ u(t, x) < \kappa_n^0(t, x) $,}
        \end{array}
        \right.
        \ \ \\ \nonumber
        \\[-2.5ex]
        & \qquad \mbox{a.e. $(t, x) \in Q$, $ n = 1, 2, 3, \dots $,} 
        \nonumber
\end{align}
so that:
\begin{equation}\label{Mk01}
    |u_n -u| \to 0, \mbox{ in the pointwise sense, a.e. in $ Q $, as $ n \to \infty $.}
\end{equation}
Also, owing to the presence of $ \bar{\kappa} \in \bigcap_{n = 1}^\infty K_n $, as in (\hyperlink{A6l}{A6}),
\begin{align*}
& -[u -\bar{\kappa}]^- \leq u_n -\bar{\kappa} \leq [u -\bar{\kappa}]^+, \mbox{ a.e. in $ Q $,}
    \\
    \mbox{i.e. }  | u & _n -\bar{\kappa}| \leq [u -\bar{\kappa}]^+ +[u -\bar{\kappa}]^- = |u -\bar{\kappa}|, \mbox{ a.e. in $ Q $,} 
\end{align*}
which leads to:
\begin{align}\label{Mk02}
    |u_n -u| \leq |u_n - & \bar{\kappa}| +|u -\bar{\kappa}| \leq 2 |u -\bar{\kappa}|~  \mbox{ a.e. in $ Q $,}
    \nonumber
    \\
    & \mbox{with $ |u -\bar{\kappa}| \in \sH $.}
\end{align}
The convergence \eqref{Mk00} will be deduced as a straightforward consequence of \eqref{Mk01}, \eqref{Mk02}, and the dominated convergence theorem \cite[Theorem 10 on page 36]{MR0492147}. \qed
}
\bigskip

Now, let $[\eta_0, \theta_0] \in D$ be the initial pair, and any constraint $K = [\![ \kappa^0, \kappa^1]\!] \in \mathfrak{K}$. Also, let us fix arbitrary forcing pair $ [\bar{u}, \bar{v}] \in \mathscr{U}_\mathrm{ad}^K $, and let us invoke the definition of the cost function $ \mathcal{J}_\varepsilon $, defined in \eqref{J}, to estimate that:
\begin{align}\label{mTh02-00}
    0 &~ \leq \underline{J}_{\varepsilon}:= \inf \mathcal{J}_{\varepsilon}(\mathscr{U}_\mathrm{ad}^K) \leq \overline{J}_{\varepsilon}:= \mathcal{J}_\varepsilon(\bar{u}, \bar{v}) < \infty, \mbox{ for all } \varepsilon \geq 0.
\end{align}
Also, for any $ \varepsilon \geq 0 $, we denote by $[\bar{\eta}, \bar{\theta}]$ the solution to (S)$_{\varepsilon}$, for the initial pair $[\eta_0, \theta_0]$ and forcing pair $ [\bar{u}, \bar{v}]$.
\bigskip

Based on these, the Main Theorem \ref{mainTh02} is proved  as follows.

\paragraph{Proof of Main Theorem \ref{mainTh02} (II-A).}
Let us fix any $\varepsilon \geq 0$. 
Then, from the estimate \eqref{mTh02-00}, we immediately find a sequence of forcing pairs $\{[u_{n}, v_{n}] \}_{n=1}^{\infty} \subset \mathscr{U}_\mathrm{ad}^K$, such that:
\begin{subequations}\label{mTh02-0102}
\begin{equation}\label{mTh02-01}
\mathcal{J}_{\varepsilon}(u_{n}, v_{n}) \downarrow \underline{J}_{\varepsilon},\ \mbox{as}\ n \to \infty,
\end{equation}
and
\begin{align}\label{mTh02-02}
    \ds\frac{1}{2}\sup_{n \in \mathbb{N}}& \bigl| [\sqrt{M_u} u_{n}, \sqrt{M_v} v_{n}] \bigr|_{[\mathscr{H}]^{2}}^{2} \leq \mathcal{J}_\varepsilon(\bar{u}, \bar{v}) < \infty.
\end{align}
\end{subequations}
Also, the estimate \eqref{mTh02-02} and the assumption (\hyperlink{A5l}{A5}) enable us to take a subsequence of $\{[u_{n}, v_{n}] \}_{n=1}^{\infty} \subset \mathscr{U}_\mathrm{ad}^K$ (not relabeled), and to find a pair of functions $[u^{*}, v^{*}] \in \mathscr{U}_\mathrm{ad}^K$, such that:
\begin{equation}\label{mTh02-03}
\begin{array}{c}
[\sqrt{M_u} u_{n}, \sqrt{M_v} v_{n}] \to [\sqrt{M_u} u^{*}, \sqrt{M_v} v^{*}]\ \mbox{weakly in}\ [\mathscr{H}]^{2}, \mbox{as}\ n \to \infty, 
\end{array}
\end{equation}

Let $[\eta^{*}, \theta^{*}] \in [\mathscr{H}]^{2}$ be the solution to (S)$_{\varepsilon}$, for the initial pair $[\eta_0, \theta_0]$ and forcing pair $[u^{*}, v^{*}]$. 
Also, for any $n \in \mathbb{N}$, let $[\eta_{n}, \theta_{n}] \in [\mathscr{H}]^{2}$ be the solution to (S)$_{\varepsilon_n}$, for the forcing pair $[u_{n}, v_{n}]$. Then, having in mind \eqref{w.i01}, \eqref{mTh02-03}, and the initial condition:
\begin{equation*}
    [\eta_{n}(0), \theta_{n}(0)] = [\eta^{*}(0), \theta^{*}(0)] = [\eta_{0}, \theta_{0}] \mbox{ in $[H]^{2}$, for $n = 1, 2, 3, \dots $,}
\end{equation*}
we can apply Main Theorem \ref{mainTh01} (\hyperlink{I-B}{I-B}), to see that:
\begin{equation}\label{mTh02-04}
[\eta_{n}, \theta_{n}] \to [\eta^{*}, \theta^{*}]\ \mbox{in}\ [C([0, T]; H)]^{2},\ \mbox{as}\ n \to \infty.
\end{equation}
On account of \eqref{mTh02-01}, \eqref{mTh02-03}, and \eqref{mTh02-04}, it is computed that:
\begin{align*}
    \mathcal{J}_{\varepsilon}(u^{*}, v^{*}) & = \frac{1}{2} \bigl| [\sqrt{M_\eta} (\eta^{*}-\eta_\mathrm{ad}), \sqrt{M_\theta} (\theta^{*}-\theta_\mathrm{ad}) ] \bigr|_{[\mathscr{H}]^{2}}^{2}\\
                                            &\qquad + \frac{1}{2} \bigl| [\sqrt{M_u} u^{*}, \sqrt{M_v} v^{*}] \bigr|_{[\mathscr{H}]^{2}}^{2}
    \\
                                            & \leq \frac{1}{2}\lim_{n \to \infty} \bigl| [\sqrt{M_\eta} (\eta_{n}-\eta_\mathrm{ad}), \sqrt{M_\theta} (\theta_{n}-\theta_\mathrm{ad}) ] \bigr|_{[\mathscr{H}]^{2}}^{2}\\
                                            &\qquad + \frac{1}{2} \varliminf_{n \to \infty} \bigl| [\sqrt{M_u} u_{n}, \sqrt{M_v} v_{n}] \bigr|_{[\mathscr{H}]^{2}}^{2}
    \\
    & = \lim_{n \to \infty}\mathcal{J}_{\varepsilon}(u_{n}, v_{n}) = \underline{J}_{\varepsilon} ~ (\leq \mathcal{J}_{\varepsilon}(u^{*}, v^{*})),
\end{align*}
and this leads to:    
\begin{equation*}
    \mathcal{J}_\varepsilon(u^{*}, v^{*}) = \min_{[u, v] \in \mathscr{U}_\mathrm{ad}^K} \mathcal{J}_{\varepsilon}(u, v).
\end{equation*}

Thus, we conclude the item (\hyperlink{II-A}{II-A}). 
\qed

\paragraph{Proof of Main Theorem \ref{mainTh02} (II-B).}
{
    Let us take $ \varepsilon \geq 0 $, $ \{ \varepsilon_n \}_{n = 1}^\infty \subset [0, \infty) $, and $ \{ [\eta_{0, n}, \theta_{0, n}] \}_{n = 1}^\infty \subset D $ as in \eqref{w.i01}. Besides, for the pair of functions $ [\bar{u}, \bar{v}] \in \mathcal{U}_\mathrm{ad}^K $ as in \eqref{mTh02-00}, let us define:
    \begin{equation*}
        \bar{u}_n := \mathrm{proj}_{K_n}(\bar{u}) = \kappa_n^0 \vee (\kappa_n^1 \wedge \bar{u}) \in K_n, ~ n = 1, 2, 3, \dots. 
    \end{equation*}
    Then, from Lemma \ref{Lem.MoscoK}, it immediately follows that:
    \begin{equation}\label{mTh02-081}
        \bar{u}_n \to \bar{u} \mbox{ in $ \mathscr{H} $, as $ n \to \infty $.}
    \end{equation}
    Here, let $ [\bar{\eta}, \bar{\theta}] \in [\mathscr{H}]^2 $ be the solution to (S)$_\varepsilon$, for the initial pair $ [\eta_0, \theta_0] $ and forcing pair $ [\bar{u}, \bar{v}] $, and let $ [\bar{\eta}_n, \bar{\theta}_n] \in [\mathscr{H}]^2 $, $ n = 1, 2, 3, \dots $, be solutions to (S)$_{\varepsilon_n}$, for the initial pairs $ [\eta_{0, n}, \theta_{0, n}] $, and forcing pairs $ [\bar{u}_n, \bar{v}] $, $ n = 1, 2, 3, \dots $, respectively. Then, invoking \eqref{w.i01} and \eqref{mTh02-081}, we can apply Main Theorem \ref{mainTh01} (I-B) to these solutions, and we can see that:
    \begin{subequations}\label{mTh02-09}
        \begin{equation}\label{mTh02-09a}
            [\bar{\eta}_n, \bar{\theta}_n] \to [\bar{\eta}, \bar{\theta}] \mbox{ in $ [C([0, T]; H)]^2 $,}
        \end{equation}
        and in particular, 
        \begin{align}\label{mTh02-09b}
            [\eta_{0, n}, \theta_{0, n}] = [ & \bar{\eta}_n(0), \bar{\theta}_n(0)] \to [\eta_0, \theta_0] = [\bar{\eta}(0), \bar{\theta}(0)]
            \nonumber
            \\
            & \mbox{in $ [H]^2 $, as $ n \to \infty $.}
        \end{align}
    \end{subequations}
    The convergences \eqref{mTh02-081} and \eqref{mTh02-09} enable us to estimate:
    \begin{equation}\label{mTh02-10}
        \overline{J}_\mathrm{sup} := \sup_{n \in \N} \mathcal{J}_{\varepsilon_n}(\bar{u}_n, \bar{v}) < \infty.
    \end{equation}

    Next, for any $ n \in \N $, let us denote by $ [\eta_n^*, \theta_n^*] \in [\mathscr{H}]^2 $ the solution to (S)$_{\varepsilon_n}$, for the initial pair $ [\eta_{0, n}, \theta_{0, n}] $, and forcing pair $ [u_n^*, v_n^*] $ of the optimal control of (\hyperlink{OP}{OP})$_{\varepsilon_n}^{K_n}$. Then, in the light of \eqref{mTh02-00} and \eqref{mTh02-10}, it is observed that:
    \begin{equation*}
        0 \leq \frac{1}{2} \bigl| [ {\ts\sqrt{M_u} u_n^*, \sqrt{M_v} v_n^*} ] \bigr|_{[\mathscr{H}]^2}^2 \leq \underline{J}_{\varepsilon_n} \leq \overline{J}_\mathrm{sup} < \infty, ~ n = 1, 2, 3, \dots. 
    \end{equation*}
    Therefore, one can find a subsequence $ \{ n_i \}_{i = 1}^\infty \subset \{n\} $, together with a limiting pair of functions $ [u^{**}, v^{**}] \in [\mathscr{H}]^2 $, such that:
}
\begin{equation}\label{mTh02-11}
\begin{array}{c}
[\sqrt{M_u} u^{*}_{n_i}, \sqrt{M_v} v^{*}_{n_i}] \to [\sqrt{M_u} u^{**}, \sqrt{M_v} v^{**}]\ \mbox{weakly in}\ [\mathscr{H}]^{2},\ \mbox{as}\ i \to \infty, 
\\[0.5ex]
\mbox{and as well as}\ [M_u u^{*}_{n_i}, M_v v^{*}_{n_i}] \to [M_u u^{**}, M_v v^{**}]\ \mbox{weakly in}\ [\mathscr{H}]^{2}, \ \mbox{as}\ i \to \infty.
\end{array}
\end{equation}
    Additionally, for every $ \ell = 0, 1 $, the convex functionals on $ L^1(Q \setminus |\kappa^\ell|^{-1}(\infty)) $, defined as:
\begin{equation*}
    \tilde{u} \in L^1(Q \setminus |\kappa^\ell|^{-1}(\infty)) \mapsto \int_{Q \setminus |\kappa^\ell|^{-1}(\infty)} [\tilde{u}]^+ \, dx dt \in [0, \infty), ~ \ell = 0, 1,
\end{equation*}
are weakly lower semi-continuous. Therefore, we can observe from \eqref{mTh02-11} and (\hyperlink{A6l}{A6}) that:
\begin{subequations}\label{mTh02-100}
    \begin{center}
    \parbox{10cm}{
    \begin{equation*}
        \begin{cases}
            [\kappa^0 -u^{**}]^+ \leq [\bar{\kappa} -u^{**}]^+ \in \mathscr{H} \subset L^1(Q),
            \\[0.5ex]
            [u^{**} -\kappa^1]^+ \leq [u^{**} -\bar{\kappa}]^+ \in \mathscr{H} \subset L^1(Q),
        \end{cases}
    \end{equation*}
        }
    \end{center}
\begin{align}\label{mTh02-100a}
    \bigl| M_u [\kappa^0 & -u^{**}]^+ \bigr|_{L^1(Q)} = \int_{{Q \setminus |\kappa^0|^{-1}(\infty)}} \hspace{-0ex} [M_u(\kappa^0 -u^{**})]^+ \, dx dt
    \nonumber
    \\
    & \leq \varliminf_{i \to \infty} \int_{{Q \setminus |\kappa^0|^{-1}(\infty)}} \hspace{-0ex} [M_u(\kappa_{\varepsilon_{n_i}}^0 -u_{n_i}^{*})]^+ \, dx dt = 0,
\end{align}
and
\begin{align}\label{mTh02-100b}
    \bigl| M_u [u^{**} & -\kappa^1]^+ \bigr|_{L^1(Q)} = \int_{{Q \setminus |\kappa^1|^{-1}(\infty)}} \hspace{-0ex} [M_u(u^{**} -\kappa^1)]^+ \, dx dt
    \nonumber
    \\
    & \leq \varliminf_{i \to \infty} \int_{{Q \setminus |\kappa^1|^{-1}(\infty)}} \hspace{-0ex} [M_u(u_{n_i}^{*} -\kappa_{\varepsilon_{n_i}}^1)]^+ \, dx dt = 0.
\end{align}
\end{subequations}
Since the limit $ u^{**} $, when $ M_u = 0 $, can be taken arbitrary, the estimates as in \eqref{mTh02-100} enable us to suppose that:
\begin{equation*}
    \kappa^0 \leq u^{**} \leq \kappa^1 \mbox{ a.e. in $ Q $, i.e. } [u^{**}, v^{**}] \in \mathscr{U}_\mathrm{ad}^K.
\end{equation*}

Now, let us denote by $[\eta^{**}, \theta^{**}] \in [\mathscr{H}]^2 $ the solution to (S)$_{\varepsilon}$, for the initial pair $ [\eta_0, \theta_0] $ and forcing pair $[u^{**}, v^{**}]$.
Then, applying Main Theorem \ref{mainTh01} (\hyperlink{I-B}{I-B}), again, to the solutions $[\eta^{**}, \theta^{**}]$ and $[\eta^{*}_{n_i}, \theta^{*}_{n_i}]$, $i = 1, 2, 3, \dots$, one can see that:
\begin{align}\label{mTh02-12}
    [\eta_{n_i}^*, \theta_{n_i}^*] & \to [\eta^{**}, \theta^{**}] \mbox{ in $ [C([0, T]; H)]^2 $, in $ \mathscr{Y} $,}
    \nonumber
    \\
    & \mbox{weakly in $ [W^{1, 2}(0, T; H)]^2 $, and }
    \nonumber
    \\
    & \mbox{weakly-$*$ in $ L^\infty(0, T; V) \times L^\infty(0, T; V_0) $, as $ i \to \infty $.}
\end{align}
As a consequence of \eqref{mTh02-081}, \eqref{mTh02-09}, \eqref{mTh02-11}, and \eqref{mTh02-12}, it is verified that:
\begin{align*}
    \mathcal{J}_{\varepsilon} (u^{**}, v^{**}) &= \frac{1}{2} \bigl| [\sqrt{M_\eta} (\eta^{**}-\eta_\mathrm{ad}), \sqrt{M_\theta} (\theta^{**}-\theta_\mathrm{ad})] \bigr|_{[\mathscr{H}]^{2}}^{2}
    \\
    &\qquad + \frac{1}{2} \bigl| [\sqrt{M_u} u^{**}, \sqrt{M_v} v^{**}] \bigr|_{[\mathscr{H}]^{2}}^{2}
    \\
    & \leq \frac{1}{2} \lim_{i \to \infty} \bigl| [\sqrt{M_\eta} (\eta^{*}_{n_i} -\eta_\mathrm{ad}), \sqrt{M_\theta} (\theta^{*}_{n_i} -\theta_\mathrm{ad})] \bigr|_{[\mathscr{H}]^{2}}^{2}
    \\
    &\qquad  + \frac{1}{2}\varliminf_{i \to \infty} \bigl| [\sqrt{M_u} u^{*}_{n_i}, \sqrt{M_v} v^{*}_{n_i}] \bigr|_{[\mathscr{H}]^{2}}^{2}
    \\
    & = \varliminf_{i \to \infty}\mathcal{J}_{\varepsilon_{n_i}}(u^{*}_{n_i}, v^{*}_{n_i}) \leq \lim_{i \to \infty}\mathcal{J}_{\varepsilon_{n_i}}(\bar{u}_{n_i}, \bar{v})
    \\
    & = \frac{1}{2}\lim_{i \to \infty} \bigl| [\sqrt{M_\eta} (\bar{\eta}_{n_i}-\eta_\mathrm{ad}), \sqrt{M_\theta} (\bar{\theta}_{n_i}-\theta_\mathrm{ad})] \bigr|_{[\mathscr{H}]^{2}}^{2}
    \\
    &\qquad  + \frac{1}{2} \lim_{i \to \infty} \bigl| [\sqrt{M_u} \bar{u}_{n_i}, \sqrt{M_v} \bar{v}] \bigr|_{[\mathscr{H}]^{2}}^{2}
    \\
    & = \mathcal{J}_{\varepsilon}(\bar{u}, \bar{v}).
\end{align*}
Since the choice of $[\bar{u}, \bar{v}] \in \mathscr{U}_\mathrm{ad}^{K}$ is arbitrary, we conclude that:
\begin{equation*}
    \mathcal{J}_\varepsilon(u^{**}, v^{**}) = \min_{[u, v] \in \mathscr{U}_\mathrm{ad}^{K}} \mathcal{J}_{\varepsilon}(u, v),
\end{equation*}
and complete the proof of Main Theorem \ref{mainTh02} (\hyperlink{II-B}{II-B}).
\qed

\section{Proof of Main Theorem \ref{mainTh03}}

Throughout this Section, we suppose the situation (\hyperlink{rs0}{r.s.0}). Let $ \varepsilon > 0 $ be a fixed constant, and let $ [\eta_0, \theta_0] \in D_0 $ be the initial pair. Let us take any forcing pair $ [u, v] \in \mathscr{X} $ $ (= L^\infty(Q) \times \mathscr{H}) $, and consider the unique solution $ [\eta, \theta] \in [\mathscr{H}]^2 $ to the state-system (S)$_\varepsilon$. Also, let us take any constant $ \delta \in (0, 1) $ and any pair of functions $ [h, k] \in \mathscr{X}$, and consider another solution $ [\eta^\delta, \theta^\delta] \in [\mathscr{H}]^2 $ to the system (S)$_\varepsilon$, for the initial pair $ [\eta_0, \theta_0] $ and a perturbed forcing pair $ [u +\delta h, v +\delta k] $. On this basis, we consider a sequence of pairs of functions $ \{ [\chi^\delta, \gamma^\delta] \}_{\delta \in (0, 1)} \subset [\mathscr{H}]^2 $, defined as:
\begin{equation}\label{set00}
    [\chi^\delta, \gamma^\delta] := \left[ \frac{\eta^\delta -\eta}{\delta}, ~ \frac{\theta^\delta -\theta}{\delta} \right] \in [\mathscr{H}]^2, \mbox{ for $ \delta \in (0, 1)$.}
\end{equation}
This sequence acts a key-role in the computation of G\^{a}teaux differential of the cost function $ \mathcal{J}_\varepsilon $, for $ \varepsilon > 0 $. 
\begin{rem}\label{Rem.GD01}
    Note that for any  $ \delta \in (0, 1) $, the pair of functions $ [\chi^\delta, \gamma^\delta] \in [\mathscr{H}]^2 $ fulfills the following variational forms:
\begin{align*}
    (\partial_{t} & \chi^\delta(t), \varphi)_{H} + (\nabla \chi^\delta(t), \nabla\varphi)_{[H]^N} 
    \\
    & +\int_{\Omega} \left(\int_{0}^{1} g'(\eta(t)+\varsigma \delta \chi^\delta(t)) \, d \varsigma \right) \chi^\delta(t) \varphi \, dx
    \\
    &  + \int_{\Omega} \left(f_{\varepsilon}(\nabla \theta(t)) \int_{0}^{1} \alpha''(\eta(t) +\varsigma \delta \chi^\delta(t)) \, d\varsigma \right) \chi^\delta(t) \varphi \, dx
    \\
    & +\int_{\Omega}\left(\alpha'(\eta^\delta(t)) \int_{0}^{1} \nabla f_{\varepsilon}(\nabla \theta(t) +\varsigma \delta \nabla \gamma^\delta(t)) \, d\varsigma \right)\cdot \nabla \gamma^\delta(t) \varphi \, dx
    \\
    = & (M_u h(t), \varphi)_{H}, \mbox{ for any $ \varphi \in V $, a.e. $ t \in (0, T) $, subject to $ \chi^\delta(0) = 0 $ in $ H $,}
\end{align*}
and
\begin{align*}
    (\alpha_{0}(t) & \partial_{t} \gamma^\delta(t), \psi)_{H} + \nu^{2}(\nabla \gamma^\delta(t), \nabla \psi)_{[H]^N} 
    \nonumber 
    \\
 & + \int_{\Omega} \left(\alpha(\eta^\delta(t)) \int_{0}^{1} \nabla^2 f_{\varepsilon}(\nabla \theta(t) + \varsigma \delta \nabla \gamma^\delta(t)) \, d \varsigma \right) \nabla \gamma^\delta(t) \cdot\nabla \psi \, dx 
    \nonumber 
    \\
& + \int_{\Omega}\left( \left(\int_{0}^{1}\alpha'(\eta(t) +\varsigma \delta \chi^\delta(t)) \, d \varsigma \right)\chi^\delta(t) \right)\nabla f_{\varepsilon}(\nabla \theta(t))\cdot \nabla \psi \, dx 
    \nonumber
    \\
    = & (M_v k(t), \psi)_{H}, \mbox{ for any $ \psi \in V_{0} $,  a.e. $ t \in (0, T) $, subject to $ \gamma^\delta(0) = 0 $ in $ H $.}
\end{align*}
    In fact, these variational forms are obtained by taking the difference between respective two variational forms for $ [\eta^\delta, \theta^\delta] $ and $ [\eta, \theta] $, as in Main Theorem \ref{mainTh01} (\hyperlink{I-A}{I-A}), and by using the following linearization formulas: 

\begin{align*}
    & \frac{1}{\delta} \bigl( g(\eta^{\delta})-g(\eta) \bigr) = \left( \int_{0}^{1}g'(\eta + \varsigma \delta \chi^\delta) \, d\varsigma \right) \chi^{\delta} \mbox{ in $ \mathscr{H} $,}
\end{align*}
\begin{align*}
    & \frac{1}{\delta} \bigl( \alpha' (\eta^{\delta})f_{\varepsilon}(\nabla \theta^{\delta}) - \alpha'(\eta)f_{\varepsilon}(\nabla \theta) \bigr)
    \\
    & \qquad = \frac{1}{\delta} \bigl( \alpha' (\eta^{\delta})-\alpha'(\eta) \bigr) f_{\varepsilon}(\nabla \theta) + \frac{1}{\delta} \alpha'(\eta^\delta) \bigl( f_{\varepsilon}(\nabla \theta^{\delta})-f_{\varepsilon}(\nabla \theta) \bigr)
    \\
    & \qquad = \left( f_{\varepsilon}(\nabla \theta) \int_{0}^{1}\alpha''(\eta+\varsigma \delta \chi^\delta) \, d\varsigma \right) \chi^{\delta}
    \\
    &  \qquad \qquad + \left( \alpha'(\eta^{\delta}) \int_{0}^{1}\nabla  f_{\varepsilon}(\nabla \theta +\varsigma \delta \nabla \gamma^\delta) \, d\varsigma \right) \cdot \nabla \gamma^{\delta} \mbox{ in $ \mathscr{H} $,}
\end{align*}
and
\begin{align*}
    & \frac{1}{\delta} \bigl( \alpha(\eta^{\delta}) \nabla f_{\varepsilon}(\nabla\theta^{\delta}) - \alpha(\eta)\nabla f_{\varepsilon}(\nabla\theta) \bigr)
    \\
    & \qquad = \frac{1}{\delta} \alpha(\eta^{\delta}) \bigl( \nabla f_{\varepsilon}(\nabla \theta^{\delta})-\nabla f_{\varepsilon}(\nabla \theta) \bigr) + \frac{1}{\delta} \bigl( \alpha (\eta^{\delta})-\alpha(\eta) \bigr) \nabla f_{\varepsilon}(\nabla \theta)
    \\
    & \qquad = \left( \alpha(\eta^{\delta}) \int_{0}^{1}\nabla^2 f_{\varepsilon}(\nabla \theta +\varsigma \delta \nabla \gamma^\delta) \, d\varsigma \right) \nabla \gamma^{\delta}
    \\
    & \qquad \qquad +\left(\left( \int_{0}^{1}\alpha'(\eta +\varsigma \delta \chi^\delta) \, d \varsigma \right) \chi^{\delta}\right) \nabla f_{\varepsilon}(\nabla \theta)\mbox{ in $ [\mathscr{H}]^N $.}
\end{align*}
Incidentally, the above linearization formulas can be verified as consequences of the assumptions (\hyperlink{A1l}{A1})--(\hyperlink{A4l}{A4}) and the mean-value theorem (cf. \cite[Theorem 5 in p. 313]{lang1968analysisI}).
\end{rem}

\begin{rem}\label{Rem.Ken01}
    Note that the situation (\hyperlink{rs0}{r.s.0}) implies $ \eta_0 \in L^\infty(\Omega) $ and $ u \in L^\infty(Q) $. Therefore, under (\hyperlink{rs0}{r.s.0}), we can suppose $  \eta \in L^\infty(Q) $ for the solution $ [\eta, \theta] \in [\mathscr{H}]^2 $ to the system (S)$_\varepsilon$. 
\end{rem}

Now, we prepare the following two Lemmas, for the proof of Main Theorem \ref{mainTh03}. 
\begin{lem}\label{Lem.GD01}
    Under the assumptions (\hyperlink{A1l}{A1})--(\hyperlink{A5l}{A5}), let us fix $ \varepsilon > 0 $, and suppose (\hyperlink{rs0}{r.s.0}) as in Main Theorem \ref{mainTh03}. Then, the restriction of the cost $\mathcal{J}_\varepsilon|_{\mathscr{X}}$ : $\mathscr{X} \longrightarrow \mathbb{R}$ is G\^{a}teaux differentiable over $\mathscr{X}$. Moreover, for any $ [u, v] \in \mathscr{X}$, the G\^{a}teaux derivative $ \left(\mathcal{J}_\varepsilon |_{\mathscr{X}} \right)'(u, v) \in \mathscr{X}^* $ admits a unique extension $ \mathcal{J}_\varepsilon'(u, v) \in ([\mathscr{H}]^2)^*  $ $ = [\mathscr{H}]^2 $, such that: 
    \begin{align}\label{GD00}
    \mathcal{J}_\varepsilon'(u, v) = \left(\mathcal{J}_\varepsilon |_{\mathscr{X}} \right)'(u, v)\ \mbox{in}\ \mathscr{X}^*,
    \end{align}
    and
    \begin{align}\label{GD01}
        \bigl( \mathcal{J}_\varepsilon'(u, v), & [h, k] \bigr)_{[\mathscr{H}]^2} = \bigl( [M_\eta (\eta -\eta_\mathrm{ad}), M_\theta (\theta -\theta_\mathrm{ad})], \bar{\mathcal{P}}_\varepsilon [M_u h, M_v k] \bigr)_{[\mathscr{H}]^2}
        \nonumber
        \\
        & +\bigl( [M_u u, M_v v], [h, k] \bigr)_{[\mathscr{H}]^2}, \mbox{ for any $ [h, k] \in \mathscr{X}$.}
    \end{align}
    In the context, $ [\eta, \theta] $ is the solution to the state-system (S)$_\varepsilon$, for the initial pair $[\eta_0, \theta_0]$ and forcing pair $ [u, v] $, and $ \bar{\mathcal{P}}_\varepsilon : [\mathscr{H}]^2 \longrightarrow \mathscr{Z} $ is a bounded linear operator, which is given as a restriction $ \mathcal{P}|_{\{[0, 0]\} \times [\mathscr{H}]^2} $ of the (linear) isomorphism $ \mathcal{P} = \mathcal{P}(a, b, \mu, \lambda, \omega, A) : [H]^2 \times \mathscr{Y}^* \longrightarrow \mathscr{Z} $, as in Proposition \ref{ASY_Cor.1}, in the case when:
    \begin{equation}\label{set01}
        \begin{cases}
            [a, b]= [\alpha_0, 0] \mbox{ in $ W^{1, \infty}(Q) \times L^\infty(Q) $,}
            \\[1ex]
            \mu = \bar{\mu}_\varepsilon := \alpha''(\eta) f_\varepsilon(\nabla \theta) \mbox{ in $ L^\infty(0, T; H) $,} 
            \\[1ex]
            \lambda = \bar{\lambda}_\varepsilon := g'(\eta) \mbox{ in $ L^\infty(Q) $,}
            \\[1ex]
            \omega = \bar{\omega}_\varepsilon := \alpha'(\eta) \nabla f_\varepsilon(\nabla \theta)  \mbox{ in $ [L^\infty(Q)]^N $,}
            \\[1ex]
            A = \bar{A}_\varepsilon := \alpha(\eta) \nabla^2 f_\varepsilon(\nabla \theta) \mbox{ in $ [L^\infty(Q)]^{N \times N} $.}
            \end{cases}
    \end{equation}
\end{lem}
\paragraph{Proof.}{Let us fix any $ [u, v] \in \mathscr{X} $, and take any $ \delta \in (0, 1)$ and any $ [h, k] \in \mathscr{X} $. 
Then, due to the assumptions (\hyperlink{rs0}{r.s.0}) and $ [u, v], [h, k] \in \mathscr{X} $, we can see that:
\begin{equation*}
    |\eta_0|_{L^\infty(\Omega)} \vee \sup_{\tilde{\delta} \in (0, 1)} |u +\tilde{\delta} h|_{L^\infty(Q)} < \infty,
\end{equation*}
and
\begin{equation*}
    [M_u(u +\delta h), M_v(v +\delta k)] \to [M_u u, M_v v] \mbox{ in $ \mathscr{X} $, as $ \delta \downarrow 0 $.}
\end{equation*}
Therefore, as a consequence of Main Theorem \ref{mainTh01} (\hyperlink{I-B}{I-B}), it is observed that:
\begin{subequations}\label{ken01}
\begin{align}\label{ken01-01}
    [\eta^\delta, \theta^\delta] & \to [\eta, \theta] \mbox{ in $ [C([0, T]; H)]^2 $, in $ \mathscr{Y} $,}
    \nonumber
    \\
    & \mbox{weakly in $ [W^{1, 2}(0, T; H)]^2 $,}
    \nonumber
    \\
    & \mbox{and weakly-$*$ in $ L^\infty(0, T; V) \times L^\infty(0, T; V_0) $,}
\end{align}
and 
\begin{align}\label{ken01-02}
    \eta^\delta \to \eta \mbox{ weakly-$*$ in $ L^\infty(Q) $, as $ \delta \downarrow 0 $.}
\end{align}
\end{subequations}

In the meantime, it is easily computed that:
\begin{align}\label{GD02}
    \frac{1}{\delta} \bigl(  \mathcal{J}_\varepsilon &( u +\delta h, v +\delta k) -\mathcal{J}_\varepsilon(u, v) \bigr)
    \nonumber
    \\
    = & \left( \frac{M_\eta}{2} (\eta^\delta +\eta -2 \eta_\mathrm{ad}), \chi^\delta \right)_{\hspace{-0.5ex}\mathscr{H}} +\left( \frac{M_\theta}{2} (\theta^\delta +\theta -2 \theta_\mathrm{ad}), \gamma^\delta \right)_{\hspace{-0.5ex}\mathscr{H}}
    \\
    & \quad +\left( \frac{M_u}{2}(2u +\delta h), h \right)_{\hspace{-0.5ex}\mathscr{H}} +\left( \frac{M_v}{2} (2v +\delta k), k \right)_{\hspace{-0.5ex}\mathscr{H}}.
    \nonumber
\end{align}
Here, let us set:
\begin{subequations}\label{sets02}
\begin{equation}\label{set02}
    \begin{cases}
        \ds \bar{\mu}_\varepsilon^\delta := f_\varepsilon(\nabla \theta) \int_0^1 \alpha''(\eta +\varsigma \delta \chi^\delta) \, d \varsigma \mbox{ in $ L^\infty(0, T; H) $,} 
        \\[2ex]
        \ds \bar{\lambda}_\varepsilon^\delta := \int_0^1 g'(\eta +\varsigma \delta \chi^\delta) \, d\varsigma \mbox{ in $ L^\infty(Q) $,}
        \\[2ex]
        \ds \bar{\omega}_\varepsilon^\delta := \alpha'(\eta^\delta) \int_0^1 \nabla f_\varepsilon(\nabla \theta +\varsigma \delta \nabla \gamma^\delta) \, d \varsigma \mbox{ in $ [L^\infty(Q)]^N $,}
        \\[2ex]
        \ds \bar{A}_\varepsilon^\delta := \alpha(\eta^\delta) \int_0^1 \nabla^2 f_\varepsilon(\nabla \theta +\varsigma \delta \nabla \gamma^\delta) \, d \varsigma \mbox{ in $ [L^\infty(Q)]^{N \times N} $,}
    \end{cases} 
\end{equation}
and
\begin{align}\label{set03}
    \ds \bar{k}_\varepsilon^\delta := M_v k + \mbox{div} & \left[ \rule{-1pt}{16pt} \right. 
    \chi^\delta \nabla f_\varepsilon(\nabla \theta) \int_0^1 \alpha'(\eta +\varsigma \delta \chi^\delta) \, d \varsigma 
    \nonumber
    \\
    & \qquad \ds  -\chi^\delta \alpha'(\eta^\delta) \int_0^1 \nabla f_\varepsilon(\nabla \theta +\varsigma \delta \nabla \gamma^\delta) \, d \varsigma 
    \left. \rule{-1pt}{16pt} \right] \mbox{ in $ \mathscr{V}_0^* $,}
    \\
    & \mbox{for all $ \delta \in (0, 1)$.}
    \nonumber
\end{align}
\end{subequations}
Then, in the light of \eqref{ken01} and  Remark \ref{Rem.GD01}, one can say that:
\begin{equation*}
    [\chi^\delta, \gamma^\delta] = \bar{\mathcal{P}}_\varepsilon^\delta [M_u h, \bar{k}_\varepsilon^\delta] \mbox{ in $ \mathscr{Z} $, for $ \delta \in (0, 1) $,}
\end{equation*}
by using the restriction $ \bar{\mathcal{P}}_\varepsilon^\delta := \mathcal{P}|_{\{[0, 0]\} \times \mathscr{Y}^*} : \mathscr{Y}^* \longrightarrow \mathscr{Z} $ of the (linear) isomorphism $ \mathcal{P} = \mathcal{P}(a, b, \mu, \lambda, \omega, A) : [H]^2 \times \mathscr{Y}^* \longrightarrow \mathscr{Z} $, as in Proposition \ref{ASY_Cor.1}, in the case when:
\begin{equation*}
    \begin{cases}
        [a, b, \lambda]= [\alpha_0, 0, \bar{\lambda}_\varepsilon^\delta] \mbox{ in $ W^{1, \infty}(Q) \times [L^\infty(Q)]^2 $,}
        \\[1ex]
        \omega = \bar{\omega}_\varepsilon^\delta \mbox{ in $ [L^\infty(Q)]^N $,}
        \\[1ex]
        A = \bar{A}_\varepsilon^\delta \mbox{ in $ [L^\infty(Q)]^{N \times N} $,}
        \\[1ex]
        \mu = \bar{\mu}_\varepsilon^\delta \mbox{ in $ L^\infty(0, T; H) $, for $ \delta \in (0, 1) $.} 
    \end{cases}
\end{equation*}
Besides, taking into account \eqref{f_eps}, \eqref{C_0*}, \eqref{sets02}, (\hyperlink{A3l}{A3}), (\hyperlink{A4l}{A4}), and Remark \ref{Rem.Prelim01}, we have:
\begin{subequations}\label{est3-01}
    \begin{align}\label{est3-01a}
        \bar{C}_0^*  := & \frac{9(1 +\nu^2)}{1 \wedge \nu^2 \wedge \delta_*} \cdot \bigl( 1 +(C_{V}^{L^4})^2 +(C_{V}^{L^4})^4 +(C_{V_0}^{L^4})^2 \bigr)
        \nonumber
        \\
        & \quad \cdot 
        \bigl( 1 +|\alpha_0|_{W^{1, \infty}(Q)} +|g'|_{L^\infty(\mathbb{R})} +|\alpha'|_{L^\infty(\mathbb{R})}^2 \bigr)
        \\
        \geq & \frac{9(1 +\nu^2)}{1 \wedge \nu^2 \wedge \inf \alpha_0(Q)}  \cdot \bigl( 1 +(C_{V}^{L^4})^2 +(C_{V}^{L^4})^4 +(C_{V_0}^{L^4})^2 \bigr)
        \nonumber
        \\
        & \quad \cdot \sup_{\delta \in (0, 1)}
        \bigl( 1 +|\alpha_0|_{W^{1, \infty}(Q)} +|\bar{\lambda}_\varepsilon^\delta|_{L^\infty(Q)} +|\bar{\omega}_\varepsilon^\delta|_{[L^\infty(Q)]^N}^2 \bigr),
        \nonumber
    \end{align}
and
    \begin{align}\label{est3-01b}
    \bigl| \bigl< [M_u\,& h(t), \bar{k}_\varepsilon^\delta(t)], [\varphi, \psi] \bigr>_{V \times V_0} \bigr| \leq |\langle M_u h(t), \varphi \rangle_{V}| +|\langle \bar{k}_\varepsilon^\delta(t), \psi \rangle_{V_0}| 
    \nonumber
        \\[0.5ex]
    \leq \, & M_u| h(t)|_H |\varphi|_H +M_v |k(t)|_H |\psi|_H +2|\alpha'|_{L^\infty(\R)}|\chi^\delta(t)|_H|\nabla \psi|_{[H]^N}
    \nonumber
        \\[0.5ex]
        \leq \, & M_u| h(t)|_H |\varphi|_V +\bigl( M_v C_{V_0}^{H} |k(t)|_H +2|\alpha'|_{L^\infty(\R)}|\chi^\delta(t)|_H \bigr) |\psi|_{V_0},
         \\[0.5ex]
    & \mbox{for a.e. $ t \in (0, T) $, any $ [\varphi, \psi] \in V \times V_0 $, and any $  \delta \in (0, 1)$,}
    \nonumber
\end{align} 
with use of the constant $C_{V_0}^{H} > 0$ of the embedding $V_0 \subset H$,
so that 
    \begin{align}\label{est3-01c}
        \bigl| [M_u h(t), \, & \bar{k}_\varepsilon^\delta(t)] \bigr|_{V^* \times V_0^*}^2 \leq \bar{C}_1^* \bigl( \bigl| [h(t), k(t)] \bigr|_{[H]^2}^2 +|\chi^\delta(t)|_H^2 \bigr),
    \nonumber
        \\
    & \mbox{for a.e. $ t \in (0, T) $, and any $ \delta \in (0, 1)$,}
\end{align}
\end{subequations}
with a positive constant $ \bar{C}_1^* := 4 \bigl( M_u^2 + M_v^2({C_{V_0}^H})^2 + |\alpha'|_{L^\infty(\R)}^2 \bigr) $.
\medskip

Now, having in mind \eqref{est3-01}, let us apply Proposition \ref{Prop(I-B)} to the case when:
\begin{equation*}
    \begin{array}{c}
    \begin{cases}
        [a^1, b^1, \mu^1, \lambda^1, \omega^1, A^1] = [a^2, b^2, \mu^2, \lambda^2, \omega^2, A^2] = [\alpha_0, 0, \overline{\mu}_\varepsilon^\delta, \bar{\lambda}_\varepsilon^\delta, \bar{\omega}_\varepsilon^\delta, \bar{A}_\varepsilon^\delta],
        \\[1ex]
        [p_0^1, z_0^1] = [p_0^2, z_0^2] = [0, 0], ~ [h^1, k^1] = [M_u h, \bar{k}_\varepsilon^\delta], ~ [h^2, k^2] = [0, 0],
        \\[1ex]
        [p^1, z^1] = [\chi^\delta, \gamma^\delta] = \bar{\mathcal{P}}_\varepsilon^\delta [M_u h, \bar{k}_\varepsilon^\delta], ~ [p^2, z^2] = [0, 0] = \bar{\mathcal{P}}_\varepsilon^\delta [0, 0], \mbox{ \ for $ \delta \in (0, 1) $.}
    \end{cases}
    \end{array}
\end{equation*}
Then, we estimate that:
\begin{align*}
    \frac{d}{dt} & \bigl( |\chi^\delta(t)|_H^2 +|\sqrt{\alpha_0(t)} \gamma^\delta(t)|_H^2 \bigr) +\bigl( |\chi^\delta(t)|_V^2 +\nu^2 |\gamma^\delta(t)|_{V_0}^2 \bigr)
    \\
    & \leq 3 \bar{C}_0^* \bigl( |\chi^\delta(t)|_H^2 +|\sqrt{\alpha_0(t)} \gamma^\delta(t)|_H^2 \bigr) +2\bar{C}_0^* \big( |M_u h(t)|_{V^*}^2 +|\bar{k}_\varepsilon^\delta(t)|_{V_0^*}^2 \bigr)
    \\
    & \leq 3 \bar{C}_0^* (1 +\bar{C}_1^*) \bigl( |\chi^\delta(t)|_H^2 +|\sqrt{\alpha_0(t)} \gamma^\delta(t)|_H^2 \bigr) +2\bar{C}_0^* \bar{C}_1^* \big( |h(t)|_{H}^2 +|k(t)|_{H}^2 \bigr), 
    \\
    & \qquad\mbox{for a.e. $ t \in (0, T) $,}
\end{align*}
and subsequently, by using Gronwall's lemma, we observe that:
\begin{description}
    \item[\textmd{$(\hypertarget{star1}{\star\,1})$}]the sequence $ \{ [\chi^\delta, \gamma^\delta] \}_{\delta \in (0, 1)} $ is bounded in $ [C([0, T]; H)]^2 \cap \mathscr{Y} $. 
\end{description}

Meanwhile, as consequences of \eqref{set00}, \eqref{set01}--\eqref{est3-01}, $(\hyperlink{star1}{\star\,1})$, (\hyperlink{A1l}{A1})--(\hyperlink{A5l}{A5}), Main Theorem \ref{mainTh01}, Remark \ref{Rem.mTh01Conv}, and the dominated convergence theorem \cite[Theorem 10 on page 36]{MR0492147}, one can find a sequence $ \{ \delta_n \}_{n = 1}^\infty \subset \R $, such that:
\begin{subequations}\label{convs3-01}
    \begin{align}\label{convs3-00a}
        0 < |\delta_n| < 1, \mbox{ and } \delta_n \to 0, \mbox{ as $ n \to \infty $,}
    \end{align}
    \begin{align}\label{convs3-01a}
        &
        \begin{cases}
            [\delta_n \chi^{\delta_n}, \delta_n \gamma^{\delta_n}] = [\eta^{\delta_n} -\eta, \theta^{\delta_n} -\theta] \to [0, 0]
            \\
            \qquad \mbox{in $ [C([0, T]; H)]^2 $, and in $ \mathscr{Y} $,}
            \\[1ex]
            [\delta_n \nabla \chi^{\delta_n}, \delta_n \nabla \gamma^{\delta_n}] = [\nabla (\eta^{\delta_n} -\eta), \nabla (\theta^{\delta_n} -\theta)] \to [0, 0]
            \\
            \qquad \mbox{in $ [L^2(0, T; [H]^N)]^2$, and in the pointwise sense a.e. in $ Q $,}
        \end{cases}
        \mbox{as $ n \to \infty $,}
    \end{align}
    \begin{align}\label{convs3-01b}
        [\bar{\lambda}_\varepsilon^{\delta_n}, & \bar{\omega}_\varepsilon^{\delta_n}, \bar{A}_\varepsilon^{\delta_n}] \to [\bar{\lambda}_\varepsilon, \bar{\omega}_\varepsilon, \bar{A}_\varepsilon]  \mbox{ weakly-$ * $ in $ L^\infty(Q) \times [L^\infty(Q)]^N \times [L^\infty(Q)]^{N \times N}$,}
        \nonumber
        \\
        & \mbox{and in the pointwise sense a.e. in $ Q $, ~ as $ n \to \infty $,}
    \end{align}
    \begin{equation}\label{convs3-01c}
        \begin{cases}
            \bar{\mu}_\varepsilon^{\delta_n} \to \bar{\mu}_\varepsilon \mbox{ weakly-$ * $ in $ L^\infty(0, T; H) $,}
            \\[1ex]
            \bar{\mu}_\varepsilon^{\delta_n}(t) \to \bar{\mu}_\varepsilon(t) \mbox{ in $ H $,}
            \mbox{ for a.e. $ t \in (0, T) $,}
        \end{cases}
        \mbox{as $ n \to \infty $,}
    \end{equation}
    and
    \begin{align}\label{convs3-01d}
        \langle \bar{k}_\varepsilon^{\delta_n} -M_v k, \psi \rangle_{\mathscr{V}_0} =  &~ -\left( \chi^{\delta_n}, ~ \nabla f_\varepsilon(\nabla \theta) \left( \rule{-1pt}{14pt} \right. \int_0^1 \alpha'(\eta +\varsigma \delta_n \chi^{\delta_n}) \, d \varsigma  \left. \rule{-1pt}{14pt} \right)\cdot \nabla \psi\right)_{\hspace{-0.5ex}\mathscr{H}}
        \nonumber
        \\[1ex]
        &~ +\left( \chi^{\delta_n}, ~ 
        \alpha'(\eta^{\delta_n}) \left( \rule{-1pt}{14pt} \right. \int_0^1 \nabla f_\varepsilon(\nabla \theta +\varsigma \delta_n \nabla \gamma^{\delta_n}) \, d \varsigma \left. \rule{-1pt}{14pt} \right)\cdot \nabla \psi \right)_{\hspace{-0.5ex}\mathscr{H}}
        \\[1ex]
        \to &~ 0, \mbox{ as $ n  \to \infty $.}
        \nonumber
    \end{align}
\end{subequations}

On account of \eqref{set00} and \eqref{set01}--\eqref{convs3-01}, we can apply Proposition \ref{ASY_Cor.2} (\hyperlink{B}{B}), and can see that:
\begin{align}\label{conv3-02}
    [\chi^{\delta_n}, \gamma^{\delta_n}] & = \bar{\mathcal{P}}_\varepsilon^{\delta_n}[M_u h, \bar{k}_\varepsilon^{\delta_n}] \to [\chi, \gamma] := \bar{\mathcal{P}}_\varepsilon[M_u h, M_v k] \mbox{ in $ [\mathscr{H}]^2 $, weakly in $ \mathscr{Y} $,}
    \nonumber
    \\
    & \mbox{and weakly in $ W^{1, 2}(0, T; V^*) \times W^{1, 2}(0, T; V_0^*) $, as $ n \to \infty $.}
\end{align}
Since the Hilbert space $ \mathscr{Y} $ is separable, and the uniqueness of the solution $ [\chi, \gamma] = \bar{\mathcal{P}}_\varepsilon [M_u h, M_v k] $ is guaranteed by Proposition \ref{Prop(I-A)}, the observations \eqref{GD02}, \eqref{convs3-01}, and \eqref{conv3-02} enable us to compute the directional derivative $ D_{[h, k]}\mathcal{J}_\varepsilon(u, v) \in \R $, as follows:
\begin{align}\label{conv3-03}
    D_{[h, k]} & \mathcal{J}_\varepsilon(u, v) := \lim_{\delta \to 0} \frac{1}{\delta} \bigl( \mathcal{J}_\varepsilon(u +\delta h, v +\delta k) -\mathcal{J}_\varepsilon(u, v) \bigr) \nonumber
    \\
    = & \bigl( [M_\eta (\eta -\eta_\mathrm{ad}), M_\theta (\theta -\theta_\mathrm{ad})], \bar{\mathcal{P}}_\varepsilon [M_u h, M_v k] \bigr)_{[\mathscr{H}]^2}
    +\bigl( [M_u u, M_v v], [h, k] \bigr)_{[\mathscr{H}]^2}, \nonumber
    \\
    & \quad \mbox{for any direction $ [h, k] \in \mathscr{X} $.}
\end{align}
Moreover, in the light of \eqref{set01}, \eqref{conv3-03}, and Proposition \ref{ASY_Cor.1}, we can observe that:
\begin{description}
    \item[\textmd{$(\hypertarget{star2}{\star\,2})$}] the mapping $[h, k] \in \mathscr{X} \mapsto D_{[h, k]}\mathcal{J}_\varepsilon(u, v) \in \mathbb{R}$ is a linear functional;
    \item[\textmd{$(\hypertarget{star3}{\star\,3})$}] there exists a constant $ M_1^{**} $, independent of $ [h, k] \in \mathscr{X} $, such that 
    \begin{align*}
        |D_{[h, k]}\mathcal{J}_\varepsilon(u, v)| \leq M_1^{**}|[h, k]|_{[\mathscr{H}]^2},\ \mbox{for any}\ [h, k] \in \mathscr{X}.
    \end{align*}
\end{description}
As a consequence of $(\hyperlink{star2}{\star\,2})$, $(\hyperlink{star3}{\star\,3})$, the continuous and dense embedding $ \mathscr{X} \subset [\mathscr{H}]^2 $, and Riesz's theorem, we can obtain the required functional $ \mathcal{J}'_\varepsilon (u, v) \in \bigl( [\mathscr{H}]^2 \bigr)^* $ $ \bigl( = [\mathscr{H}]^2 \bigr) $, satisfying \eqref{GD00} and \eqref{GD01}, as the unique extension of the G\^{a}teaux differential $\left(\mathcal{J}_\varepsilon|_{\mathscr{X}} \right)'(u, v) \in \mathscr{X}^*$ at $[u, v] \in \mathscr{X}$. 
\medskip

Thus, we complete the proof of this lemma. 
\hfill \qed
}

\begin{lem}\label{Lem.GD02}
    Under the assumptions (\hyperlink{A1l}{A1})--(\hyperlink{A5l}{A5}) with (\hyperlink{rs0}{r.s.0}), let $ [u_\varepsilon^*, v_\varepsilon^*] \in \mathscr{U}_\mathrm{ad}^K $ be an optimal control of the problem (\hyperlink{OP}{OP})$_\varepsilon^K$, and let $ [\eta_\varepsilon^*, \theta_\varepsilon^*] $ be the solution to the system (S)$_\varepsilon$, for the initial pair $[\eta_0, \theta_0] $ and forcing pair $ [u_\varepsilon^*, v_\varepsilon^*] $. Also, let $ \mathcal{P}_\varepsilon^* : [\mathscr{H}]^2 \longrightarrow \mathscr{Z} $ be the bounded linear operator, defined in Remark \ref{Rem.mTh03}, with use of the solution $ [\eta_\varepsilon^*, \theta_\varepsilon^*] $. Let $ \mathcal{P}_\varepsilon : [\mathscr{H}]^{2} \longrightarrow \mathscr{Z} $ be a bounded linear operator, which is defined as a restriction $ \mathcal{P}|_{\{[0, 0]\} \times [\mathscr{H}]^2} $ of the linear isomorphism $ \mathcal{P} = \mathcal{P}(a, b, \mu, \lambda, \omega, A) : [H]^2 \times \mathscr{Y}^* \longrightarrow \mathscr{Z} $, as in Proposition \ref{ASY_Cor.1}, in the case when:
    \begin{equation}\label{set_GD02-01}
        \begin{cases}
            [a, b]= [\alpha_0, 0] \mbox{ in $ W^{1, \infty}(Q) \times L^\infty(Q) $,}
            \\[1ex]
            \mu = \alpha''(\eta_\varepsilon^*) f_\varepsilon(\nabla \theta_\varepsilon^*) \mbox{ in $ L^\infty(0, T; H) $,} 
            \\[1ex]
            \lambda = g'(\eta_\varepsilon^*) \mbox{ in $ L^\infty(Q) $,}
            \\[1ex]
            \omega = \alpha'(\eta_\varepsilon^*) \nabla f_\varepsilon(\nabla \theta_\varepsilon^*) \mbox{ in $ [L^\infty(Q)]^N $,}
            \\[1ex]
            A = \alpha(\eta_\varepsilon^*) \nabla^2 f_\varepsilon(\nabla \theta_\varepsilon^*) \mbox{ in $ [L^\infty(Q)]^{N\times N} $.}
            \end{cases}
    \end{equation}
    Then, the operators $ \mathcal{P}_\varepsilon^* $ and $ \mathcal{P}_\varepsilon $  have a conjugate relationship, in the following sense:
    \begin{align*}
        \bigl( \mathcal{P}_\varepsilon^*[u, v] & , [h, k] \bigr)_{[\mathscr{H}]^2} = \bigl( [u, v], \mathcal{P}_\varepsilon [h, k] \bigr)_{[\mathscr{H}]^2}, 
        \\
        & \mbox{for all $ [h, k], [u, v] \in [\mathscr{H}]^2 $.}
    \end{align*}
\end{lem}
\paragraph{Proof.}{
    Let us fix arbitrary pairs of functions $ [h, k], [u, v] \in [\mathscr{H}]^2 $, and let us put:
    \begin{equation*}
    [\chi_\varepsilon, \gamma_\varepsilon] := \mathcal{P}_\varepsilon [h, k] \quad \mbox{ and } \quad [p_\varepsilon, z_\varepsilon] := \mathcal{P}_\varepsilon^* [u, v], \mbox{ in $ [\mathscr{H}]^2 $.} 
\end{equation*}
Then, invoking Proposition \ref{Prop(I-A)}, and the settings as in \eqref{setRem4} and \eqref{set_GD02-01}, we compute that:
\begin{align*}
    \bigl( \mathcal{P}_\varepsilon^* & [u, v], [h, k] \bigr)_{[\mathscr{H}]^2} = \int_0^T \bigl( p_\varepsilon(t), h(t) \bigr)_{H} \, dt +\int_0^T \bigl( z_\varepsilon(t), k(t) \bigr)_{H} \, dt
    \hspace{19.5ex}
    \\
    = & \int_0^T \langle h(t), p_\varepsilon(t) \rangle_{V} \, dt +\int_0^T \langle k(t), z_\varepsilon(t) \rangle_{V_0} \, dt
\end{align*}
\begin{align*}
    = & \int_0^T \left[ \rule{-1pt}{16pt} \right. \bigl< \partial_t \chi_\varepsilon(t), p_\varepsilon(t) \bigr>_{V} +\bigl( \nabla \chi_\varepsilon(t), \nabla p_\varepsilon(t) \bigr)_{[H]^N} 
    \\
    & \qquad +\bigl( \alpha''(\eta_\varepsilon^*(t)) f_\varepsilon(\nabla \theta_\varepsilon^*(t)) \chi_\varepsilon(t),  p_\varepsilon(t) \bigr)_H 
    \\
    & \qquad +\bigl( g'(\eta_\varepsilon^*(t)) \chi_\varepsilon(t), p_\varepsilon(t) \bigr)_H +\bigl( \alpha'(\eta_\varepsilon^*(t))\nabla f_\varepsilon(\nabla \theta_\varepsilon^*(t)) \cdot \nabla \gamma_\varepsilon(t), p_\varepsilon(t) \bigr)_{H} \left. \rule{-1pt}{16pt} \right] \, dt
    \\
    & +\int_0^T \left[ \rule{-1pt}{16pt} \right. \bigl<\partial_t \gamma_\varepsilon(t),  \alpha_0(t) z_\varepsilon(t) \bigr>_{V_0} +\bigl( \alpha'(\eta_\varepsilon^*(t)) \chi_\varepsilon(t) \nabla f_\varepsilon(\nabla \theta_\varepsilon^*(t)), \nabla z_\varepsilon(t) \bigr)_{[H]^N} 
   \\
    & \qquad +\bigl( \alpha(\eta_\varepsilon^*(t)) \nabla^2 f_\varepsilon(\nabla \theta_\varepsilon^*(t)) \nabla \gamma_\varepsilon(t), \nabla z_\varepsilon(t) \bigr)_{[H]^N} +\nu^2 \bigl( \nabla \gamma_\varepsilon(t), \nabla z_\varepsilon(t) \bigr)_{[H]^N} \left. \rule{-1pt}{16pt} \right] \, dt
\end{align*}
\begin{align*}
    = & \bigl( p_\varepsilon(T), \chi_\varepsilon(T) \bigr)_H -\bigl( p_\varepsilon(0), \chi_\varepsilon(0) \bigr)_H +\int_0^T \left[ \rule{-1pt}{16pt} \right. \bigl< -\partial_t p_\varepsilon(t), \chi_\varepsilon(t) \bigr>_{V} 
    \\
    & \qquad +\bigl( \nabla p_\varepsilon(t), \nabla \chi_\varepsilon(t) \bigr)_{[H]^N} +\bigl( \alpha''(\eta_\varepsilon^*(t)) f_\varepsilon(\nabla \theta_\varepsilon^*(t)) p_\varepsilon(t), \chi_\varepsilon(t) \bigr)_H 
    \\
& \qquad +\bigl( g'(\eta_\varepsilon^*(t)) p_\varepsilon(t), \chi_\varepsilon(t) \bigr)_H +\bigl( \alpha'(\eta_\varepsilon^*(t))\nabla f_\varepsilon(\nabla \theta_\varepsilon^*(t)) \cdot \nabla z_\varepsilon(t), \chi_\varepsilon(t) \bigr)_{H} \left. \rule{-1pt}{16pt} \right] \, dt
    \\
    & +\bigl( \alpha_0(T) z_\varepsilon(T), \gamma_\varepsilon(T) \bigr)_H -\bigl( \alpha_0(0) z_\varepsilon(0), \gamma_\varepsilon(0) \bigr)_H 
    \\
    & \qquad +\int_0^T \left[ \rule{-1pt}{16pt} \right. \bigl< -\partial_t \bigl(\alpha_0 z_\varepsilon)(t),  \gamma_\varepsilon(t) \bigr>_{V_0} +\bigl( \alpha'(\eta_\varepsilon^*(t)) p_\varepsilon(t) \nabla f_\varepsilon(\nabla \theta_\varepsilon^*(t)), \nabla \gamma_\varepsilon(t) \bigr)_{[H]^N} 
    \\
    & \qquad +\bigl( \alpha(\eta_\varepsilon^*(t)) \nabla^2 f_\varepsilon(\nabla \theta_\varepsilon^*(t)) \nabla z_\varepsilon(t), \nabla \gamma_\varepsilon(t) \bigr)_{[H]^N} +\nu^2 \bigl( \nabla z_\varepsilon(t), \nabla \gamma_\varepsilon(t) \bigr)_{[H]^N} \left. \rule{-1pt}{16pt} \right] \, dt
    \\
    = & ( u, \chi_\varepsilon )_{\mathscr{H}} +( v, \gamma_\varepsilon )_{\mathscr{H}} = \bigl( [u, v], \mathcal{P}_\varepsilon [h, k]  \bigr)_{[\mathscr{H}]^2}.
\end{align*}
This finishes the proof of Lemma \ref{Lem.GD02}. \qed
}
\begin{rem}\label{Rem.GD02}
    Note that the operator $ \mathcal{P}_\varepsilon \in \mathscr{L}([\mathscr{H}]^2; \mathscr{Z}) $, as in Lemma \ref{Lem.GD02}, corresponds to the operator $ \bar{\mathcal{P}}_\varepsilon \in \mathscr{L}([\mathscr{H}]^2; \mathscr{Z}) $, as in the previous Lemma \ref{Lem.GD01}, under the special setting \eqref{set_GD02-01}. 
\end{rem}
\bigskip

Now, we are ready to prove the Main Theorem \ref{mainTh03} (\hyperlink{III-A}{III-A}).

\paragraph{Proof of Main Theorem \ref{mainTh03} (III-A).}{
    Let $ [u_\varepsilon^*, v_\varepsilon^*] \in \mathscr{U}_\mathrm{ad}^K $ be the optimal control of (\hyperlink{OP}{OP})$_\varepsilon^K$, with the solution $ [\eta_\varepsilon^*, \theta_\varepsilon^*] \in [\mathscr{H}]^2 $ to the system (S)$_\varepsilon$ for the initial pair $[\eta_0, \theta_0] \in D_0 $, as in (\hyperlink{rs0}{r.s.0}), and forcing pair $ [u_\varepsilon^*, v_\varepsilon^*] $, and let $ \mathcal{P}_\varepsilon, \mathcal{P}_\varepsilon^* \in \mathscr{L}([\mathscr{H}]^2; \mathscr{Z}) $ be the two operators as in Lemma \ref{Lem.GD02}. 
    In addition, let us put $[p_\varepsilon^* , z_\varepsilon^*]:=\mathcal{P}_\varepsilon^*[M_\eta (\eta_\varepsilon^* -\eta_\mathrm{ad}), M_\theta (\theta_\varepsilon^* -\theta_\mathrm{ad})]$.
    Then, on the basis of the previous Lemmas \ref{Lem.GD01} and \ref{Lem.GD02}, we compute that:
    \begin{align}\label{nec01}
        0 \leq &~ \left(\mathcal{J}_\varepsilon'(u_\varepsilon^*, v_\varepsilon^*), [h, k]\right)_{[\mathscr{H}]^2} 
        \nonumber
        \\[1ex]
        = &~ \lim_{\delta \downarrow 0} \frac{1}{\delta} \bigl( \mathcal{J}_\varepsilon(u_\varepsilon^* +\delta (h-u_\varepsilon^*), v_\varepsilon^* +\delta k) -\mathcal{J}_\varepsilon(u_\varepsilon^*, v_\varepsilon^*) \bigr) 
        \nonumber
        \\[1ex]
        = &~ \bigl( [M_\eta (\eta_\varepsilon^* -\eta_\mathrm{ad}), M_\theta (\theta_\varepsilon^* -\theta_\mathrm{ad})], \mathcal{P}_\varepsilon [M_u(h-u_\varepsilon^*), M_v k] \bigr)_{[\mathscr{H}]^2} 
        \nonumber
        \\[0ex]
        & \hspace{4ex} +\bigl( [M_u u_\varepsilon^*, M_v v_\varepsilon^*], [h-u_\varepsilon^*, k] \bigr)_{[\mathscr{H}]^2} 
        \nonumber
        \\[1ex]
        = &~ \bigl( \mathcal{P}_\varepsilon^* [M_\eta (\eta_\varepsilon^* -\eta_\mathrm{ad}), M_\theta (\theta_\varepsilon^* -\theta_\mathrm{ad})], [M_u (h-u_\varepsilon^*), M_v k] \bigr)_{[\mathscr{H}]^2} 
        \nonumber
        \\[0ex]
        & \hspace{4ex}+\bigl( [M_u u_\varepsilon^*, M_v v_\varepsilon^*], [h-u_\varepsilon^*, k] \bigr)_{[\mathscr{H}]^2} 
        \nonumber
        \\[1ex]
        = & \bigl(M_u(p_\varepsilon^* +u_\varepsilon^*), h-u_\varepsilon^* \bigr)_{\mathscr{H}} + \bigl(M_v(z_\varepsilon^* + v_\varepsilon^*), k \bigr)_{\mathscr{H}}, 
        \\[1ex]
        & \hspace{4ex} \mbox{for any $ [h, k] \in \mathscr{U}_\mathrm{ad}^K  $.}
        \nonumber
    \end{align}

    Now, in \eqref{nec01}, let us consider the case when $ [h, k] = [h, 0] \in \mathscr{U}_\mathrm{ad}^K $ with arbitrary $ h \in K $. Then, we have:
    \begin{align}\label{Thm.5-100}
        0 \leq  \bigl( M_u(p_\varepsilon^* +u_\varepsilon^*), h & -u_\varepsilon^* \bigr)_{\mathscr{H}} = -M_u \bigl( -p_\varepsilon^* -u_\varepsilon^*, h -u_\varepsilon^* \bigr)_{\mathscr{H}}
        \nonumber
        \\
        & \mbox{for any $ h \in K $.}
    \end{align}
    It is equivalent to \eqref{Thm.5-001}. Indeed, if $ M_u > 0 $, then the equivalence of \eqref{Thm.5-001} and \eqref{Thm.5-100} is a straightforward consequence of (\hyperlink{Fact2l}{Fact\,2}). Also, if $ M_u = 0 $, then the both \eqref{Thm.5-001} and \eqref{Thm.5-100} coincides with the tautology ``$0 = 0$''. 

    In the meantime, putting $ [h, k] = [u_\varepsilon^*, k] \in \mathscr{U}_\mathrm{ad}^K $ with arbitrary $ k \in \mathscr{H} $, one can see that:
    \begin{equation*}
        \bigl( M_v(v^* +z_\varepsilon^*), k \bigr)_\mathscr{H} \geq 0 \mbox{ for any $ k \in \mathscr{H} $.}
    \end{equation*}
    This implies the equality \eqref{Thm.5-002}.
    \medskip

    Thus, we conclude  Main Theorem \ref{mainTh03} (\hyperlink{III-A}{III-A}). \qed
}
\bigskip

Next, before the proof of Main Theorem \ref{mainTh03} (\hyperlink{III-B}{III-B}), we prepare the following lemma.
\begin{lem}\label{u.comp}
Let us assume (\hyperlink{A5l}{A5}) and (\hyperlink{A6l}{A6}), and fix a constraint $K = \jump{\kappa^0, \kappa^1} \in \mathfrak{K}$. Also, let us assume that: 
\begin{subequations}\label{u.comp01}
    \begin{align}\label{u.comp01a}
        & \tilde{p} \in \mathscr{H}, ~ \{ \tilde{p}_n\}_{n=1}^\infty \subset \mathscr{H}, \mbox{ and } \tilde{p}_n \to \tilde{p} \mbox{ in } \mathscr{H} \mbox{ as } n \to \infty,
    \end{align}
    and let us put:
    \begin{align}\label{u.comp01b}
        & \begin{cases}
            \tilde{u} := \mathrm{proj}_{K}(\tilde{p}) \mbox{ in $ \mathscr{H} $, }
            \\[1ex]
            \tilde{u}_n := \mathrm{proj}_{K_n}(\tilde{p}_n) \mbox{ in $ \mathscr{H} $}, \mbox{ for } n = 1, 2, 3, \ldots.
        \end{cases}
    \end{align}
\end{subequations}
Then, it holds that: 
\begin{align}
 \tilde{u}_n \to \tilde{u} \mbox{ in } \mathscr{H}, \mbox{ as } n \to \infty.
\end{align}
\end{lem}
\paragraph{Proof.}{
By using the assumptions as in \eqref{u.comp01}, Remark \ref{Rem.proj01}, and Lemma \ref{Lem.MoscoK}, this Lemma is easily verified as follows. 
\begin{align*}
 |\tilde{u}_n - \tilde{u}|_{\mathscr{H}} & \leq |\mathrm{proj}_{K_n}(\tilde{p}_n) - \mathrm{proj}_{K_n}(\tilde{p}) |_{\mathscr{H}} + | \mathrm{proj}_{K_n}(\tilde{p}) - \mathrm{proj}_{K}(\tilde{p}) |_{\mathscr{H}} \nonumber
 \\
    & \leq |\tilde{p}_n - \tilde{p} |_{\mathscr{H}} + | \mathrm{proj}_{K_n}(\tilde{p}) - \mathrm{proj}_{K}(\tilde{p}) |_{\mathscr{H}} \to 0,  \mbox{ as $ n \to \infty $}.
\end{align*}

\qed
}

Now, we are on the stage to prove Main Theorem \ref{mainTh03} (\hyperlink{III-B}{III-B}).

\paragraph{Proof of Main Theorem \ref{mainTh03} (III-B).}{
Let us notice that the assumptions \eqref{w.i01} and \eqref{IIIB00} guarantee that:
\begin{itemize}
    \item the sequence of initial pairs $ \{ [\eta_{0, n}, \theta_{0, n}] \}_{n = 1}^\infty $ is bounded in $|SK{D_0} = \bigl( V \cap L^\infty(\Omega) \bigr) \times V_0 $;
    \item the sequence  $ \{ u_n^* \in K_n \}_{n = 1}^\infty $, consisting of the first components of optimal controls $ [u_n^*, v_n^*] \in \mathscr{U}_\mathrm{ad}^{K_n} $ of (\hyperlink{OP}{OP})$_{\varepsilon_n}^{K_n}$, $ n = 1, 2, 3, \dots $, is bounded in $ L^\infty(Q) $.
\end{itemize}
Hence, with the compact embedding $|SK{D =} V \times V_0 \subset [H]^2 $ and Alaoglu's theorem in mind, we may suppose that:
\begin{equation}\label{IIIB10}
    \begin{cases}
        [\eta_{0, n_i}, \theta_{0, n_i}] \to [\eta_0, \theta_0] \mbox{ in $ [H]^2 $, and weakly in $ V \times V_0 $,}
        \\[0.5ex]
        \eta_{0, n_i} \to \eta_0 \mbox{ weakly-$*$ in $ L^\infty(\Omega) $,}
        \\[0.5ex]
        u_{n_i}^* \to u^{**} \mbox{ weakly-$*$ in $ L^\infty(Q) $, as $ i \to \infty $,}
    \end{cases}
\end{equation}
for the subsequence $ \{ n_i \}_{i = 1}^\infty \subset \{n\} $ and the limiting optimal control $ [u^{**}, v^{**}] \in \mathscr{U}_\mathrm{ad}^K $, as in Main Theorem \ref{mainTh02} (\hyperlink{II-B}{II-B}). 

By \eqref{w.i01} and \eqref{IIIB10}, we can apply Main Theorem \ref{mainTh01} (\hyperlink{I-B}{I-B}), to the solutions $ [\eta^{**}, \theta^{**}] \in [\mathscr{H}]^2 $ and $ [\eta_{n_i}^*, \theta_{n_i}^*] \in [\mathscr{H}]^2 $, $ i = 1, 2, 3, \dots $, as in \eqref{mTh02-12}, and can deduce that:
\begin{equation}\label{IIIB11}
    \eta_{n_i}^* \to \eta^{**} \mbox{ weakly-$*$ in $ L^\infty(Q) $, as $ i \to \infty $.}
\end{equation}

Meanwhile, by taking more subsequence(s) if necessary, one can see from \eqref{setRem4}, \eqref{mTh02-12}, and \eqref{IIIB11} that:
\begin{subequations}\label{IIIB12}
    \begin{align}\label{IIIB12a}    
        \lambda_i^* := ~& \mathcal{R}_T \bigl[ g'(\eta_{n_i}^*) \bigr] \to \lambda^{**} := \mathcal{R}_T \bigl[ g'(\eta^{**}) \bigr]\mbox{ weakly-$*$ in $ L^\infty(Q) $,}
        \nonumber
        \\
        & \mbox{and in the pointwise sense a.e. in $ Q $,}
    \end{align}
    \begin{align}\label{IIIB12b}
        \omega_i^* := ~& \mathcal{R}_T \bigl[ \alpha'(\eta_{n_i}^*) \nabla f_{\varepsilon_{n_i}}(\nabla \theta_{n_i}^*) \bigr] \to \omega^{**} := \mathcal{R}_T \bigl[ \alpha'(\eta^{**}) \nabla f_{\varepsilon}(\nabla \theta^{**}) \bigr] 
        \nonumber
        \\
        & \mbox{ weakly-$*$ in $ [L^\infty(Q)]^N $, and in the pointwise sense a.e. in $ Q $,}
    \end{align}
    \begin{align}\label{IIIB12c}
        A_i^* := ~& \mathcal{R}_T \bigl[ \alpha(\eta_{n_i}^*) \nabla^2 f_{\varepsilon_{n_i}}(\nabla \theta_{n_i}^*) \bigr] \to A^{**} := \mathcal{R}_T \bigl[ \alpha(\eta^{**}) \nabla^2 f_{\varepsilon}(\nabla \theta^{**}) \bigr] 
        \nonumber
        \\
        & \mbox{ weakly-$*$ in $ [L^\infty(Q)]^{N \times N} $, and in the pointwise sense a.e. in $ Q $,}
    \end{align}
    \begin{align}\label{IIIB12d}
        \mu_i^* := ~& \mathcal{R}_T \bigl[ \alpha''(\eta_{n_i}^*) f_{\varepsilon_{n_i}}(\nabla \theta_{n_i}^*) \bigr] \to \mu^{**} := \mathcal{R}_T \bigl[ \alpha''(\eta^{**}) f_{\varepsilon}(\nabla \theta^{**}) \bigr] 
        \nonumber
        \\
        & \mbox{ weakly-$*$ in $ L^\infty(0, T; H) $, 
        }
    \end{align}
    and
    \begin{align}\label{IIIB12e}
        \mu_i^*(t) \to \mu^{**}(t) \mbox{ in $ H $, for a.e. $ t \in (0, T) $, as $ i \to \infty $.}
    \end{align}
\end{subequations}

Now, let us denote by $ \mathcal{P}_0^{**} $ and $ \mathcal{P}_i^{**} $, $ i = 1, 2, 3, \dots $, the operators $ \mathcal{P}_\varepsilon^* $, as in Remark \ref{Rem.mTh03}, in the cases when:
\begin{align*}
    [\mu, \lambda, &\, \omega, A] = [\mu^{**}, \lambda^{**}, \omega^{**}, A^{**}]
    \\
    & \mbox{in $ L^\infty(0, T; H) \times L^\infty(Q) \times [L^\infty(Q)]^N \times [L^\infty(Q)]^{N \times N} $,}
\end{align*}
\begin{align*}
    \varepsilon & = \varepsilon_i, \mbox{ and }[\mu, \lambda, \omega, A] = [\mu_i^{*}, \lambda_i^{*}, \omega_i^{*}, A_i^{*}]
    \\
    & \mbox{in $ L^\infty(0, T; H) \times L^\infty(Q) \times [L^\infty(Q)]^N \times [L^\infty(Q)]^{N \times N} $, $ i = 1, 2, 3, \dots $,}
\end{align*}
respectively. Then, as a consequence of Proposition \ref{ASY_Cor.2}, Main Theorem \ref{mainTh03} (\hyperlink{III-A}{III-A}), and Remark \ref{Rem.mTh03}, we can derive from \eqref{mTh02-12} and \eqref{IIIB12} that:
\begin{align}\label{IIIB13}
    [p_i^*, & z_i^*] := \mathcal{P}_i^{**}\bigl[ M_\eta(\eta_{n_i}^{*} -\eta_\mathrm{ad}), M_\theta(\theta_{n_i}^* -\theta_\mathrm{ad}) \bigr] 
    \nonumber
    \\
    & \to [p^{**}, z^{**}] := \mathcal{P}_0^{**}\bigl[ M_\eta(\eta^{**} -\eta_\mathrm{ad}), M_\theta(\theta^{**} -\theta_\mathrm{ad}) \bigr]\mbox{ in $ [C([0, T]; H)]^2 $, } 
    \nonumber
    \\
    & \mbox{in $ \mathscr{Y} $, and weakly in $ W^{1, 2}(0, T; V^*) \times W^{1, 2}(0, T; V_0^*) $, as $ i \to \infty $.}
\end{align}
Furthermore, taking into account \eqref{IIIB13} and Lemma \ref{u.comp}, one can infer that:
\begin{subequations}\label{IIIB14}
\begin{align}\label{IIIB14a}
    M_u u_{n_i}^* ~& = M_u \mathrm{proj}_{K_{n_i}}(-p_i^*) \to M_u u^{**} = M_u \mathrm{proj}_{K}(-p^{**}) \mbox{ in $ \mathscr{H} $, as $ i \to \infty $,} 
\end{align}
and
\begin{align}\label{IIIB14b}
    M_v v_{n_i}^* ~& = -M_v z_i^* \to M_v v^{**} = -M_v z^{**} \mbox{ in $ C([0, T]; H) $, in $\mathscr{Y}$,}
    \nonumber
    \\
    & \mbox{and weakly in $ W^{1, 2}(0, T; V_0^*) $, as $ i \to \infty $.} 
\end{align}
\end{subequations}

\eqref{IIIB10} and \eqref{IIIB14} are sufficient to verify the convergences as in \eqref{IIIB01}, and to conclude Main Theorem \ref{mainTh03} (\hyperlink{III-B}{III-B}). \qed
}

\section{Proof of Main Theorem \ref{mainTh04}}

Under the assumptions (\hyperlink{A1l}{A1})--(\hyperlink{A5l}{A5}) and the situation $ \neg $(\hyperlink{rs0}{r.s.0}), let us set:
\begin{align*}
    & \begin{cases}
        \varepsilon_n := \varepsilon +\frac{1}{n},
        \\[1ex]
        \eta_{0, n} := (-n) \vee (n \wedge \eta_0) \mbox{ a.e. in $ \Omega $,}
        \\[1ex]
        \theta_{0, n} := \theta_0 \mbox{ a.e. in $ \Omega $,}
        \\[1ex]
        \kappa_n^\ell := (-n) \vee (n \wedge \kappa^\ell) \mbox{ a.e. in $ Q $, $ \ell = 0, 1 $,}
    \end{cases}
    n = 1, 2, 3, \dots.
\end{align*}
Then, we immediately see that:
\begin{description}
    \item[\hypertarget{star4}{($\star\,4$)}]$ \{ \varepsilon_n \}_{n = 1}^\infty \subset (\varepsilon, \infty) $, $ \{ [\eta_{0, n}, \theta_{0, n}] \}_{n = 1}^\infty \subset D_0 $, and $ \{ K_n \}_{n = 1}^\infty := \{ \jump{\kappa_n^0, \kappa_n^1} \} \subset \mathfrak{K}_0 $, and these sequences fulfill the assumptions \eqref{w.i01} and (\hyperlink{A6l}{A6}), as in Main Theorems \ref{mainTh01}--\ref{mainTh03}. 
\end{description}
Additionally,  we can apply Main Theorem \ref{mainTh01} (\hyperlink{I-A}{I-A}) and Main Theorem \ref{mainTh02} (\hyperlink{II-A}{II-A}), and can take sequences of functional pairs $ \{ [u_n^\circ, v_n^\circ] \}_{n = 1}^\infty $ and $ \{ [\eta_n^\circ, \theta_n^\circ] \}_{n = 1}^\infty $, such that:
\begin{itemize}
    \item for any $ n \in \N $, $ [u_n^\circ, v_n^\circ] \in \mathscr{U}_\mathrm{ad}^{K_n} $ is an optimal control of (\hyperlink{OP}{OP})$_{\varepsilon_n}^{K_n}$;
    \item for any $ n \in \N $, $ [\eta_n^\circ, \theta_n^\circ] \in [\sH]^2 $ is the solution to (S)$_{\varepsilon_n}$, for the initial pair $ [\eta_{0, n}, \theta_{0, n}] $ and forcing pair $ [u_n^\circ, v_n^\circ] $.
\end{itemize}
Also, applying Main Theorem \ref{mainTh01} (\hyperlink{I-B}{I-B}) and Main Theorem \ref{mainTh02} (\hyperlink{II-B}{II-B}), we can find subsequences of $ \{[u_n^\circ, v_n^\circ]\}_{n = 1}^\infty $ and $ \{[\eta_n^\circ, \theta_n^\circ]\}_{n = 1}^\infty $ (not relabeled), together with limiting pairs $ [u^\circ, v^\circ] \in [\sH]^2 $ and $ [\eta^\circ, \theta^\circ] \in [\sH]^2 $, and a limiting function $ \sigma^\circ \in [L^\infty(Q)]^N $, such that:
\begin{subequations}\label{sec7-01}
\begin{align}\label{sec7-01a}
    & [M_u u_n^\circ, M_v v_n^\circ] \to [M_u u^\circ, M_v v^\circ] \mbox{ weakly in $ [\sH]^2 $,}
\end{align}
\begin{align}\label{sec7-01b}
    [\eta_n^\circ, \theta_n^\circ] &~  \to [\eta^\circ, \theta^\circ] \mbox{ in $ [C([0, T]; H)]^2 $, in $ \mathscr{Y} $,}
    \nonumber
    \\
     &~ \mbox{and weakly-$*$ in $ L^\infty(0, T; V) \times L^\infty(0, T; V_0) $,}
\end{align}
\begin{align}\label{sec7-01c}
    [\nabla \eta_n^\circ,&~  \nabla \theta_n^\circ] \to [\nabla \eta^\circ, \nabla \theta^\circ] \mbox{ in $ [L^2(0, T; [H]^N)]^2 $,} 
    \nonumber
    \\
    &~ \mbox{and in the pointwise sense a.e. in $ Q $,}
\end{align}
\begin{align}\label{sec7-01d}
    & 
    \begin{cases}
        \mu_n^\circ := \alpha''(\eta_n^\circ) f_{\varepsilon_n}(\nabla \theta_n^\circ) \to \mu^\circ := \alpha''(\eta^\circ) f_\varepsilon(\nabla \theta^\circ) 
        \\
        \quad~ \mbox{weakly-$ * $ in $ L^\infty(0, T; H) $,}
        \\
        \quad~ \mbox{and in the pointwise sense a.e. in $ Q $,}
        \\[1ex]
        \mu_n^\circ(t) \to \mu^\circ(t) \mbox{ in $ H $, for a.e. $ t \in (0, T) $,}
    \end{cases}
\end{align}
\begin{align}\label{sec7-01e}
    \lambda_n^\circ := ~&  g'(\eta_n^\circ) \to \lambda^\circ := g'(\eta^\circ) \mbox{ in $ \mathscr{H} $, weakly-$ * $ in $ L^\infty(Q) $,}
    \nonumber
    \\
    & \mbox{and in the pointwise sense a.e. in $ Q $,}
\end{align}
\begin{align}\label{sec7-01f}
    \nabla f_{\varepsilon_n}(\nabla \theta_n^\circ) \to \sigma^\circ \mbox{ weakly-$*$ in $ [L^\infty(Q)]^N $,}
\end{align}
and
\begin{align}\label{sec7-01g}
    \omega_n^\circ := \alpha'(\eta_n^\circ) &~ \nabla f_{\varepsilon_n}(\nabla \theta_n^\circ) \to \alpha'(\eta^\circ)\sigma^\circ 
    \nonumber
    \\
    & \mbox{ weakly-$*$ in $ [L^\infty(Q)]^N $, as $ n \to \infty $.} 
\end{align}
\end{subequations}
Additionally, from \eqref{sec7-01c}, \eqref{sec7-01f}, Remark \ref{Rem.MG}, and \cite[Proposition 2.16]{MR0348562}, one can observe that:
\begin{equation}\label{conv3B-07}
    \sigma^\circ \in \partial f_\varepsilon(\nabla \theta^\circ) = 
    \begin{cases}
        \{ \nabla f_\varepsilon(\nabla \theta^\circ) \}, \mbox{ if $ \varepsilon > 0 $,}
        \\[1ex]
        \mathrm{Sgn}^N(\nabla \theta^\circ), \mbox{ if $ \varepsilon = 0 $,}
    \end{cases} \mbox{a.e. in $ Q $.}
\end{equation}

Next, for any $ n \in \N $, let us put:
\begin{equation}\label{A_n^circ}
    A_n^\circ := \alpha(\eta_n^\circ) \nabla^2 f_{\varepsilon_n}(\nabla \theta_n^\circ) \mbox{ in $ [L^\infty(Q)]^{N\times N} $,}
\end{equation}
and let us denote by $ \mathcal{P}_n^\circ  $ the operator $ \mathcal{P}_\varepsilon^* \in \mathscr{L}([\sH]^2; \mathscr{Z}) $, as in Remark \ref{Rem.mTh03}, in the case when 
\parbox{15.75cm}{
    the constant $ \varepsilon > 0 $ (in Remark \ref{Rem.mTh03}) and the sextuplet $ [a, b, \mu, \lambda, \omega, A] \in \mathscr{S} $ is replaced by $ \varepsilon_n > 0 $ and $ \mathcal{R}_T [\alpha_0, -\partial_t \alpha_0, \mu_n^\circ, \lambda_n^\circ, \omega_n^\circ, A_n^\circ] \in \mathscr{S} $, respectively.}

On this basis, let us set:
\begin{equation*}
    [p_n^\circ, z_n^\circ] := \mathcal{P}_n^\circ \bigl[ M_\eta(\eta_n^\circ -\eta_\mathrm{ad}), M_\theta(\theta_n^\circ -\theta_\mathrm{ad}) \bigr] \mbox{ in $ \mathscr{Z} $, for $ n = 1, 2, 3, \dots $.}
\end{equation*}
Then, from Main Theorem \ref{mainTh03} (\hyperlink{III-A}{III-A}), it is inferred that:
\begin{subequations}\label{sec7-02}
\begin{equation}\label{sec7-02a}
    (M_u (p_n^\circ + u_n^\circ), h-u_n^\circ)_\mathscr{H} \geq 0,\ \mbox{ for any $ h \in K_n = \jump{\kappa_n^0, \kappa_n^1} $,} 
\end{equation}
\begin{equation}\label{sec7-02b}
    M_v (z_n^\circ +  v_n^\circ) = 0\ \mbox{ in $\mathscr{H} $,}
\end{equation}
\begin{align}\label{sec7-02c}
    \bigl< -\partial_t & p_n^\circ, \varphi \bigr>_{\mathscr{V}} +\bigl( \nabla p_n^\circ, \nabla \varphi \bigr)_{[\sH]^N} +\bigl< \mu_n^\circ p_n^\circ, \varphi \bigr>_{\mathscr{V}} 
    \nonumber
    \\
    & +\bigl( \lambda_n^\circ p_n^\circ +\omega_n^\circ  \cdot \nabla z_n^\circ, \varphi \bigr)_{\mathscr{H}} = \bigl( M_\eta  (\eta_n^\circ -\eta_\mathrm{ad}), \varphi \bigr)_{\mathscr{H}}, 
    \mbox{for any $ \varphi \in \mathscr{V} $,}
\end{align}
\begin{align}\label{sec7-02d}
    \bigl< -\alpha_0 \partial_t z_n^\circ, \psi \bigr>_{\mathscr{V}_0} +\bigl( & ( -\partial_t \alpha_0) z_n^\circ, \psi \bigr)_{\mathscr{H}} +\bigl( A_n^\circ \nabla z_n^\circ   +\nu^2 \nabla z_n^\circ +p_n^\circ \omega_n^\circ , \nabla \psi \bigr)_{[\sH]^N}
    \nonumber
    \\
    & = \bigl( M_\theta (\theta_n^\circ -\theta_\mathrm{ad}), \psi \bigr)_{\mathscr{H}}, \mbox{ for any $ \psi \in \mathscr{V}_0 $,}
\end{align}
and
\begin{equation}\label{sec7-02e}
    [p_n^\circ(T), z_n^\circ(T)] = [0, 0] \mbox{ in $ [H]^2 $, $ n = 1, 2, 3, \dots $.}
\end{equation}
\end{subequations}
Also, having in mind \eqref{sec7-01}--\eqref{A_n^circ}, and applying Proposition \ref{Prop(I-B)} to the case when:
\begin{equation*}
    \begin{array}{c}
    \begin{cases}
        [a^1, b^1, \mu^1, \lambda^1, \omega^1, A^1] = [a^2, b^2, \mu^2, \lambda^2, \omega^2, A^2] 
        \\
        \hspace{7ex}= \mathcal{R}_T[\alpha_0, -\partial_t \alpha_0, \mu_n^\circ, \lambda_n^\circ, \omega_n^\circ, A_n^\circ],
        \\[1ex]
        [p_0^1, z_0^1] = [p_0^2, z_0^2] = [0, 0],
        \\[1ex]
        [h^1, k^1] = \mathcal{R}_T[M_\eta(\eta_n^\circ -\eta_\mathrm{ad}), M_\theta(\theta_n^\circ - \theta_\mathrm{ad})], 
        \\[1ex]
        [h^2, k^2] = [0, 0],
        \\[1ex]
        [p^1, z^1] = \mathcal{R}_{T} [p_n^\circ, z_n^\circ], ~ [p^2, z^2] = [0, 0],
    \end{cases}
    \mbox{for $ n = 1, 2, 3, \dots $.}
    \end{array}
\end{equation*}
we deduce that:
\begin{align}\label{sec7-03}
    \frac{d}{dt} & \bigl( \bigl| (\mathcal{R}_T p_n^\circ)(t) \bigr|_H^2 +\bigl| \mathcal{R}_T \bigl( \sqrt{\alpha_0} z_n^\circ \bigr)(t) \bigr|_H^2 \bigr)
    \nonumber
    \\
    & \qquad +\bigl( \bigl|(\mathcal{R}_T p_n^\circ)(t)\bigr|_V^2 +\nu^2 \bigl|(\mathcal{R}_T z_n^\circ)(t)\bigr|_{V_0}^2 \bigr)
    \nonumber
    \\
    & \leq 3 \bar{C}_0^* \bigl( \bigl| (\mathcal{R}_T p_n^\circ)(t) \bigr|_H^2 +\bigl| \mathcal{R}_T \bigl( \sqrt{\alpha_0} z_n^\circ \bigr)(t) \bigr|_H^2 \bigr) 
    \nonumber
    \\
    & \qquad +2\bar{C}_0^* \big( \bigl| \mathcal{R}_T \bigl( M_\eta (\eta_n^\circ -\eta_\mathrm{ad}) \bigr)(t) \bigr|_{V^*}^2 +\bigl|\mathcal{R}_T \bigl( M_\theta (\theta_n^\circ -\theta_\mathrm{ad}) \bigr) (t) \bigr|_{V_0^*}^2 \bigr), 
    \\
    & \mbox{for a.e. $ t \in (0, T) $, $ n = 1, 2, 3, \dots $,}
    \nonumber
\end{align}
with use of the constant $ \bar{C}_0^* $ as in \eqref{est3-01a}. As a consequence of \eqref{sec7-01b}, \eqref{sec7-03}, (\hyperlink{A4l}{A4}), and Gronwall's lemma, it is observed that:
\begin{description}
    \item[\textmd{(\hypertarget{*5}{$\star$\,5})}]the sequence $ \{ [p_n^\circ, z_n^\circ] \}_{n = 1}^\infty $ is bounded in $ [C([0, T]; H)]^2 \cap \mathscr{Y} $.
\end{description}

In the meantime, from \eqref{f_eps}, \eqref{sec7-01b}--\eqref{sec7-01g}, \eqref{sec7-02c}, \eqref{sec7-02d}, (\hyperlink{A3l}{A3}),  and  Remark \ref{Rem.Prelim01}, we can derive the following estimates:
\begin{subequations}\label{sec7-04}
\begin{align}\label{sec7-04a}
    \bigl| \bigl< \partial_t p_n^\circ, \varphi \bigr>_{\mathscr{V}} \bigr| \leq &~ \bigl| \bigl< \mu_n^\circ p_n^\circ, \varphi \bigr>_{\mathscr{V}} \bigr| +\bigl| \bigl( \nabla p_n^\circ, \nabla \varphi \bigr)_{[\sH]^N} \bigr| +\bigl| \bigl( \lambda_n^\circ p_n^\circ +\omega_n^\circ \cdot \nabla z_n^\circ, \varphi \bigr)_{\mathscr{H}} \bigr|
    \nonumber
    \\
    &~  +\bigl| \bigl( M_\eta  (\eta_n^\circ -\eta_\mathrm{ad}), \varphi \bigr)_{\mathscr{H}} \bigr| \leq C_1^\circ |\varphi|_{\mathscr{V}}, \mbox{ for any $ \varphi \in \mathscr{V} $,}
\end{align}
and 
\begin{align}\label{sec7-04b}
    \bigl| \bigl< -\mbox{div} & (A_n^\circ \nabla z_n^\circ), \psi \bigr>_{\mathscr{W}_{0}} \bigr| 
    = \bigl| \bigl( A_n^\circ \nabla z_n^\circ, \nabla \psi \bigr)_{[\sH]^N} \bigr| \leq \bigl| \bigl( \alpha_0 z_n^\circ, \partial_t \psi \bigr)_{\mathscr{H}} \bigr|  
    \nonumber
    \\
    ~& +\bigl| \bigl( \nu^2 \nabla z_n^\circ +p_n^\circ \omega_n^\circ , \nabla \psi \bigr)_{[\sH]^N} \bigr| +\bigl| \bigl( M_\theta (\theta_n^\circ -\theta_\mathrm{ad}), \psi \bigr)_{\mathscr{H}} \bigr|
    \\
    & \leq C_2^\circ |\psi|_{\mathscr{W}_0}, \mbox{ for any $ \psi \in C_\mathrm{c}^\infty(Q) $, $ n = 1, 2, 3, \dots $,}
    \nonumber
\end{align}
\end{subequations}
with $ n $-independent positive constants:
\begin{align*}
    C_1^\circ &:= \sup_{n \in \N} \left\{ \begin{array}{c}
            (1 +(C_V^{L^4})^2|\mu_n^\circ|_{L^\infty(0, T; H)} + |\lambda_n^\circ|_{L^\infty(Q)} + |\omega_n^\circ|_{[L^\infty(Q)]^N})  
        \\[1ex]
        \cdot \bigl( \bigl| [p_n^\circ, z_n^\circ] \bigr|_{\mathscr{Y}} +\bigl| M_\eta (\eta_n^\circ -\eta_\mathrm{ad}) \bigr|_{\mathscr{H}} \bigr) 
    \end{array} \right\} ~(< \infty),
\end{align*}
and
\begin{align*}
    C_2^\circ &:= \sup_{n \in \N} \left\{ \begin{array}{c}
        (1 +\nu^2 +C_{V_0}^H|\alpha_0|_{L^\infty(Q)} +|\omega_n^\circ|_{[L^\infty(Q)]^N})  
        \\[1ex]
        \cdot \bigl( \bigl| [p_n^\circ, z_n^\circ] \bigr|_{\mathscr{Y}} +\bigl| M_\theta (\theta_n^\circ -\theta_\mathrm{ad}) \bigr|_{\mathscr{H}} \bigr) 
    \end{array} \right\} ~(< \infty),
\end{align*}
where $ C_{V}^{L^4} > 0 $ and $ C_{V_0}^H > 0 $ are the constants of embeddings $ V \subset L^4(\Omega) $ and $ V_0 \subset H $, respectively. 
\medskip

Due to \eqref{sec7-01d}--\eqref{sec7-01g}, \eqref{sec7-04}, (\hyperlink{*5}{$\star$\,5}), Lemma \ref{Lem.ax01}, and the compactness theory of Aubin's type (cf. \cite[Corollary 4]{MR0916688}), we can find subsequences of $ \{ [p_n^\circ, z_n^\circ] \}_{n = 1}^\infty \subset \mathscr{Y} $, $ \{ \omega_n^\circ \cdot \nabla z_n^\circ \}_{n = 1}^\infty \subset \mathscr{H} $, and $ \{ -\mbox{div} (A_n^\circ \nabla z_n^\circ) \}_{n = 1}^\infty \subset \mathscr{W}_0^* $ (not relabeled), together with the respective limits $ [p^\circ, z^\circ] \in \mathscr{Y} $, $ \xi^\circ \in \mathscr{H} $, and $ \zeta^\circ \in \mathscr{W}_0^* $, such that:
\begin{subequations}\label{sec7-05}
\begin{align}\label{sec7-05a}
    & \begin{cases}
        [p_n^\circ, z_n^\circ] \to [p^\circ, z^\circ] \mbox{ weakly in $ \mathscr{Y} $,}
        \\[1ex]
        p_n^\circ \to p^\circ \mbox{ in $ \mathscr{H} $, weakly in $ W^{1, 2}(0, T; V^*) $,}
        \\
        \quad \mbox{and in the pointwise sense a.e. in $ Q $,}
    \end{cases}
\end{align}
\begin{align}\label{conv3B-14}
    \mu_n^\circ p_n^\circ \to \mu^\circ p^\circ \mbox{ weakly in $ \mathscr{V}^* $,}
\end{align}
\begin{align}\label{sec7-05b}
    \lambda_n^\circ p_n^\circ \to \lambda^\circ p^\circ \mbox{ in $ \mathscr{H} $,}
\end{align}
\begin{align}\label{sec7-05c}
    p_n^\circ \omega_n^\circ \to p^\circ  \alpha'(\eta^\circ) \sigma^\circ \mbox{ weakly in $ \mathscr{H} $,}
\end{align}
\begin{align}\label{sec7-05d}
    \begin{cases}
        \nabla f_{\varepsilon_n}(\nabla \theta_n^\circ) \cdot \nabla z_n^\circ \to \xi^\circ ~\mbox{weakly in $ \mathscr{H} $,}
        \\[1ex]
        \omega_n^\circ \cdot \nabla z_n^\circ = \alpha'(\eta_n^\circ)\nabla f_{\varepsilon_n}(\nabla \theta_n^\circ) \cdot \nabla z_n^\circ \to \alpha'(\eta^\circ) \xi^\circ 
        \\
        \hspace{22.5ex}\mbox{weakly in $ \mathscr{H} $,}
    \end{cases}
\end{align}
and
\begin{align}\label{conv3B-13}
    -\mbox{div} (A_n^\circ \nabla z_n^\circ) \to \zeta^\circ \mbox{ weakly in $ \mathscr{W}_0^* $, as $ n \to \infty $.}
\end{align}
\end{subequations}

Now, the properties \eqref{Thm.5-10}--\eqref{Thm.5-12} will be verified through the limiting observations for \eqref{sec7-02}, as $ n \to \infty $, with use of \eqref{sec7-01} and \eqref{sec7-05}. 

Finally, we verify the properties as in \eqref{Thm.5-141}, under the situation (\hyperlink{rs1}{r.s.1}). To this end, we first invoke \eqref{exM01} and \eqref{sec7-01c}, and confirm that:
\begin{equation}\label{sec7-07-01}
    \begin{cases}
        \bigl| \varphi\bigl( \nabla f_\varepsilon(\nabla \theta_n^\circ) -\nabla f_\varepsilon(\nabla \theta^\circ) \bigr) \bigr| \to 0, \mbox{ in the pointwise sense,} 
        \\
        \qquad \mbox{a.e. in $ Q $, as $ n \to \infty $,}
        \\[1ex]
        \bigl| \varphi\bigl( \nabla f_\varepsilon(\nabla \theta_n^\circ) -\nabla f_\varepsilon(\nabla \theta^\circ) \bigr) \bigr| \leq 2 |\varphi|, \mbox{ a.e. in $ Q $, $ n = 1, 2, 3, \dots, $}
        \\[1ex]
        \qquad \mbox{for any $ \varphi \in \mathscr{H} $.}
    \end{cases}
\end{equation}
With \eqref{sec7-07-01}, (\hyperlink{rs1}{r.s.1}), (\hyperlink{star4}{$\star\,4$}), and (\hyperlink{A4l}{A4}) in mind, using \eqref{exM01} and \eqref{sec7-01b}, and applying the dominated convergence theorem \cite[Theorem 10 on page 36]{MR0492147} yield that:
\begin{align}\label{sec7-08}
    \bigl| \varphi \bigl( & \alpha'(\eta_n^\circ) \nabla f_{\varepsilon_n}(\nabla \theta_n^\circ) -\alpha'(\eta^\circ) \nabla f_\varepsilon(\nabla \theta^\circ) \bigr) \bigr|_{[\mathscr{H}]^N}
    \nonumber
    \\
    & \leq \bigl| \varphi \bigl( \alpha'(\eta_n^\circ) -\alpha'(\eta^\circ) \bigr) \bigr|_{\mathscr{H}} 
    \nonumber
    \\
    & \qquad +|\alpha'|_{L^\infty(\R)} \bigl| \varphi \bigl( \nabla f_{\varepsilon}(\nabla \theta_n^\circ) -\nabla f_{\varepsilon}(\nabla \theta^\circ) \bigr) \bigr|_{[\mathscr{H}]^N}
    \nonumber
    \\
    & \qquad +|\alpha'|_{L^\infty(\R)} \bigl| \varphi \bigl( \nabla f_{\varepsilon_n}(\nabla \theta_n^\circ) -\nabla f_{\varepsilon}(\nabla \theta_n^\circ) \bigr) \bigr|_{[\mathscr{H}]^N}
    \nonumber
    \\
    & \leq \bigl| \varphi \bigl( \alpha'(\eta_n^\circ) -\alpha'(\eta^\circ) \bigr) \bigr|_{\mathscr{H}} 
    \nonumber
    \\
    & \qquad +|\alpha'|_{L^\infty(\R)} \bigl| \varphi \bigl( \nabla f_{\varepsilon}(\nabla \theta_n^\circ) -\nabla f_{\varepsilon}(\nabla \theta^\circ) \bigr) \bigr|_{[\mathscr{H}]^N}
   \nonumber
    \\
    & \qquad +\frac{2|\alpha'|_{L^\infty(\R)}}{\varepsilon} |\varepsilon_n -\varepsilon| |\varphi|_{\mathscr{H}}
    \\
    & \to 0, \mbox{ as $ n \to \infty $, for any $ \varphi \in \mathscr{H} $. }
   \nonumber
\end{align}
Owing to \eqref{sec7-05a} and \eqref{sec7-08}, one can further observe that:
\begin{align}\label{sec7-09}
    & \bigl( \alpha'( \eta_n^\circ) \nabla f_{\varepsilon_n}(\nabla \theta_n^\circ) \cdot \nabla z_n^\circ, \varphi \bigr)_{\mathscr{H}} = \bigl( \nabla z_n^\circ, \varphi \alpha'(\eta_n^\circ) \nabla f_{\varepsilon_n}(\nabla \theta_n^\circ) \bigr)_{[\mathscr{H}]^N}
    \nonumber
    \\
    & \qquad \to   \bigl( \alpha'(\eta^\circ) \nabla f_{\varepsilon}(\nabla \theta^\circ) \cdot \nabla z^\circ, \varphi \bigr)_{\mathscr{H}} = \bigl( \nabla z^\circ, \varphi \alpha'(\eta^\circ) \nabla f_{\varepsilon}(\nabla \theta^\circ) \bigr)_{[\mathscr{H}]^N}
    \\
    & \qquad \qquad \mbox{as $ n \to \infty $, for any $ \varphi \in \mathscr{H} $. }
    \nonumber
\end{align}
Meanwhile, from \eqref{exM02}, \eqref{sec7-01b}, \eqref{sec7-01c}, (\hyperlink{star4}{$\star\,4$}), and (\hyperlink{A4l}{A4}), it is inferred that:
\begin{align*}
    \bigl| \bigl( \alpha & (\eta_n^\circ)\nabla^2 f_{\varepsilon_n}(\nabla \theta_n^\circ) -\alpha(\eta^\circ) \nabla^2 f_{\varepsilon}(\nabla \theta^\circ) \bigr) \nabla \psi \bigr|_{[\mathscr{H}]^N}
    \\
    & \leq |\alpha(\eta_n^\circ) -\alpha(\eta^\circ)|_\mathscr{H} \bigl| \nabla^2 f_{\varepsilon_n}(\nabla \theta_n^\circ) \bigr|_{L^\infty(Q; \R^{N \times N})} |\nabla \psi|_{C(\overline{Q}; \R^N)}
    \\
    & \qquad +|\alpha(\eta^\circ)|_{\mathscr{H}} \bigl| \nabla^2 f_{\varepsilon_n}(\nabla \theta_n^\circ) -\nabla^2 f_\varepsilon(\nabla \theta^\circ) \bigr|_{[\mathscr{H}]^{N \times N}} |\nabla \psi|_{C(\overline{Q}; \R^N)}
    \\
    & \leq \frac{N +1}{\varepsilon} |\alpha'|_{L^\infty(\R)} |\nabla \psi|_{C(\overline{Q}; \R^N)} |\eta_n^\circ -\eta|_{\mathscr{H}}
    \\
    & \qquad +\frac{3(N +1)^2}{\varepsilon^2} |\alpha(\eta^\circ)|_{\mathscr{H}} |\nabla \psi|_{C(\overline{Q}; \R^N)} \left(|\nabla(\theta_n^\circ -\theta)|_{[\mathscr{H}]^N} + \frac{1}{n}\right)
    \\
    & \to 0, \mbox{ as $ n \to \infty $, for any $ \psi \in C_\mathrm{c}^\infty(Q) $,}
\end{align*}
and therefore,
\begin{align}\label{sec7-11}
    & \bigl< -\mathrm{div} \bigl( \alpha(\eta_n^\circ) \nabla^2 f_{\varepsilon_n}(\nabla \theta_n^\circ) \nabla z_n^\circ \bigr), \psi \bigr> = \bigl( \nabla z_n^\circ, \alpha(\eta_n^\circ) \nabla^2 f_{\varepsilon_n}(\nabla \theta_n^\circ) \nabla \psi \bigr)_{[\mathscr{H}]^N}
    \nonumber
    \\
    & \quad  \to \bigl< -\mathrm{div} \bigl( \alpha(\eta^\circ) \nabla^2 f_{\varepsilon}(\nabla \theta^\circ) \nabla z^\circ \bigr), \psi \bigr> = \bigl( \nabla z^\circ, \alpha(\eta^\circ) \nabla^2 f_{\varepsilon}(\nabla \theta^\circ) \nabla \psi \bigr)_{[\mathscr{H}]^N},
    \\
    & \qquad\qquad  \mbox{as $ n \to \infty $, for any $ \psi \in C_\mathrm{c}^\infty(Q) $.}
    \nonumber
\end{align}

The fine properties as in \eqref{Thm.5-141} will be a consequence of \eqref{sec7-01f}, \eqref{sec7-01g}, \eqref{conv3B-07}, \eqref{sec7-05c}--\eqref{conv3B-13}, \eqref{sec7-09}, and \eqref{sec7-11}. 
\medskip

Thus, we complete the proof of Main Theorem \ref{mainTh04}. \qed


\end{document}